%% file: single-site-dilute-defects.tex
\documentclass[11pt,letterpaper]{amsart}
\usepackage[svgnames]{xcolor}
\usepackage[colorlinks,allcolors={DarkBlue}]{hyperref}

\usepackage{comment}

\usepackage{esint}
\usepackage{amssymb}
\usepackage{tikz}
\usetikzlibrary{calc}
\usepackage{graphicx}
\usepackage{amsrefs,enumitem}

\usepackage{manfnt}
\usepackage{graphicx} 

\DeclareRobustCommand{\cube}{%
  \mathchoice
    {\mbox{\hspace{.1em}\mancube}}
    {\mbox{\hspace{.1em}\mancube}}
    {\mbox{\scalebox{0.75}{\hspace{.1em}\mancube}}}
    {\mbox{\scalebox{0.6}{\hspace{.1em}\mancube}}}
}

\newtheorem*{theorem*}
{Theorem}
\newtheorem{theorem}{Theorem}
\newtheorem{lemma}[theorem]{Lemma}
\newtheorem{corollary}[theorem]{Corollary}
\newtheorem{proposition}[theorem]{Proposition}
\theoremstyle{definition}

\newtheorem{remark}[theorem]{Remark}

\newtheorem{definition}[theorem]{Definition}

\newcommand{\eref}[1]{(\ref{e.#1})}
\newcommand{\tref}[1]{Theorem \ref{t.#1}}
\newcommand{\lref}[1]{Lemma \ref{l.#1}}
\newcommand{\pref}[1]{Proposition \ref{p.#1}}
\newcommand{\cref}[1]{Corollary \ref{c.#1}}
\newcommand{\fref}[1]{Figure \ref{f.#1}}
\newcommand{\sref}[1]{Section \ref{s.#1}}

\newcommand{\partref}[1]{\ref{part.#1}}
\newcommand{\dref}[1]{Definition \ref{d.#1}}
\newcommand{\rref}[1]{Remark \ref{r.#1}}

\numberwithin{theorem}{section}
\numberwithin{equation}{section}

\newcommand{\Z}{\mathbb{Z}}
\newcommand{\R}{\mathbb{R}}

\newcommand{\grad}{\nabla}

\def\XXint#1#2#3{{\setbox0=\hbox{$#1{#2#3}{\int}$ }
\vcenter{\hbox{$#2#3$ }}\kern-.6\wd0}}

\newcommand{\ep}{\varepsilon}
\newcommand{\e}{\varepsilon}
\newcommand{\osc}{\mathop{\textup{osc}}}

\definecolor{darkgreen}{rgb}{0,0.4,0}
\def\strikethrough#1{\setbox0\hbox{#1}\rlap{#1}\hbox to \wd0{\hss\strikebox\hss}}
\def\strikebox{\vrule height 0.6\ht0 depth -0.4\ht0 width 1.1\wd0}

\begin{document}

\title{The pinning effect of dilute defects}
\author{William M Feldman}
\address{Department of Mathematics, The University of Utah, Salt Lake City}
\author{Inwon C Kim}
\address{Department of Mathematics, University of California, Los Angeles}
\begin{abstract}
We consider the Bernoulli free boundary problem with ``periodic defects", or a periodic array of compactly supported inhomogeneities. First we study the problem with a single defect site in $\R^d$. In contrast to the standard perforated domain theory, for our problem a given defect is associated not with a single capacitory potential but with a family of pinned solutions. We describe this family completely. Building on the single-site analysis, we establish an asymptotic expansion for the interval of pinned slopes for our original problem. The expansion decomposes the pinning effects into the nonlocal effect of the distribution of the defects and the local effect of the defect shape. As a consequence we prove that the pinning interval is nontrivial for a generic family of defects. Our work is motivated by the contact angle hysteresis phenomena in capillary contact lines. In particular our expansion provides the first mathematically rigorous justification of the formal analysis featured in \cite{Joanny-deGennes}, for general periodic defects.

\end{abstract}

\maketitle

\setcounter{tocdepth}{1}
\tableofcontents

\section{Introduction}

In this paper we study the {contact angle hysteresis} phenomena \cite{eral2013contact}, arising in capillary droplets on rough surfaces. On a perfectly smooth surface the balance of surface tension between the liquid, air and the surface yields a constant contact angle between the drop and the surface along the contact line, yielding axisymmetric profiles as the only stationary droplets on a flat surface. However, in practice, small scale inhomogeneities (or defects) are present on the surface. This small scale heterogeneity leads to intriguing non-round stationary droplet shapes and a range of large-scale contact angles \cites{caffarelli-mellet,feldman-smart,feldman2024obstacleapproachrateindependent,Bico_2001,Courbin_2007,Chu_2010,Raj_2014,eral2013contact}. Dynamically, a slowly evaporating or spreading drop will get pinned along some directions of propagation, leading to the aforementioned non-round geometries \cite{kim-contact}. Energetically, the energy landscape is wiggly at the microscopic scale, and thus there are many local minimizers. The challenge is thus to describe the range of macroscopic contact angles consistent with micro-scale local energy minimality.
This puts us into the realm of non-minimizing stationary solutions of free boundary problems, which is relatively less studied and understood in the literature.

We will work with the Bernoulli free boundary model
\begin{equation}\label{e.model}
\begin{cases}
\Delta u = 0 & \hbox{in } \ \{u>0\};\\
|\grad u| =Q(x) & \hbox{on } \partial \{u>0\}.
\end{cases}
\end{equation}
Here  $u:\R^d\to [0,\infty)$ denotes the height of the liquid, $\{u>0\}$ is the wetted set, and the free boundary $\partial\{u>0\}$ corresponds to the contact line. While we will work in all $d \geq 2$,  the case $d=2$ is most important since it is the physical dimension for the capillary problem. The slope condition $|\grad u| = Q(x)$ at the free boundary corresponds to Young's contact angle condition in capillarity.  

The interval of pinned slopes in the setting of \eqref{e.model} has been introduced and studied via a homogenization approach in \cites{kim-contact,feldman2021limit,CaffarelliMellet-StickingDrop}. However the pinning mechanism is both nonlinear and highly dependent on the distribution of $Q$. For this reason it has only been computed analytically in the setting of laminar media \cites{schwartz1985contact,CaffarelliMellet-StickingDrop,kim-contact,feldman2021limit,DeSimoneGrunewaldOtto}, and an exactly solvable discrete model \cite{feldman-smart}. Similar issues arise in other models, for example in the quenched Edwards-Wilkinson model and non-local variants \cites{CourteBattacharyaDondl,CourteDondlStefanelli,DirrYip} and curvature flows \cite{BraidesScilla}.

Thus, our aim in this paper is to understand the nonlocal and local effects of the medium $Q$ on the pinning interval in a general periodic media. The goal is to understand the anisotropic nature of the pinning phenomena beyond specific exactly solvable scenarios. As a first step in this program we study the ``dilute defect" setting, where the $Q$ has a sparse distribution of compactly supported inhomogeneities. This setting was originally introduced and studied in the physics literature by Joanny and De Gennes \cite{Joanny-deGennes}, and has been used widely in physics and engineering literature \cites{raphael1989dynamics,ConstantinLeDoussalWiese,joanny1990motion,quere2008wetting,reyssat2009contact}. \par

While \cite{Joanny-deGennes} argues via formal linearization, the PDE determining the contact line in \eqref{e.model} is highly nonlinear, non-local, and non-graphical. Topological singularities do occur during evolution, for example when two components of the wetted set $\{u>0\}$ advance and meet the local blow-up profile is expected to be a ``two-plane solution" $|x \cdot e|$ (see more discussion on this below \dref{proper} and at the end of \sref{literature} below). A mathematically rigorous justification of the asymptotic expansion in \cite{Joanny-deGennes} requires understanding the geometry of the free boundary both near and far from the defects, with attention to the possibility of non-graphical singularities arising at intermediate length scales. This is what we accomplish, for dilute and periodic distributions of defects.

We carry out our analysis in two parts. First we study the {\it local pinning effect} of a single defect (\sref{intro-single-site}): a complete description is given for both the family of pinned solutions of the single-site problem and the interval of ``nonlinear capacities" they generate. We further prove that the pinning interval is nontrivial for a generic family of defects. Second we study {\it the nonlocal geometry } of the free boundary in the region between the defects in periodic media (\sref{periodic-media-intro}). We compute the asymptotic expansion of the pinning interval, whose local factor is exactly the nonlinear capacity from the first part. We address each of these issues via respective auxiliary cell problems. These cell problems are adopted from the theory of perforated domain homogenization, but their application in the study of pinning is a new contribution of this work (see \sref{literature} for more discussion).
Let us point out that our analysis builds on the far-field expansion of solutions of the single-site problem, which is proved in the companion paper \cite{FKexterior}.

 Our analysis will have consequences that are of independent interest, to be discussed below. In particular it leads to the proof that the pinning phenomena is generic (see Theorem~\ref{t.main-pinning}), addressing a question that was open for a while.

\subsection{Statement of main results and ideas}
Now we make precise the setting and main results of the paper. 

\subsubsection{The single-site problem}\label{s.intro-single-site}
 Suppose $q: \R^d \to (-1,\infty), q\in C_0(\overline{B_1(0)})$. The single site problem is given as $u: \R^d \to \R^+ $ solving
\begin{equation}\label{e.single-site}
     \begin{cases}
         \Delta u = 0 & \hbox{in } \{u>0\}  \\
         |\grad u| = 1+q(x) &\hbox{on } \partial \{u>0\} .
     \end{cases}
\end{equation}
The local pinning effect of the defect will appear through the far-field asymptotic expansion of  solutions of \eref{single-site}. While this is reminiscent of capacities from the study of elliptic PDEs exterior domains, there is a significant difference: since we have a free boundary problem, there is not just one but instead a family of capacitory potentials and capacities associated with a given $q$. Describing that family is the purpose of the first part of this paper.

As in the companion paper \cite{FKexterior}, our analysis will be focused on {\it proper solutions} of \eqref{e.single-site}:
\begin{definition}\label{d.proper}
 $u \in C(\R^d)$ is \emph{proper} with direction $e$ if the sequence $u_r(x) : = r^{-1}u(rx)$ converges to $(x \cdot e)_+$ as $r\to 0$ in the following sense:
\begin{equation}\label{e.s2-blow-down-cond}
(u_r(x) - x\cdot e){\bf 1}_{\{u_r>0\}} \to 0 \ \hbox{ locally uniformly in $\R^d$ as } r\to 0.
\end{equation}
\end{definition}

This rather strong definition is introduced so that we can rule out solutions that blow down to singular solutions. We will show that it is enough to consider proper solutions to obtain the expansion \eqref{formula}. The main concern is the two-plane solution $|x \cdot e|$, a well-known singular profile in the Bernoulli problem \cites{CS,KriventsovWeiss,JerisonKamburov2}. Ruling out this scenario is most challenging in $d=2$ due to the far-field logarithmic growth of $u$ away from the planar profile.  An important step in achieving this is \pref{ukR-supersolution-extension}, which rules out the occurrence of two-plane blow-downs for advancing solutions in $d=2$.  

 By the far-field expansion proved in the companion paper \cite{FKexterior}*{Theorem 4.1}, every proper solution has a well-defined capacity $k(u)$ and capacitory-type expansion at $\infty$. Our first pair of main results address the attainable capacities.

In all theorems presented in this section, we fix $q$ and normalize $e = e_d$ by rotation. We begin with $d=2$.

\begin{theorem}[\tref{monotone-family2}] \label{t.main-2d-capacities}
    Let $d=2$. For any  $e \in S^{d-1}$,  there exists $-\infty<k_{\textup{rec}}(e) \leq 0 \leq k_{\textup{adv}}(e) \leq +\infty$ such that  for any proper solution $u$ of \eref{single-site}, there is $k\in [k_{\textup{rec}},k_{\textup{adv}}]$ such that
   \begin{equation}\label{e.main-2d-capacities-e1}
       u(x)  = x\cdot e +k\log|x| +O(1)\ \hbox{ as } \ \{u>0\}\ni x\to \infty.
   \end{equation} 
   Conversely, for any $k\in(k_{\textup{rec}},k_{\textup{adv}})$ there is a proper solution of \eqref{e.single-site} with  \eref{main-2d-capacities-e1}. 
\end{theorem}
To prove the above theorem, we identify the interval of capacities with the limit of a ``pulling experiment": the free boundary of a single-site problem solution is pulled through the defect by boundary datum $(x_d + k \log R)_+$ on $\partial B_R$. Denoting $\kappa^R_{\textup{adv}/\textup{rec}}$ of $k$ as the critical values of capacities at which the free boundary detaches from the support of $q$, we show that $\kappa^R_{\textup{adv}/\textup{rec}} \to k_{\textup{adv}/\textup{rec}}$ as $R \to \infty$.
 Interestingly, it is unclear whether $k_{\textup{adv}} < +\infty$ in general, with the exception of small $\max q$.
This challenge, specific to the advancing capacity in $d=2$, is tied to the  logarithmic tail of solutions.

Next we turn to $d\geq 3$. Here the behavior of proper solutions is a bit different due to the thin tail of the fundamental solution: there is a monotone family of  solutions that goes through the defect and varying capacities can be achieved at varying heights.

\begin{theorem}[\tref{monotone-family}]\label{t.main-3d-capacities} 
Let $d\geq 3$. There exist functions $\kappa_{\textup{rec}}(\cdot;e), \kappa_{\textup{adv}}(\cdot;e): \R \to \R$ with $\kappa_{\textup{rec}} \leq 0 \leq \kappa_{\textup{adv}}$ such that for any proper solution $u$ of \eref{single-site} there is $s \in \R$ and $k\in [\kappa_{\textup{rec}}(s;e), \kappa_{\textup{adv}}(s;e)]$ such that
\begin{equation}\label{expansion:d>2}
   u(x)  = x\cdot e +s-k|x|^{2-d} +O(|x|^{1-d})\ \hbox{ as } \ \{u>0\} \ni x \to \infty. 
   \end{equation}
  Moreover, for all $s \in \R$ there is a minimal supersolution $u_{\textup{adv}}(s)$ with $k=\kappa_{\textup{adv}}(s)$, and for all except countably many $s \in \R$ there is a maximal subsolution  $u_{\textup{rec}}(s)$ achieving \eqref{expansion:d>2} with $k=\kappa_{\textup{rec}}(s)$. 
\end{theorem}

See \tref{capacity-properties} for further properties of $\kappa_{\textup{rec}}$ and $\kappa_{\textup{adv}}$. For $d\geq 3$, the extremal solutions $u_{\textup{adv}/\textup{rec}}(\cdot;s)$ form monotone families, laminations of the whole space passing through the defect at every height $s$, which sandwich every proper solution. The lamination is genuinely ``defective": the maps $s \mapsto u_{\textup{adv}/\textup{rec}}(\cdot;s)$ have countably many jump discontinuities, and the jumps are exactly the pinned gaps, regions that the free boundary skips as it is pushed through the defect.

\subsubsection{Periodic media}\label{s.periodic-media-intro}
 As a model example for a surface with dilute defects, we will consider $Q$ of the form
\[Q(x) = Q_\delta(x) = 1+\sum_{z \in \Z^d} q(\tfrac{x-z}{\delta}), \quad \delta\ll 1, \]
where the defect function $q$ is as given before.
 To state our results, we briefly recall the definition of the pinning interval for \eref{model} from \cite{feldman2021limit} (see Section \ref{s.pinning-interval-recall} for a full discussion): $\alpha>0$ is a {\it pinned slope} at normal direction $e\in S^{d-1}$  if there is a solution $v$ of
\begin{equation}\label{e.model2}
\begin{cases}
\Delta v = 0 & \hbox{in } \ \{v>0\};\\
|\grad v| = Q_\delta(x) & \hbox{on } \partial \{v>0\};\\
\sup_{x \in \{v>0\}} |v(x) - \alpha (e \cdot x)| < + \infty.
\end{cases}
\end{equation}
The pinned slopes at direction $e$ then form a non-empty, closed interval, which we call the {\it pinning interval}: $[Q^\delta_{\textup{rec}}(e),Q^\delta_{\textup{adv}}(e)]$.   Here $Q^\delta_{\textup{rec}}(e)$ is called the \emph{receding slope} and $Q^\delta_{\textup{adv}}(e)$ is called the \emph{advancing slope}. This terminology arises from the dynamic manifestation of the pinning phenomenon. The slope is smaller than $Q^\delta_{\textup{rec}}(e)$ where the contact line recedes, and is larger than $Q^\delta_{\textup{adv}}(e)$ where the contact line advances: see for example \cites{kim-contact,AlbertiDeSimone}.

To obtain an asymptotic expansion for the pinning interval endpoints, the nonlinear capacities from \sref{intro-single-site}, will play an important role. In dimension $d=2$ $k_{\textup{adv}/\textup{rec}}(e)$ are defined above in \tref{main-2d-capacities}. In dimension $d \geq 3$ we define
 in terms of $\kappa_{\textup{adv}/\textup{rec}}(s)$ from \tref{main-3d-capacities},
\[k_{\textup{adv}}(e) := \max_{s \in \R} \kappa_{\textup{adv}}(s;e) \ \hbox{ and } \ k_{\textup{rec}}(e) := \min_{s \in \R} \kappa_{\textup{rec}}(s;e).\]

\begin{theorem}[Propositions \ref{p.asymptotic-exp-lower-bound} and \ref{p.asymptotic-exp-upper-bound}]\label{t.asymptotics}
    Let $\xi \in \Z^d$ irreducible and $e := \frac{\xi}{|\xi|}$. Then as $\delta \to 0$
     \begin{equation}\label{formula}
        \begin{split}
            &Q^\delta_{\textup{adv}}(e) = 1 +\gamma_d\delta^{d-1}|\xi|^{-1}k_{\textup{adv}}(e)+o(\delta^{d-1}), \hbox{ and } \\
    &Q^\delta_{\textup{rec}}(e) = 1 +\gamma_d\delta^{d-1}|\xi|^{-1}k_{\textup{rec}}(e)+o(\delta^{d-1}).
        \end{split}
    \end{equation}
Here $\gamma_2 = \pi$ and $\gamma_d = \frac{1}{2}d(d-2)|B_1|$ for $d \geq 3$.

    \end{theorem}
 
 The expansion \eqref{formula} has a natural product form: the width of the hysteresis interval equals the density of defects met by the contact line, times a pinning strength per defect, up to $o(\delta^{d-1})$. This product structure is already present in the formal analysis of \cite{Joanny-deGennes}. What is new in this theorem is  the identification of the local factor as the nonlinear capacity of the single-site problem, the anisotropic lattice factor $|\xi|^{-1}$, and control of the free boundary at the intermediate scales where the formal matching could fail. In this sense \eqref{formula}  separates the nonlocal singular anisotropy caused by the lattice structure from the local effect of the individual defects as represented by $k_{\textup{adv}/\textup{rec}}$.

Using \eqref{formula} and our results on single site problem, we can show that the pinning interval is nontrivial for all rational directions, for a generic family of $q$.  In particular it follows that $k_{\textup{adv}}$ and $k_{\textup{rec}}$ are generically non-zero.
Here, by the word generic, we mean that it occurs at least on an open set of defect coefficients $q \in C_c(B_1)$. While physical intuition suggests that the pinning interval should be nontrivial at all directions for typical media, our result gives the first rigorous confirmation of this in a general setting (see \cite{feldman-smart} for the corresponding result in an exactly solvable discrete model). This is because the main computable example in the continuum, laminar media, is exceptional in this regard and has pinning only at the laminar direction. 

\begin{theorem}[\tref{nonzero} and \cref{nonzero-ks}]\label{t.main-pinning}
The pinning interval $[Q^{\delta}_{\textup{rec}}(e), Q^{\delta}_{\textup{adv}}(e)]$ is non-trivial for all rational $e\in \mathcal{S}^{d-1}$ if (a) $q$ from \eqref{e.model} is nontrivial and has sign,  or if (b) $q$ changes sign with $\int q dx \neq 0$ and with a small Lipschitz constant.
\end{theorem}

Below we describe the main ideas leading to \tref{asymptotics}, as well as two consequences of our analysis that are of independent interest.

 It is previously known that $Q^\delta_{\textup{adv}}(e)$ and $Q^\delta_{\textup{rec}}(e)$ are achieved by strong Birkhoff solutions $u^{\delta}$ of the periodic problem \eqref{e.model} (see \sref{pinning-interval-recall}). Hence we study the behavior of $u^{\delta}$.

 When $\delta$ is small, we show that the free boundary of $u^\delta$ is close to a hyperplane $H$ with normal $e$. This means that $u^\delta$ only feels the effect of the defects $\Z^d \cap H$, and near these defects we expect $u^\delta$ to look like a solution of the single site problem. Namely we expect, for each $z \in \Z^d \cap H$,
\begin{equation}\label{e.heuristic-inner}
    \frac{1}{\delta}u^\delta(z+\delta x) \approx u_{\textup{in}}(x) \ \hbox{ which is a solution of \eref{single-site}.}
\end{equation}
We will show that $u_{\textup{in}}$ is ``proper", morally this is happening in \tref{monotone-family2} below: then Theorem \ref{t.main-2d-capacities} or \ref{t.main-3d-capacities} applies to conclude that $u_{\textup{in}}$ has a far-field capacity $k \in [k_{\textup{rec}},k_{\textup{adv}}]$.

To complete the asymptotic expansion of $u^\delta$ we must describe the behavior of $u^\delta$ in the outer region away from the holes. We will show that, away from the defects,
\begin{equation}\label{e.heuristic-outer}
    u^{\delta}(x) \approx (1+\gamma_d|\xi|^{-1}k\delta^{d-1}) x_d + \delta^{d-1} k \omega(x)
\end{equation}
Here the ``corrector" $\omega(x)$ is the solution of an auxiliary cell problem 
\begin{equation}\label{cell:eq}
    \Delta \omega = 0 \,\, \hbox{ in } \{x\cdot e >0\}; \quad
    \partial_{x_d} \omega = 1 - \sum_{z \in \Z^d \cap H}|\xi|\delta_z \,\, \hbox{ on } H:=\{x\cdot e=0\}.
\end{equation}
Recall that we have written $e = \xi/|\xi|$ for an irreducible vector $\xi \in \Z^d$, so $|\xi|$ is the area of the fundamental domain of the lattice $\Z^d \cap H$. This means that the average of $\partial_d\omega$ is $0$ on $H$ and so the slope of $\omega$ at infinity is zero as well. The fundamental solution type singularities at the lattice sites arises from matching the far-field behavior of the inner solution $u_{\textup{in}}$.

The above heuristic expansions \eqref{e.heuristic-inner} -\eqref{e.heuristic-outer} are used in two ways to deliver \tref{asymptotics}: first to construct barriers to obtain a lower bound on the pinning interval (Section \ref{s.homogenization-lower-bound}), and second to study the asymptotics of arbitrary plane-like solutions $u^\delta$, in order to obtain a matching upper bound on the pinning interval (Section \ref{s.homogenization-upper-bound}). As a byproduct we obtain a decomposition of linear and nonlinear part of $u^{\delta}$ as $\delta\to 0$, which appears to be new in the analysis of free-boundary problems:

\begin{theorem}\label{t.udelta-expansion}
 (see \sref{homogenization-capacity}) For a rational $e\in\mathcal{S}^{d-1}$, suppose $k_{\textup{adv}}(e) \in (0,\infty)$. If $u^\delta$ are strong Birkhoff plane-like solutions with  slope $Q^\delta_{\textup{adv}}(e)$, then there are $s_\delta \to 0$ so that, modulo a period translation of $u^{\delta}$ to re-center the free boundary near $x_d=0$,
\begin{align*}
    \lim_{\delta \to 0} \dfrac{1}{\delta^{d-1}}[Q^\delta_{\textup{adv}}(e)\,x\cdot e+s_\delta -u^\delta(x)]= \gamma_d|\xi|^{-1}k_{\textup{adv}}(e)\omega(x),
\end{align*}
where the convergence is locally uniform. Moreover $\omega$ solves \eqref{cell:eq}. A parallel result holds for strong Birkhoff plane-like solutions with the minimal slope, with $Q^{\delta}_{\textup{adv}}$ and $k_{\textup{adv}}$ replaced by $Q^\delta_{\textup{rec}}$ and $k_{\textup{rec}}$.
\end{theorem}

\begin{remark} 
We will reduce the problem to the case $e = e_d$ for the rest of the paper. For the single site problem this normalization can be achieved by rotating $q$ accordingly. In Sections \ref{s.periodic-prelim}-\ref{s.homogenization-upper-bound} where we consider the full periodic setting, we will reduce to the case $e=e_d$ by rotation, at the cost of considering a rotated periodicity lattice. 
\end{remark}

\subsection{Comparison with previous literature}\label{s.literature}

Below we discuss the contribution and relevance of the present work in the context of perforated domain homogenization, the rigorous pinning literature, and free boundary regularity theory.

\subsubsection*{A family of capacities, not just one} In perforated domain homogenization the single-site problem behind the ``strange term" of Cioranescu and Murat \cite{cioranescu1997strange}, and behind dilute limits going back to \cite{PapanicolaouVaradhan1980}, is a linear exterior problem: the hole carries exactly one capacitory potential: see e.g. \cites{CaffarelliMellet-RandomObstacle,heida2020,Cioranescu-Donato,CioranescuUnfolding}. The uniqueness of the capacitory potential is what makes the classical theory linear at the cell level. Since we have a free boundary problem, the defect supports a genuine continuum of pinned states, and the physically meaningful quantity is the \emph{spread} $[k_{\textup{rec}},k_{\textup{adv}}]$ of their far-field coefficients. Theorems \ref{t.main-2d-capacities} and \ref{t.main-3d-capacities} describe this family completely. We are not aware of an antecedent for such a description in either the perforated domain or the exterior Bernoulli literature. Viewing the periodic array of defects as a perforated domain, our analysis also establishes, for the first time, a homogenization result in which the perforations are localized in the height variable: the finite height of the ``perforation" means that the drop is allowed to sit on the defects, a new feature even compared to other free boundary problems in perforated domains \cites{guillen-kim,AbedinFeldman}. In our setting each hole is governed by the whole family of pinned states rather than a single capacitory potential, and selecting the correct member of the family couples the local and nonlocal problems; the analysis of arbitrary plane-like solutions in \sref{homogenization-upper-bound} is a correspondingly new feature caused by the pinning phenomenon.

\subsubsection*{The rigorous pinning literature} 
\tref{asymptotics} is, to our knowledge, the first mathematically rigorous justification of the dilute expansion of Joanny and De Gennes \cite{Joanny-deGennes} for general defects in the periodic setting; the rigorous results known before are the qualitative hysteresis interval of \cite{caffarelli-mellet}, the laminar computations \cites{schwartz1985contact,CaffarelliMellet-StickingDrop,kim-contact,feldman2021limit,DeSimoneGrunewaldOtto}, and the discrete model of \cite{feldman-smart}. In this context \tref{main-pinning} resolves a question left open in \cite{caffarelli-mellet}, where a nontrivial homogenized contact angle interval was proved in a chosen direction and pinning in all directions was flagged as likely out of reach of the methods there: we obtain all rational directions at once, on an open set of defect profiles, with a rate.  Let us also point out that our formula \eqref{formula} shows that the \emph{static} interval of pinned capacities of a single defect is nontrivial at first order in the defect strength, which is a distinctive statement not present in \cite{Joanny-deGennes}.

\subsubsection*{Free boundary regularity theory and non-minimizing solutions} 
 Sections \ref{s.prelim-single-site} and \ref{s.strong-max} develop the stability and rigidity tools which play a key role in our homogenization results.
The challenge in establishing stability property for extremal solutions is the lack of uniform non-degeneracy for maximal subsolutions in $d\geq 3$, and
general viscosity solutions \cite{KamburovWang}. Without it the class of viscosity subsolutions is not closed under uniform limits. 
To overcome this we introduce a relaxation in  \sref{prelim-single-site} by a new notion of \emph{relaxed subsolutions}, which is closed under uniform limits but allows for degeneracy. It is inspired by similar  ideas in the class of inner variational solutions \cite{Weiss}.  This new notion allows us to make several kinds of compactness arguments in the relaxed viscosity solution class. 
In particular going through the relaxed notion allows us to obtain a strong maximum principle with a defect (\tref{defect-strong-max}) for viscosity solutions. This result will be an important ingredient in Section 4. The closest result to \tref{defect-strong-max} is the effective removable singularities theorem of Jerison and Kamburov \cite{JerisonKamburov2} for classical solutions in $d=2$, proved via conformal mapping and the planar classification of global solutions.

\subsection{ Open questions} 
We finish this section with some open questions.  

\subsubsection*{Two-plane solutions} While it is sufficient for our purposes to describe the pinning interval only in terms of ``flat" (or proper) solutions, there could also be two-phase like solutions of \eref{single-site} that are pinned on the defect. The challenge of understanding two-plane solutions is a serious issue in general for non-minimizing Bernoulli solutions \cites{JerisonKamburov2,KriventsovWeiss}. It also is an issue in the minimal surface capillary model, posing an obstacle for the description of hydrophilic advancing angles as well as hydrophobic receding angles: 
see e.g. \cite{DeRosa25}.

\subsubsection*{Non-periodic Media} Some important non-periodic settings include the irrational directions $e\in\mathcal{S}^{d-1}$ or the case of dilute defects on a Poisson point cloud. In both of these cases there is no longer a $(d-1)$ dimensional periodic structure on the hyperplane $\{x: x\cdot e=0\}$. Instead the hyperplane sees the defects in an $O(\delta)$ neighborhood. Due to the sparsity of lattice sites near the irrational hyperplane we expect the leading order in \eqref{formula} to be order $\delta^{d}$ instead of $\delta^{d-1}$, and we have a similar expectation in the case of Poisson clouds. We expect that to address these cases the theory of the single defect problem developed in this paper should be combined with new ergodic theory or probabilistic ideas.

\subsubsection*{Capillary free boundary problem} Our problem \eqref{e.model} arises from the capillary free boundary problem in the limit as the contact angle approaches $0$ or $\pi$ capillary energy. In that limit $u$ represents a normalized droplet height. For more details see \cite{Chodosh-Li}, which obtains the first rigorous application of this limit procedure, also see \cite{FK}. We expect parallel results to hold for the capillary free boundary problem, but leave this as an open problem.

\subsection{Acknowledgments}
We thank Norbert Po\v{z}\'ar for providing the code used for the simulations shown in \fref{2d-pinning-family}. W.F.'s research is partially supported by NSF DMS-2407235. Part of this work was completed during W.F.'s visit at the ESI, he thanks the ESI for hosting him. I.K.'s research is partially supported by NSF DMS 2452649. Part of this work was completed during I.K’s visit at KIAS, and she thanks KIAS’s hospitality.

\subsection{Declaration of AI use} The original version of this paper \cite{FKoriginal} was written without AI assistance. For this revised version we used AI tools for editing and improving the presentation.  In particular \pref{strong-max-flat} is strengthened based on a strategy suggested by Claude Fable 5. The remainder of the mathematical ideas in the paper were generated in the old-fashioned human way.

\section{Setting and viscosity solution theory}\label{s.prelim-single-site}

We recall the theory of viscosity solutions of the Bernoulli free boundary problem in an open set $U$:
\begin{equation}\label{e.bernoulli-basic}
     \begin{cases}
         \Delta u = 0 & \hbox{in } \{u>0\} \cap U \\
         |\grad u| = Q(x) &\hbox{on } \partial \{u>0\} \cap U. 
     \end{cases}
\end{equation}
Here we assume that $Q, u \in C(\bar{U})$, and $Q>0$, $u\geq 0$, and $Q$ is bounded. Often  we will be able to reduce to the homogeneous version 
\begin{equation}\label{e.bernoulli}
    \begin{cases}
        \Delta u = 0 & \hbox{ in } \{u>0\} \cap U\\
        |\grad u| = 1 &\hbox{ on } \partial \{u>0\} \cap U.
    \end{cases}
    \end{equation}

    Viscosity solutions in the context of energy minimizers and minimal supersolutions for the above problems have been well-studied, see for example \cite{CS}. On the other hand much less is known on general viscosity solutions. Maximal subsolutions, in particular, play a central role in studying pinning in the context of Bernoulli and capillary free boundary conditions. This physical motivation suggests they should be better understood from a pure mathematical perspective.  The main difficulty in this study lies in the lack of uniform non-degeneracy estimates and the consequential loss of uniform stability.
    
    In this section we introduce a generalized notion of subsolution which allows uniform stability. This relaxed notion and the corresponding stability turns out to work well for the applications we consider.

\subsection{Bernoulli free boundary problems and solution notions}
We will begin with the following definition, which is classical when $E=U$.
\begin{definition}
    For $E$ relatively closed in $U$, we say $\psi$ \emph{touches $u$ from above in $E$ (strictly) at $x_0$} if $u-\psi$ has a (strict) local maximum zero in $E$ at $x_0$. Touching from below is defined analogously with the minimum. When no set $E$ is specified it means $E = U$.
\end{definition}

Now we define viscosity sub- and supersolutions.

\begin{definition}\label{d.supersolution}
    $u \in C(U)$  is a \emph{supersolution} of \eref{bernoulli-basic} if $u\geq 0$ and whenever $\varphi \in C^\infty(U)$ touches $u$ from below at $x \in U$, either $\Delta \varphi(x) \leq 0$, or $\varphi(x) = 0$ and $|\grad \varphi(x)| \leq Q(x)$.
\end{definition}

\begin{definition}\label{d.subsolution}
    $u \in C(U)$ is a \emph{subsolution} of \eref{bernoulli-basic} if $u\geq 0$ and whenever $\varphi \in C^\infty(U)$ touches $u$ from above in $\overline{\{u>0\}} \cap U$ at $x$, either $\Delta \varphi(x) \geq 0$, or $\varphi(x) = 0$ and $|\grad \varphi(x)| \geq Q(x)$.
\end{definition}

\begin{definition}\label{d.viscosity-solution}
   $u \in C(U)$ is a \emph{viscosity solution} (or just \emph{solution}) of \eref{bernoulli-basic} if $u$ is a supersolution and a subsolution of \eref{bernoulli-basic}.
\end{definition}

Note that for viscosity subsolution we allow test functions which only touch from above in $\overline{\{u>0\}} \cap U$ instead of globally: this is due to the fact that smooth test functions cannot touch $u$ from above globally at free boundary points. As a consequence stability of the subsolutions involve convergence of the positivity sets, which however may degenerate in the limit: for instance  $u_n(x) :=\alpha(x_d)_+ + \frac{1}{n}(-x_d)_+$ is a sequence of subsolutions for any $0<\alpha<1$ that fails to converge to a subsolution. In order to obtain unconditional stability of the subsolution notion we must introduce a relaxation.

\begin{definition}\label{d.degenerate}
    Let $u \geq 0$ be in  $C(U)$. Consider a set $E$ containing $\{u>0\}$ such that $E \cap U$ is relatively closed in $U$. The pair $(u,E)$ is a \emph{{relaxed subsolution}} of \eref{bernoulli-basic} if whenever $\varphi \in C^\infty(U)$ touches $u$ from above in $E \cap U$ at $x$, either $\Delta \varphi(x) \geq 0$, or $\varphi(x) = 0$ and $|\grad \varphi(x)| \geq Q(x)$.

    Say that $u$ is a \emph{{relaxed solution}} of \eref{bernoulli-basic} if it is a supersolution and a {relaxed subsolution}.
\end{definition}

Note that $u$ is a subsolution if and only if $(u, \overline{\{u>0\}})$ is a {relaxed subsolution}.
Moreover the set of subsolutions are a strict subset of the {relaxed subsolution}s: for instance $(\alpha (x_d)_+, \R^d)$ is a {relaxed subsolution} of \eref{bernoulli-basic} for every $\alpha>0$.

\begin{remark}\label{r.subsolutions-subharmonic}
    If $(u,E)$ is a {relaxed subsolution}, then $u$ is subharmonic in $U$:  a smooth $\varphi$ touching $u$ from above at $x_0 \in \{u=0\}$ has a local minimum there, and thus $\Delta \varphi(x_0) \geq 0$. If the touching is in $\{u>0\}$ then $u(x_0) = \varphi(x_0)>0$ and we can apply \dref{degenerate}.
\end{remark}

\begin{remark}\label{r.remove-strictness-issue}
    The slope conditions at the free boundary only need to be checked when the interior PDE condition fails. However by a standard perturbation argument one can show a stronger property. For example say that $(u,E)$ is a {relaxed subsolution} in $U$. Then one can check that whenever $\varphi \in C^\infty(U)$ touches $u$ from above in $E$ at $x$ with $\varphi(x) = 0$ and $|\grad \varphi(x)|>0$ then $|\grad \varphi(x)| \geq Q(x)$.  A similar conclusion holds for supersolutions. This is a folk-lore result in the literature, the argument is the same from, for instance, \cite{feldman-kim-pozar}*{Lemma A.4}.
\end{remark}

\subsection{Uniform stability} Here we prove the stability of {relaxed subsolution}s  under local uniform limits.

Recall the Kuratowski upper limit relative to an open set $U$: for $E_n \subset U$,
\begin{equation}\label{limsup_set} E^*:={\limsup}^* E_n := \{x \in U : x_k\to x \hbox{ for some } x_k \in E_{n_k}, n_k \to \infty\},
\end{equation}
a relatively closed subset of $U$.

\begin{theorem}[Stability]\label{t.stability}
    Let $U$ be an open set, and $u_k \in C(U)$ with $u_k \to u$ locally uniformly in $U$.
    \begin{enumerate}[label = \textup{(\roman*)}]
        \item\label{part.stability-super} If each $u_k$ is a supersolution of \eref{bernoulli-basic} in $U$, then so is $u$.
        \item\label{part.stability-pair} If each $(u_k, E_k)$ is a {relaxed subsolution} of \eref{bernoulli-basic} in $U$, then so is $(u, E^*)$, where $E^*$ is as given in \eqref{limsup_set}.
    \end{enumerate}
\end{theorem}
In particular a uniform limit of relaxed solutions, in the sense just explained for relaxed subsolutions, is a relaxed solution.
\begin{proof}
The supersolution part is standard (see e.g.  \cite{FeldmanPozar2024}*{Lemma 3.1}).

We proceed to (ii).  Note that $\{u>0\} \subset E^*$ since $u(x)>0$ implies that $u_k(x)>0$ for $k$ sufficiently large. Suppose that $\varphi$ smooth touches $u$ from above in $E$ at $x_0=0$. If $u(0)>0$ then the local uniform convergence implies that $u_k(x) >0$ in a neighborhood of $0$ for sufficiently large $k$, thus they are subharmonic in that neighborhood, and then $\Delta \varphi(0) \geq 0$ follows from a classical argument. So we reduce to the case $\varphi(0) = 0$ and $\Delta \varphi(0) <0$. 

 Next we strictify the touching by replacing $\varphi$ with $\psi(x):=\varphi(x) +\delta |x|^2$ with $0<\delta \leq -\frac{1}{4d}\Delta \varphi(0)$. Note that $\psi(0)=0$, $|\grad \psi(0)| = |\grad \varphi(0)|$,  and $\Delta \psi(0)<0$ still. Now $\psi$ touches $u$ from above in $E^*$ strictly at $0$: namely for a small $r>0$, 
 \begin{equation}\label{max_strict}
 u(0)=\psi(0) \ \hbox{ and } \ (u(x)-\psi(x))<0 \ \hbox{ for } \ x \in (E^* \cap \overline{B_r(0)}) \setminus \{0\}.
 \end{equation}
By definition of the $E^*$, along a subsequence, there is a sequence of $x_k \in E_k$ with $x_k \to 0$. {Let us choose $r$ small enough such that $\overline{B_r(0)}\subset U$.} Since $E_k \cap \overline{B_r(0)}$ is compact there is $y_k \in E_k \cap \overline{B_r(0)}$ so that
\[c_k:=(u_k-\psi)(y_k)= \max_{E_k \cap \overline{B_r(0)}} (u_k-\psi). \]
Then, since $x_k\to 0$,
\[ \liminf_{k \to \infty} c_k = \liminf_{k \to \infty}(u_k-\psi)(y_k) \geq \limsup_{k \to \infty}(u_k-\psi)(x_k) =0.\]
If  $y_{k_j} \to y$ as $j\to\infty$ along a subsequence $k_j$, then $y \in E^*$ by \eqref{limsup_set}.  So 
\[0 \leq \lim_{j \to \infty} c_{k_j} = \lim_{j \to \infty} (u_{k_j}-\psi)(y_{k_j}) = (u-\psi)(y). \]
Hence from \eqref{max_strict} it follows that $y = 0$. Thus we conclude that the full sequence $y_k$ converges to $0$ and $c_k\to 0$.

Now $\psi+c_k$ touches $u_k$ from above at $y_k \in E_k$. Since $\Delta \psi(y_k) \to \Delta \psi(0) <0$ we conclude from \dref{degenerate} that $u_k(y_k) = 0$ and
\[|\grad \psi(y_k)| \geq Q(y_k).\]
Taking the limit on both sides and using continuity of $\grad \psi$ and $Q$ and $y_k \to 0$, we conclude that $|\nabla \varphi|(0) \geq Q(0)$, and thus (ii) holds.
 
\end{proof}

The degeneracy can be removed in some situations.  For example when the sequence is monotonically increasing.
\begin{corollary}\label{c.monotone-stability}
    If $u_k$ is a monotonically increasing sequence of  subsolutions of \eref{bernoulli-basic} in a domain $U$ and $u_k \nearrow u$ locally uniformly in $U$ then $u$ is a subsolution of \eref{bernoulli-basic}.
\end{corollary}
\begin{proof}
    Let $E^*:= {\limsup}^*\overline{\{u_k>0\}}$. By \tref{stability} applied to the pairs $(u_k, \overline{\{u_k>0\}})$, $(u,E^*)$ is a {relaxed subsolution}. Since each $\overline{\{u_k>0\}} \subset \overline{\{u>0\}}$ we conclude $E^* \subset \overline{\{u>0\}}$ and therefore $u$ is a subsolution.
\end{proof}
Also limits of subsolutions are subsolutions when the sequence satisfies a uniform non-degeneracy hypothesis.
\begin{definition}\label{d.nondegen}
    Say that $u \in C(U)$, non-negative, is $\ell$-non-degenerate if for all $x \in \overline{\{u>0\}}$ and all $\overline{B_r(x)} \subset U$
    \[\max_{\overline{B_r(x)}} u \geq \ell r.\]
\end{definition}
\begin{lemma}\label{l.nondegen-stability}
    If $u_k$ is a sequence of  subsolutions of \eref{bernoulli-basic} in a domain $U$ which are all $\ell$-non-degenerate and $u_k \rightarrow u$ locally uniformly in $U$, then ${\limsup}^*\overline{\{u_k>0\}} \subset \overline{\{u>0\}}$ and therefore $u$ is a subsolution of \eref{bernoulli-basic}.
\end{lemma}
\begin{proof}
    By \tref{stability} it suffices to show that $E^*:={\limsup}^*\overline{\{u_k>0\}} \subset \overline{\{u>0\}}$. Let $x \in E^*$ and $r>0$ for any $N$ there is $k \geq N$ such that there is $x_k \in \{u_k>0\} \cap B_{r/2}(x)$. By the $\ell$-non-degeneracy of $u_k$ we conclude that $\max_{\overline{B_r(x)}} u_k \geq \frac{1}{2}\ell r$. Sending $N$ and hence $k(N) \to \infty$ and using the uniform convergence we find $\max_{\overline{B_r(x)}} u \geq \frac{1}{2}\ell r$. Since $r>0$ was arbitrary $x \in \overline{\{u>0\}}$ and we are done.
\end{proof}
Note that the weaker version of \dref{nondegen} which only produces non-degeneracy at $x \in \partial \{u>0\}$ instead of all $x \in \overline{\{u>0\}}$ does not suffice for the conclusion of \lref{nondegen-stability} as evidenced by the counter-example $u_k(x)=\frac{1}{2}(x_d)_+ + \frac{1}{k}(x_d)_-$.

\subsection{Lipschitz estimate and non-degeneracy properties}

Here we present some Lipschitz and non-degeneracy estimates. First the Lipschitz estimate, which only requires the viscosity solution property.
    \begin{lemma}[Lipschitz estimate]\label{l.lipschitz-estimate}
        There is $C(d) \geq 1$ so that if $u\in C(B_2)$ is harmonic in its positivity set and a supersolution of
\[\Delta u = 0 \ \hbox{ in } \ \{u>0\} \cap B_2 \ \hbox{ and } \ |\grad u| \leq \Lambda \ \hbox{ on } \ \partial \{u>0\} \cap B_2\]
 then
\[\|\grad u\|_{L^\infty(B_1)} \leq C(\Lambda+\osc_{B_2} u).\]
    \end{lemma}

We note the following non-degeneracy property at outer regular points of $\partial E$. 
\begin{lemma}\label{l.linear-growth-outer-reg}
    Let $(u,E)$ be a {relaxed subsolution} of \eref{bernoulli} in $B_r(0)$, with $0 \in \partial E$. In addition suppose that $E$ has an exterior tangent ball at $0$ of radius $1$. Then $0\in\partial\{u>0\}$ and 
    \[ \sup_{B_\rho(0)} u \geq c(d) \rho \ \hbox{ for all } \ 0 < \rho \leq \min\{r,1\}.\]
\end{lemma}
This follows from a standard barrier argument, with a radially increasing symmetric barrier which vanishes in the exterior tangent ball of radius $\rho/4$.

\begin{corollary}\label{c.bdry-E-contain}
    Suppose that $(u,E)$ is a {relaxed subsolution} of \eref{bernoulli} in $U$. Then $\partial E \subset \partial \{u>0\}$.
\end{corollary}
\begin{proof}
    Let $x_0 \in \partial E$, then $x_n\to x_0$ for some $x_n \in U \setminus E$. Let $y_n$ be the closest point in $E$ to $x_n$ with $r_n = |x_n-y_n|>0$. Then $y_n \in \partial E$ and has an exterior tangent ball of radius $r_n \to 0$. Then, by a rescaling of \lref{linear-growth-outer-reg}, $y_n \in \partial \{u>0\}$ and $y_n \to x_0$ as well. We conclude that $x_0\in\partial\{u>0\}$.
\end{proof}

\subsection{Perron extremal solutions}\label{s.perron-method}

The existence of minimal supersolutions and maximal subsolutions can be proven by Perron's method: see e.g. \cite{CS}*{Theorem 6.1}.

Next we introduce notions of minimal supersolution and maximal subsolution. 
\begin{definition}\label{d.weakly-minimal}
    Say that $u \in C(U)$ non-negative is a \emph{minimal supersolution} of \eref{bernoulli-basic} if 
    \begin{itemize}
        \item $u$ is a supersolution of \eref{bernoulli-basic} in $U$
        \item for all balls $B \subset\subset U$ and $v \in C(\overline{B})$ a supersolution of \eref{bernoulli-basic} in $B$ with $v > u$ on $\overline{\{u>0\}} \cap \partial B$ then $v \geq u$ in $B$. 
    \end{itemize}
\end{definition}
\begin{definition}\label{d.weakly-maximal}
    Say that $u \in C(U)$ non-negative is a \emph{maximal subsolution} of \eref{bernoulli-basic} if 
    \begin{itemize}
        \item $u$ is a subsolution of \eref{bernoulli-basic} in $U$
        \item for all balls $B \subset\subset U$ and $v \in C(\overline{B})$ a subsolution of \eref{bernoulli-basic} in $B$ with $v < u$ on $\overline{\{v>0\}} \cap \partial B$ then $v \leq u$ in $B$.
    \end{itemize}
\end{definition}
A relaxed maximal subsolution is defined analogously by requiring that $(u,E)$ is a {relaxed subsolution} instead of a subsolution in \dref{weakly-maximal}.

Note that these notions are somewhat weaker than the solutions which are typically obtained from Perron's method, due to the strict boundary ordering required.  Again this is a kind of relaxation which behaves better under uniform limits, as we will see below in \lref{extremal-stability}. 

First we make a clarifying note.
\begin{lemma}
    If $u$ is a minimal supersolution (resp. maximal subsolution) of \eref{bernoulli-basic} in $U$ then $u$ is a solution of \eref{bernoulli-basic} in $U$.
\end{lemma}
The proof is a typical argument from Perron's method which we skip, e.g. if a smooth strict supersolution $\varphi$ touched $u$ from above at some $x_0$ we could perturb to obtain strict ordering on some $\partial B_r(x_0)$ and then take $\min \{u,\varphi - \delta\}$ for small $\delta>0$ to get a smaller supersolution contradicting the minimality property in \dref{weakly-minimal}.

\begin{lemma}\label{l.extremal-stability}
Let $U$ be an open set, and $u_k \in C(U)$ with $u_k \to u$ locally uniformly in $U$.
\begin{enumerate}[label=(\roman*)]
    \item If each $u_k$ is a minimal supersolution in $U$ then so is $u$.
    \item If each $(u_k,E_k)$ is a relaxed maximal subsolution then so is $(u,E^*)$, where $E^* = {\limsup}^*E_k$.
\end{enumerate}
\end{lemma}
To prove the  lemma, we state the well-known non-degeneracy property of minimal supersolutions
\begin{lemma}[see for example \cite{CS}]\label{l.minimal-nondegen}
    If $u$ is a minimal supersolution of \eref{bernoulli-basic} in $U$ then $u$ is $\ell$-non-degenerate in the sense of \dref{nondegen} with $\ell=c(d)\min Q$.
\end{lemma}
\begin{proof}[Proof of \lref{extremal-stability}]
    First we do the minimal supersolution case. By \tref{stability} $u$ is a supersolution. By \lref{minimal-nondegen} and \lref{nondegen-stability} we get $E^* :={\limsup}^*\{u_k>0\} \subset\overline{\{u>0\}}$.   Let $B \subset\subset U$ and $v \in C(\overline{B})$ a supersolution of \eref{bernoulli-basic} in $B$ with $v>u$ on $\overline{\{u>0\}} \cap \partial B$.   By continuity $v - u>0$ in a neighborhood of $\overline{\{u>0\}} \cap \partial B$, since ${\limsup}^*\overline{\{u_k>0\}} \subset \overline{\{u>0\}}$ we get $v>u_k$ on $\overline{\{u_k>0\}} \cap \partial B$ for $k$ sufficiently large. Thus, by weak minimality of $u_k$, $v \geq u_k$ in $B$ and therefore, passing to the limit, $v \geq u$ in $B$. Thus $u$ is weakly minimal.

    Next consider the case $(u_k,E_k)$ are relaxed maximal subsolutions, \tref{stability} gives that $(u,E^*)$ is a relaxed subsolution. Let $B \subset\subset U$ and $v \in C(\overline{B})$ a subsolution of \eref{bernoulli-basic} in $B$ with $v<u$ on $\overline{\{v>0\}} \cap \partial B$.   By local uniform convergence and compactness of $\overline{\{v>0\}} \cap \partial B$ we have $v < u_k$ for sufficiently large $k$. Thus, by  maximality of $u_k$, $v \leq u_k$ in $B$ and therefore, passing to the limit, $v \leq u$ in $B$. Thus $(u,E^*)$ is a relaxed maximal subsolution.
\end{proof}

    \subsection{Sliding comparison} 
The standard comparison principle in bounded domains does not hold for solutions of \eqref{e.bernoulli-basic}. Pinning on heterogeneities is one source of such non-uniqueness, but non-uniqueness even occurs in the homogeneous problem \eref{bernoulli}. Instead, we will often use the following {sliding comparison principle}. We only state the supersolution version, the subsolution version is symmetric.  Note that $u$ does not need to be regular, only the sliding family.

\begin{lemma} [Sliding comparison \cite{CS}*{Theorem 2.2}]\label{l.sliding-comparison}
Let $u \in C(\bar{U})$ be a viscosity supersolution of \eqref{e.bernoulli-basic}.  Let $v_0 \in C(\bar{U})$ satisfy the following:
\begin{enumerate}[label = (\roman*)]
	\item $v_t(x):=v_0(x+te_d)$ are classical subsolutions of \eqref{e.bernoulli-basic} for all $t \in [0,T]$;
    \item $v_0\leq u$ in $U$;
    \item $v_t \leq u$ on $\partial U$ and $v_t <u$ on $\overline{\{v_t>0\}} \cap \partial U$ for $0\leq t\leq T$.
\end{enumerate}
    Then $v_t \leq u$ in $U$  for  $0\leq t\leq T$.
\end{lemma}

\section{A quantitative removable singularity result}\label{s.strong-max}

In this section we show that flatness improves in smaller scales for solutions of the exterior problem (\tref{defect-strong-max} below). Essentially this is a free boundary version of the Harnack inequality for the exterior problem. If the $u$ solves the Bernoulli problem without a defect, then this kind of Harnack inequality follows by combining several significant results from the literature \cite{DeSilva,kinderlehrer1977regularity,Lieberman}, see \cite{FKexterior}*{Section 2} for further discussion. The challenge we face is the presence of the defect in $B_1$, where we only assume upper and lower bounds on the slope at the free boundary.  A crucial ingredient for the proof is the strong maximum principle Proposition~\ref{p.strong-max-flat}, which states that flat super/subsolutions which touch a half-plane solution in $B_1$ are exactly half-plane solutions in a smaller ball. In particular it rules out the possibility of small-scale two-plane profiles. Proposition~\ref{p.strong-max-flat}, stable under uniform limits, is a new observation even for the homogeneous Bernoulli problem. 

\tref{defect-strong-max} will play a central role in Section \ref{s.pinned-solutions-2d} to rule out two-plane profiles appearing at intermediate scales for the single-site solutions that we consider.

In what follows, the slope conditions of $u$ on its free boundary are understood in the sense of viscosity sub- and supersolutions as discussed in Section~\ref{s.prelim-single-site}.

\begin{theorem}[Defect strong maximum principle]\label{t.defect-strong-max}
    Let $L \geq 1$ and let $u\in C(B_1)$, harmonic in its positive set, $\|\grad u\|_{L^\infty(B_1)} \leq L$ and $0 \in \partial\{u>0\}$.
    \begin{enumerate}[label = \textup{(\roman*)}]
        \item\label{part.main-super} For every $\delta > 0$ and $\ell \in (0,1]$ there is $r_0 = r_0(d,\delta,L,\ell) > 0$ with the following property. If
        \begin{equation}\label{e.flatness}
            (x_d)_+ \leq u \leq (x_d + \eta_0)_+ \ \hbox{ in } B_1, \qquad \eta_0 = \eta_0(d,\ell) \hbox{ from \pref{strong-max-flat}},
        \end{equation}
       and 
        \[\ell \leq |\grad u| \leq 1 \ \hbox{ on } \ \partial \{u>0\} \cap B_1 \setminus \overline{B_{r_0}},\] 
        then $u \leq (x_d + \delta)_+$ in $B_{1/2}$.
        \item\label{part.main-sub} For every $\delta > 0$ and $\ell \in (0,1)$ there is $r_0 = r_0(d,\delta,L,\ell) > 0$ with the following property. If
        \begin{equation}\label{e.nondegen-hyp}
            u \leq (x_d)_+ \ \hbox{ in } B_1 \qquad \hbox{and} \qquad |\grad u| \geq \ell \ \hbox{ on } \partial \{u>0\} \cap B_1,
        \end{equation}
        and  
        \[ |\grad u| \geq 1  \ \hbox{ on } \ \partial \{u>0\} \cap B_1 \setminus \overline{B_{r_0}},\]
        then $u \geq (x_d - \delta)_+$ in $B_{1/2}$. 
    \end{enumerate}
\end{theorem}
The proof of this result is via a compactness argument using the stability property \tref{stability} in the class of supersolutions and {relaxed subsolution}s.

\subsection{Removable singularities} We first show that non-degenerate isolated singularities of the Bernoulli problem are removable.

\begin{lemma}\label{l.removable-super}
    Let $u\in C(B_1)$ non-negative, Lipschitz, with $0 \in \partial\{u>0\}$. If $u$ is a supersolution of \eref{bernoulli} in $B_1 \setminus \{0\}$, then $u$ is a supersolution in $B_1$.
\end{lemma}

\begin{lemma}\label{l.removable-sub}
    Let $u\in C(B_1)$ non-negative, Lipschitz, with $0 \in \partial\{u>0\}$, and
    \begin{equation}\label{e.nondegen-origin}
        \limsup_{r \to 0^+}\ \frac1r \max_{\partial B_r} u > 0.
    \end{equation}
    If $(u,E)$ is a {relaxed subsolution} of \eref{bernoulli} in $B_1 \setminus \{0\}$, then $(u, E \cup \{0\})$ is a {relaxed subsolution} in $B_1$. In particular, if $u$ is a subsolution of \eref{bernoulli} in $B_1 \setminus \{0\}$ then $u$ is a subsolution in $B_1$. 
\end{lemma}

These lemmas are also a useful indicator that the effect of the dilute holes will vanish at zeroth order in the limit as $\delta \to 0$ for the periodic problem \eref{model}. This idea is later made precise in \lref{initial-flatness}.

\begin{proof}[Proof of \lref{removable-super}]
     We only need to check the supersolution property at $x=0\in\partial\{u>0\}$. Let $\varphi$ touch $u$ from below at $0$, then $\varphi(0)=0$. If $\grad\varphi(0) = 0$ we are done, so suppose $\grad\varphi(0) = \alpha e$, with $\alpha > 0$ and $|e|=1$. By \cite{CS}*{Lemma 11.17} there is $\alpha_1 \geq \alpha$ with $r^{-1}u(rx) \to \alpha_1(x\cdot e)$ in $\{x\cdot e > 0\}$. Since the blow-up sequence $r^{-1}u(rx)$ is uniformly Lipschitz and zero at the origin, it converges locally uniformly to some $w$ as $r\to 0$, along a subsequence. \tref{stability} then yields that $w$ is a supersolution of \eref{bernoulli} in $\R^d \setminus \{0\}$. Since $w(x)=\alpha_1(x\cdot e)$ in $\{x\cdot e\geq 0\}$, By the supersolution condition evaluated at any nonzero $x$ orthogonal to $e$  (see \rref{remove-strictness-issue}) we conclude that $\alpha_1 \leq 1$ and so $\alpha \leq 1$.
\end{proof}
\begin{proof}[Proof of \lref{removable-sub}]
      Let $E_0 := E \cup \{0\}$, which is relatively closed in $B_1$ and contains $\{u>0\}$. As before only touchings at $x=0$ need to be checked. Let $\varphi$ touch $u$ from above in $E_0$ at $0$ with $\Delta\varphi(0) < 0$. We must show $|\grad\varphi(0)| \geq 1$.  Since $\{u>0\} \subset E_0$, we have
    \begin{equation}\label{e.phi-plus}
        u \leq \varphi_+ \ \hbox{ in } B_{1}.
    \end{equation}
     Thus if $\grad\varphi(0) = 0$ then we have $\max_{\partial B_r} u \leq C_\varphi r^2$, contradicting \eref{nondegen-origin}, and it remains to consider the case $\grad\varphi(0) = \alpha e$ with $\alpha > 0$ and $|e|=1$.  Now, since $\overline{\{u>0\}}\subset E_0\subset \{\varphi\geq 0\}$ with $\varphi(0)=0$ and $0\in \partial\{u>0\}$, we conclude that both $E_0$ and $\overline{\{u>0\}}$ have an exterior tangent ball at $x=0$.

    As before, since the sequence $r^{-1}u(rx)$ is uniformly Lipschitz and vanishes at the origin, along a subsequence it converges locally uniformly. By \eqref{e.phi-plus} and \cite{CS}*{Lemma 11.17} the limit is $ w(x) := \alpha_1(x\cdot e)_+$ for some $\alpha_1\geq 0$. By \eref{nondegen-origin} we conclude $\alpha_1 > 0$. By \eqref{e.phi-plus} we conclude $\alpha \geq \alpha_1$. We will show that $\alpha_1\geq 1$.
    
    Let $E_0^*: = {\limsup}^*_{r \to 0} \, r^{-1}E_0$. By \tref{stability}, applied to $(r^{-1}u(rx), r^{-1}E_0)$, $(w, E^*_0)$ is a {relaxed subsolution} of \eref{bernoulli} in $\R^d \setminus \{0\}$. So $\{x \cdot e >0\} = \{w>0\} \subset E_0^*$ and since $E_0$ had a positive radius exterior tangent ball at $0$ we also conclude that $E_0^* \subset \{x \cdot e \geq 0\}$. Since $E_0^*$ is closed, it follows that $E_0^* = \{x \cdot e \geq 0\}$.
    
    Since $\overline{\{w>0\}} = E_0^*$ we conclude that $w$ is a (standard) subsolution in $\R^d \setminus \{0\}$. Testing this condition at any $x \cdot e =0$ with $x \neq 0$ (see \rref{remove-strictness-issue}) we conclude $\alpha_1 \geq 1$ and so also $\alpha \geq 1$ and we conclude.

The last statement on subsolutions follow from the fact that $u$ is a subsolution in $U$ if and only if $(u, \overline{\{u>0\}}\cap U)$ is a {relaxed subsolution}.
\end{proof}

\subsection{Strong maximum principles}

Next we show a strong maximum / minimum principle for the Bernoulli problem.  The difference between the upper and lower flatness conditions manifests already here.

\begin{lemma}[Strong maximum principle]\label{l.strong-max}
    Let $u\in C(B_1)$ non-negative and harmonic in $\{u>0\}$ with $0 \in \partial\{u>0\}$.
    \begin{enumerate}[label = \textup{(\roman*)}]
    \item\label{part.strongmax-sub} Suppose $(u,E)$ is a {relaxed subsolution} of \eref{bernoulli} in $B_1$ with $u \leq (x_d)_+$ in $B_1$ and $E \subset \{x_d \geq 0\}$. Then $u \equiv (x_d)_+$ in $B_1$.
        \item\label{part.strongmax-super} Suppose $u$ is a supersolution of \eref{bernoulli} in $B_1$ and $u \geq (x_d)_+$ in $B_1$. Then $u \equiv x_d$ in $\{x_d\geq0\} \cap B_1$.
    \end{enumerate}
\end{lemma}
Note that in \ref{part.strongmax-super} $u$ may not equal $0$ in $\{x_d<0\}$: counterexamples include the two-plane solution $u(x) = |x_d|$.

\begin{proof}
Let us show \ref{part.strongmax-sub}. Suppose $u \not \equiv (x_d)_+$ in $B_1$. Since $u(x)\leq (x_d)_+$, $u(x) < (x_d)_+$ for some $x=x_0$ with $0<r:=|x_0|<1$. Let $\psi$ be the harmonic function in $B_r^+:=B_r \cap \R^d_+$ with $\psi(x) = (x_d)_+-u(x)$ on $\partial (B_r^+)$. Note that since $u=0$ on $\{x_d=0\}\cap \partial B_r^+$,  $\psi$ is zero there as well. Hence by odd reflection through $x_d = 0$ we can extend $\psi$ to be smooth and harmonic in $B_r$.

Since $\psi(x_0) >0$ and $x_0\in\partial (B_r^+)$, harmonicity yields $\psi>0$ in $B_r^+$, and $\partial_{x_d}\psi(0)>0$ by Hopf's Lemma. Since $u$ is subharmonic in $B_r^+$ and $x_d - \psi$ is harmonic in $B_r^+$, by comparison $u \leq x_d - \psi$ in $\overline{B_r^+}$. Since $E\subset \{x_d\geq 0\}$, we conclude that $x_d - \psi(x)$ touches $u$ from above in $E^*$ at $0$. By maximum principle in $B_r^+$, $x_d - \psi>0$ in $B_r^+$. Since $0 \in \partial \{u>0\}$ and $u \leq x_d - \psi$ in $B_r^+$ we can't have $x_d - \psi \equiv 0$. So $x_d - \psi>0$ in $B_r^+$. Hopf lemma again at $0$ implies $\partial_{x_d}(x_d-\psi)(0)>0$. 

Since  $x_d-\psi(x)=0$ on $\{x_d=0\}\cap B_r$, we have $|\grad (x_d-\psi)|(0) =\partial_{x_d}(x_d - \psi)(0)\in(0,1)$, yielding a contradiction with the {relaxed subsolution} condition (see \rref{remove-strictness-issue}). 

  For part \ref{part.strongmax-super}, we have $u(x) \geq (x_d)_+$, and thus  $w(x) := u(x) - x_d \geq 0$ in $\{x_d>0\}$ and $w(0) = 0$. Since $w$ is continuous and harmonic in $\{x_d>0\}$, the argument of Hopf's Lemma implies that either $u(x) \equiv x_d$ in $x_d>0$ or
  \[u(x) - x_d \geq cx_d \ \hbox{ in a small neighborhood of } \ x = 0.\]
  In the latter case $(1+c)x_d$ touches $u$ from below at $0 \in \partial \{u>0\}$, contradicting the supersolution condition (see \rref{remove-strictness-issue}).  Thus we conclude. 
\end{proof}

Now we prove a key rigidity result that rules out small-scale two-plane profiles for almost-solutions of \eref{bernoulli} which are ``flat".  Note that we work in a more general setting than solutions of Bernoulli problem, in a relaxed class which is closed under uniform limits. In this class, this rigidity result is new even for the case $\ell=1$ and is of independent interest.
\begin{proposition}[Flat strong maximum principle]\label{p.strong-max-flat}

        For any $\ell \in (0,1]$, there is $\eta_0 = \eta_0(d,\ell) \in (0,\tfrac14]$ such that the following holds: suppose that
    \begin{enumerate}[label = \textup{(\alph*)}]
        \item $(u,E)$ is both a {relaxed subsolution} and a supersolution of
        \begin{equation}
            \begin{cases}
                \Delta u = 0 & \hbox{ in } \{u>0\} \cap B_1\\
        \ell \leq |\grad u| \leq 1 &\hbox{ on } \partial \{u>0\} \cap B_1;
            \end{cases}
        \end{equation}
        \item $0\in\partial\{u>0\}$ and, for some $\eta \in (0, \eta_0]$,
        \[ (x_d)_+ \leq u \leq (x_d + \eta)_+  \ \hbox{ and } \ \{x_d\geq 0\} \subset E \subset \{x_d+\eta\geq 0\} \hbox{ in } B_1.\]
    \end{enumerate}
    Then
    \begin{equation}\label{e.strong-max-flat-conclusion}
        u \equiv (x_d)_+ \ \hbox{ and } \ E  = \{x_d \geq 0\} \ \hbox{ in } \ \overline{B_{1/2}}.
    \end{equation}
    \end{proposition}
    \begin{proof}

        By \lref{strong-max}\ref{part.strongmax-super}, we already have $u \equiv x_d$ in  $\{x_d\geq0\} \cap B_1$. 
    It remains to show that 
    \begin{equation}\label{goal}
    E = \{x_d \geq 0\}\hbox{ in } B_{1/2}.
    \end{equation}
  We will proceed with a barrier argument.  Let $S :=\{-\eta \leq x_d \leq 0\}$ and $S_r:=S \cap \{|x'| \leq r\}$. Note that $u$ is subharmonic (\rref{subsolutions-subharmonic}), $u=0$ on  $\partial S \cap B_1$, and $u \leq \eta$ on  $S\cap B_1$. We will first show that $u$ decays exponentially with respect to $\eta$ in $S_{8/10}$. Intuitively this is clear due to the geometry of $O(\eta)$-thin strip domain as well as the doubling property of harmonic measures, but we carefully proceed with quantitative estimates.  
  
   Let us define \[ b(x) := \sqrt2\,\eta\,\frac{\cosh(\mu|x'|)}{\cosh(\tfrac{9}{10} \mu)}\,\cos\Bigl(\frac{\pi x_d}{4\eta}\Bigr), \qquad \mu := \frac{\pi}{4\eta\sqrt{d-1}}.\]
    When $d=2$, $b$ is harmonic and we claim that $b$ is superharmonic in $S$ for general $d$.  We compute, using $\sinh s \leq s\cosh s$ for the inequality,
    \begin{align*}
        \Delta_{x'}\cosh(\mu|x'|) &= \mu^2\cosh(\mu|x'|) + \frac{d-2}{|x'|}\,\mu\sinh(\mu|x'|) \\
        &\leq (d-1)\,\mu^2\cosh(\mu|x'|) = \Bigl(\frac{\pi}{4\eta}\Bigr)^2\cosh(\mu|x'|),
    \end{align*}
    On $S$ we have $\cos\Bigl(\frac{\pi x_d}{4\eta}\Bigr) \in [\cos\tfrac\pi4, 1]$ so, in particular, $b \geq 0$ and we can get
    \[\Delta b = \left[\frac{\Delta_{x'}\cosh(\mu|x'|)}{\cosh(\mu|x'|)}+\frac{\partial_{x_d}^2\cos\Bigl(\frac{\pi x_d}{4\eta}\Bigr)}{\cos\Bigl(\frac{\pi x_d}{4\eta}\Bigr)}\right] b \leq \left[\Bigl(\frac{\pi}{4\eta}\Bigr)^2-\Bigl(\frac{\pi}{4\eta}\Bigr)^2\right] b = 0.\] 
    Thus $b$ is superharmonic in $S$. 
    
    Since $b \geq \sqrt2\,\eta\cos\tfrac\pi4 = \eta \geq u$ on the lateral boundary $\partial S_{9/10} \cap \{|x'| = \frac{9}{10}\}$ and $u = 0$ on the top and bottom boundaries $\partial S\cap \{|x'| < \frac{9}{10}\}$, comparison principle yields that
    \[
    u \leq  b \leq \sqrt2\,\eta\,\frac{\cosh(\tfrac{8}{10}\mu)}{\cosh(\tfrac{9}{10} \mu)} \leq  C\eta\,e^{-c_d/\eta} =: \ep' \qquad \hbox{on } S_{8/10}.
    \]

    Next we refine the upper bound in $\{x_d\leq 0\}$. Let us denote $w$   as the harmonic replacement of $u$ in $S_{8/10}$. Then $u \leq w \leq \ep'$. Since $w=0$ on $x_d = 0$, the odd reflection of $w$ across $\{x_d=0\}$ is harmonic in $\{|x'| < \frac{8}{10}\} \times (-\eta,\eta)$ and bounded by $\ep'$. The interior gradient estimate for harmonic functions, applied to $O(\eta)$- neighborhoods of $\{x_d=0\}$, then yields 
    \[u \leq w \leq \tfrac{C\ep'}{\eta}\,|x_d| \ \hbox{ on } x \in S_{7/10} \ \hbox{ with }  \ x_d \geq -\eta/2 .\]
    A similar bound holds for $|x'| \leq \frac{7}{10}$ and $x_d \leq -\eta/2$ since $u \leq \ep' \leq \tfrac{2\ep'}{\eta}|x_d|$. Hence we conclude
    \begin{equation}\label{e.plane-decay}
        u \leq C e^{-c_d/\eta}\,|x_d| \ \hbox{ on } \ \{|x'| \leq 7/10\} \cap \{x_d \leq 0\}.
    \end{equation}

    Lastly we perform a version of sliding comparison.  For any $a\in B_{1/2}(0)\cap \{x_d=0\}$,  consider the sliding family $D_t(a) := \overline{B_1(a-(1+t)e_d)}$.  We aim to show that $D_t(a) \subset E^C$ for all $t > 0$. This will show that $E \cap B_{1/2} \subset \{x_d \geq 0\}$, completing the proof.
    
     For large $t$, $D_t$ lies outside of $E$. Let $t_*$ be the largest $t$ such that $D_{t} \cap E \neq \emptyset$ (so that $D_{t_*}$ is an exterior tangent ball of $E$ at each point of $D_{t_*} \cap E$). We claim that $t_* = 0$.
   Note that $D_t \cap \{|x'| \geq 6/10\} \subset \{x_d < -a_0\}$ for all $t \geq 0$, for some $a_0>0$.  Hence, since $\eta\leq \eta_0 \leq a_0$, if our claim is false and if $t_*>0$, then there is a point $z_0=(z',z_d) \in D_{t_*}$ that touches $E$ from outside with $z_d<0$ and $(z',z_d) \in S_{6/10}$. Note that in particular $z_0$ is an outer regular point of $\partial E$.

Now a rescaling of \lref{linear-growth-outer-reg} (with $z_0$ and $|z_d|$ replacing $0$ and $\rho$ respectively), as well as  \eref{plane-decay}, gives
    \[ c(d)\ell |z_d| \leq \sup_{B_{|z_d|/4}(z_0)}u \leq C(d) e^{-c(d)/\eta}|z_d|. \]
     Now choosing $\eta_0 := \min\{a_0, c(d)/\log (1+1/\ell)\}$, and recalling $\eta \leq \eta_0$, we obtain a contradiction.

    \end{proof}

    Now we are ready to prove \tref{defect-strong-max} via a compactness argument.
\begin{proof}[Proof of \tref{defect-strong-max}]
    Both parts are compactness arguments and begin the same way. If the statement fails for every $r_0$, then with $r_k \downarrow 0$ there are $u_k$ satisfying the hypotheses with $r_0 = r_k$ and points $x^k \in B_{1/2}$ where the conclusion fails. The $u_k$ are uniformly Lipschitz with $u_k(0) = 0$, so along a subsequence $u_k \to u_\infty$ locally uniformly in $B_1$, $E_\infty:= {\limsup}^*\{u_k>0\}$, and $x^k \to x^\infty \in \overline{B_{1/2}}$. The limit is $L$-Lipschitz, harmonic in its positivity set, and $u_\infty(0)=0$. 

    \emph{Part \ref{part.main-super}.} Passing to the limit in \eref{flatness}, and using that \eref{flatness} implies $\{x_d>0\} \subset \{u_k>0\} \subset \{x_d+\eta_0>0\}$ for every $k$, gives, in $B_1$,
    \[(x_d)_+ \leq u_\infty(x) \leq (x_d+\eta_0)_+ \ \hbox{ and } \ \{x_d \geq 0\} \subset E_\infty \subset \{x_d+\eta_0 \geq 0\}.\]
    In particular the first inequality gives $0 \in \partial\{u_\infty>0\}$. 
    
    By \tref{stability}, \lref{removable-super} and \lref{removable-sub} $(u_\infty,E_\infty)$ is, respectively, a {relaxed subsolution} and supersolution of
    \[\ell \leq |\grad u_\infty| \leq 1 \ \hbox{ on } \ \partial \{u_\infty>0\} \cap B_1.\]
    Thus \pref{strong-max-flat} applies:
    \begin{equation}\label{e.limit-rigidity}
        u_\infty \equiv (x_d)_+ \ \hbox{and} \ E_\infty = \{x_d \geq 0\}  \ \hbox{ in } \ \overline{B_{1/2}}.
    \end{equation}
    Now we aim for a contradiction at the $x^k$. Recall $u_k(x^k) > (x^k_d + \delta)_+$, so in the limit $u_\infty(x^\infty) \geq (x^\infty_d + \delta)_+$. If $x^\infty_d > -\delta$ this contradicts $u_\infty(x^\infty) = (x^\infty_d)_+$ directly. If $x^\infty_d \leq -\delta$, then for large $k$ the witnesses satisfy $x^k \in \{u_k>0\}$, so $x^\infty \in E_\infty$ by definition of the Kuratowski limit, with $x^\infty_d < 0$ this contradicts \eqref{e.limit-rigidity}.

    \emph{Part \ref{part.main-sub}.} Here $u_k(x^k) < (x^k_d - \delta)_+$ forces $x^k_d > \delta$. The hypotheses \eref{nondegen-hyp} pass to the limit and give $u_\infty \leq (x_d)_+$, $E_\infty \subset \{x_d \geq 0\}$, and the non-degeneracy property $\max_{\partial B_r} u_\infty \geq c(d)\ell r$ for all $r \in (0,1)$ (using also \lref{linear-growth-outer-reg}). In particular $u_\infty$ is non-degenerate at $0$ in the sense of \eref{nondegen-origin} and $0 \in \partial\{u_\infty>0\}$. By \tref{stability} and \lref{removable-sub} the pair $(u_\infty, E_\infty)$ is a {relaxed subsolution} of \eref{bernoulli} in $B_1$. Then \lref{strong-max}\ref{part.strongmax-sub} gives $u_\infty \equiv (x_d)_+$ in $B_1$. But the witnesses give 
    \[u_\infty(x^\infty) =\lim_{k \to \infty} u_k(x^k)  \leq x^\infty_d - \delta < (x^\infty_d)_+,\]
    a contradiction.
\end{proof}

\subsection{Blow-downs} Finally we apply the same circle of ideas to blow-downs of exterior {relaxed solution}s $(u,E)$ with a one-sided flatness condition, namely either
\begin{equation}\label{e.above-bound}
    u(x) \leq (x_d +t)_+ \ \hbox{ and } \ E \subset \{x_d + t \geq 0\} \  \hbox{ for some } \ t \in \R 
\end{equation}
or
\begin{equation}\label{e.below-bound}
    u(x) \geq (x_d +t)_+ \ \hbox{ for some } \ t \in \R.
\end{equation}

\begin{proposition}[Blow-downs of one-sided flat exterior solutions]\label{p.blow-down}
    Let $(u,E)$ be a {relaxed solution} of \eref{bernoulli} in $U := \R^d \setminus \overline{B_1}$, and suppose $u$ is $L$-Lipschitz on $\R^d$.
    \begin{enumerate}[label = \textup{(\roman*)}]
        \item\label{part.blow-down-1} If $(u,E)$ also satisfies \eref{above-bound} and $\max_{\partial B_\rho} u \geq \ell\rho$ for some $\ell>0$ and all large $\rho$, then  \[\lim_{r \to \infty} \frac{1}{r}u(rx) = (x_d)_+ \ \hbox{ locally uniformly on $\R^d$.}\]
        \item\label{part.blow-down-2} If $u$ also satisfies \eref{below-bound} then
        \[\lim_{r \to \infty} \frac{1}{r}u(rx)  = (x_d)_+ \ \hbox{ locally uniformly on $\{x_d \geq 0\}$.}\]
    \end{enumerate}
\end{proposition}

As usual the conclusion is weaker in the flat from below case, the blow-down may, for example, be a two-plane solution.

\begin{proof}
    Fix any sequence $r_k \to \infty$ and set $u_k(x) := \frac{1}{r_k}u(r_k x)$ and $E_k = \frac{1}{r_k}E$. The blow-down sequence is bounded and uniformly Lipschitz so a subsequence converges locally uniformly in $\R^d \setminus \{0\}$ to some $u_\infty$ with $u_\infty \leq L|x|$, which therefore extends continuously to $\R^d$ with $u_\infty(0) = 0$.  Call $E_\infty:= {\limsup}^*E_k$. By \tref{stability} $(u_\infty,E_\infty)$ is a supersolution and {relaxed subsolution} of \eref{bernoulli} in $\R^d \setminus \{0\}$.

    Part \ref{part.blow-down-1}. Taking the limit of \eref{above-bound} we find $u_\infty \leq (x_d)_+$ and $E_{\infty} \subset \{x_d \geq 0\}$. From the limit of the non-degeneracy hypothesis $\max_{\partial B_s} u_\infty \geq \ell s$ for every $s>0$; in particular $0 \in \partial \{u_\infty>0\}$. So \lref{removable-super} and \lref{removable-sub} imply that $(u_\infty,E_\infty)$ is a {relaxed solution} in $\R^d$. Then \lref{strong-max}\ref{part.strongmax-sub} implies that $u_\infty \equiv (x_d)_+$. Since the limit is the same on every subsequence we conclude convergence of the full blow-down sequence.

    Part \ref{part.blow-down-2}. Taking the limit \eref{below-bound} we find $u_\infty \geq (x_d)_+$, so $0 \in \partial\{u_\infty>0\}$ and $u_\infty$ is a supersolution of \eref{bernoulli} in $\R^d$ by \lref{removable-super}. So by \lref{strong-max}\ref{part.strongmax-super} $u_\infty \equiv x_d$ in $\{x_d\geq0\}$.  Since the limit is the same on every subsequence we conclude convergence of the full blow-down sequence.
\end{proof}

\section{\texorpdfstring{Pinned solutions on a single site in $d=2$}{Pinned solutions on a single site in two dimensions}}\label{s.pinned-solutions-2d}

In this section we will analyze the pinning force of a single-site defect in $d=2$ via a limit of finite domain approximations and we will prove \tref{main-2d-capacities}.

More specifically, for $R$ large to be sent to $\infty$, we consider the pinning problem in $B_R$
\begin{equation}\label{e.2d-k-R-setup}
    \begin{cases}
        \Delta u = 0 & \hbox{ in } \ \{u>0\} \cap B_R(0)\\
        |\grad u| = 1+q(x) &\hbox{on } \partial \{u>0\} \cap B_R(0)\\
        u(x) = (x_d+k\log R)_+ & \hbox{on } \partial B_R(0).
    \end{cases}
\end{equation}
Intuitively speaking we think of pulling the free boundary through the defect by slowly moving the boundary data via increasing (advancing case) or decreasing (receding case) the parameter $k$.  For a certain range of $k$ there will exist solutions that are pinned on the defect, i.e. their free boundary intersects the support of $q$, and then past a certain critical value the free boundary jumps and the solution is half-planar and does not see the defect. See \fref{2d-pinning-family} for a simulation.  

Our plan in this section is to first make a precise definition of this critical pinning force $\kappa^R_{\textup{adv}}$ and $\kappa^R_{\textup{rec}}$ in $B_R$, and then prove a limit theorem as $R \to \infty$ which includes a classification of the $R \to \infty$ limit in terms of the capacities of global {proper} solutions.

\begin{figure}
    \centering

\newcommand{\LabelW}{1em}
\newcommand{\ColSep}{4mm}
\newcommand{\RowSep}{4mm}

\newlength{\ImgW}
\setlength{\ImgW}{\dimexpr(\linewidth-\LabelW-2\ColSep)/2\relax}

\begin{tikzpicture}
\matrix (M) [
  nodes={inner sep=0pt, outer sep=0pt, anchor=center},
  column sep=\ColSep,
  row sep=\RowSep
]{
\node{}; &
\node {Positive bump}; &
\node {Negative bump}; &
\node[minimum width=\LabelW, align=left] {}; \\

\node[minimum width=\LabelW, align=left,rotate=90] {Advancing}; &
\node{\includegraphics[width=\ImgW]{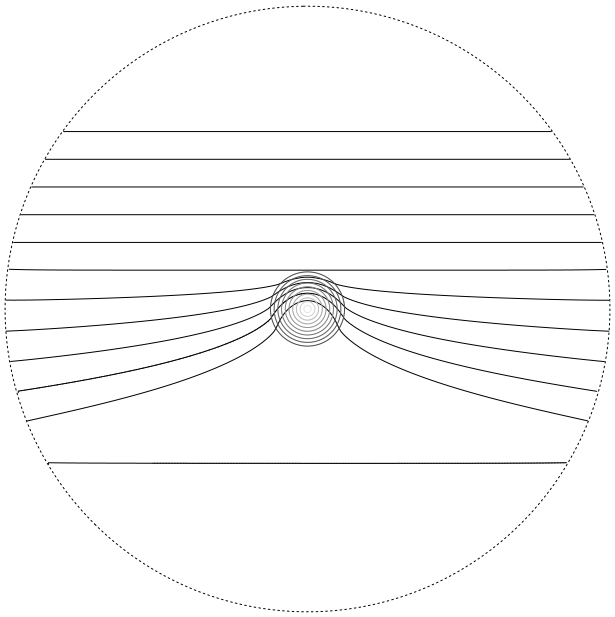}}; &
\node{\includegraphics[width=\ImgW]{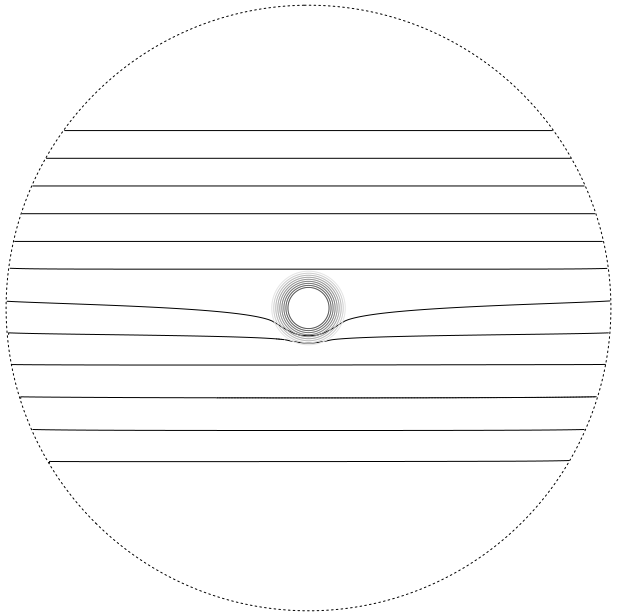}};  &
\node[minimum width=4mm, minimum height=\dimexpr\ImgW/2\relax] (a1) {};
\\

\node[minimum width=\LabelW, align=left,rotate=90] {Receding}; &
\node{\includegraphics[width=\ImgW]{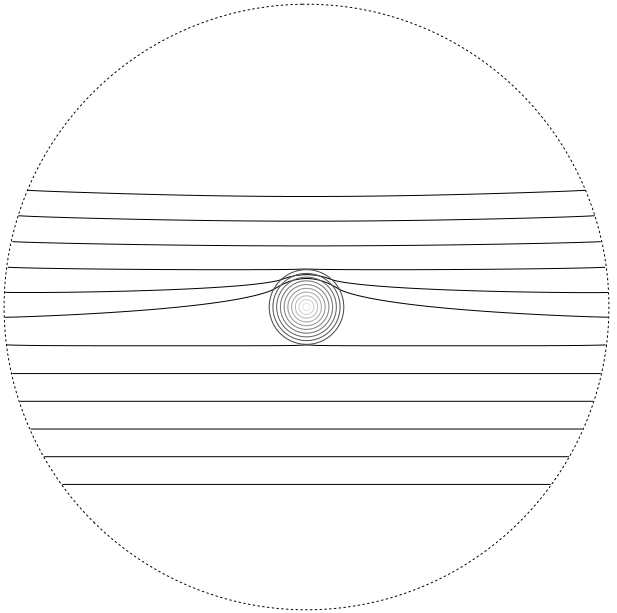}}; &
\node{\includegraphics[width=\ImgW]{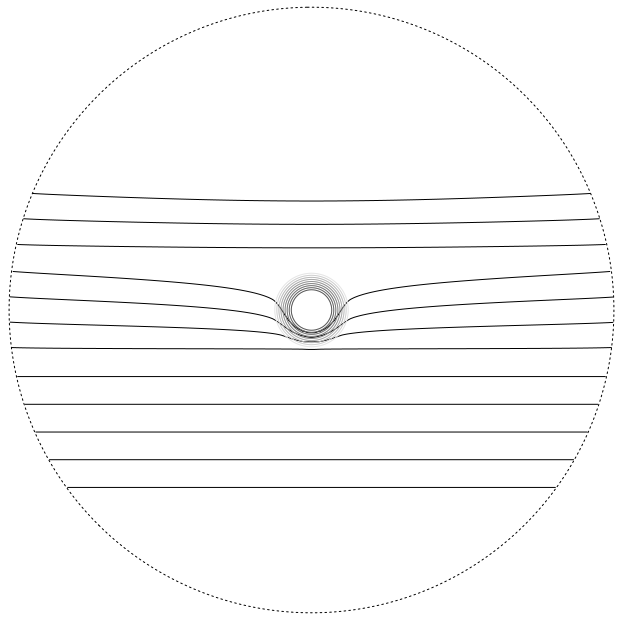}}; &
\node[minimum width=4mm, minimum height=\dimexpr\ImgW/2\relax] (a2) {}; \\
};
\draw[->] ($(a1.north)+(0,-1mm)$) -- ($(a1.south)+(0,1mm)$);
\draw[<-] ($(a2.north)+(0,-1mm)$) -- ($(a2.south)+(0,1mm)$);
\end{tikzpicture}
    \caption{ Free boundaries pinned on the defect in $\R^2$. Propagation directions are shown with arrows, with vertical direction $e_2$. Black curves are the free boundaries of the $u^{k,R}_{\textup{adv}}$ (top row) and $u^{k,R}_{\textup{rec}}$ (bottom row) for varying $k$, positive sets are above the free boundaries. Colored lines are contours of $q$: on the left column $q(x) = 3(1-|x|)_+$, on the right column $q(x) = -0.95[(2-|x|)_+ \wedge 1]$. } 
    \label{f.2d-pinning-family}
\end{figure}

\begin{definition}\label{d.adv-2d-R}
     Define $u^{k,R}_{\textup{adv}}$ to be the minimal supersolution of \eref{2d-k-R-setup} above $(x_d-1)_+$. Define
\begin{equation}\label{e.kappaR-adv}
    \kappa_{\textup{adv}}^R:= \sup \{k : \ \{u^{k,R}_{\textup{adv}}=0\} \cap \overline{B_1(0)} \neq \emptyset\} \geq 0.
\end{equation}
\end{definition}
\begin{definition}\label{d.rec-2d-R}
     Define $u^{k,R}_{\textup{rec}}$ to be the maximal subsolution of \eref{2d-k-R-setup} below $(x_d+1)_+$. Define
\begin{equation}\label{e.kappaR-rec}
    \kappa_{\textup{rec}}^R:= \inf \{k : \ \overline{\{u^{k,R}_{\textup{rec}}>0\}} \cap \overline{B_1(0)} \neq \emptyset\} \leq 0.
\end{equation}
\end{definition}

\begin{definition}\label{d.cap-2d}
    Define
\[k_{\textup{adv}}:= \sup \left\{k : \exists \hbox{ a {proper} solution $u$ of \eref{single-site} with capacity $k$}\right\} \]
and
\[k_{\textup{rec}} :=\inf \left\{k : \exists  \hbox{ a {proper} solution $u$ of \eref{single-site} with capacity $k$}\right\}.\]
\end{definition}
Naturally $k_{\textup{adv}} \geq 0 \geq k_{\textup{rec}}$ since half-plane solutions whose free boundary avoids the defect do exist.

Now we are ready to state our main convergence result.

\begin{theorem}\label{t.monotone-family2}
\begin{enumerate}[label=(\roman*)]
    \item\label{part.monotone-family2-p1} The limits hold
\[\lim_{R \to \infty} \kappa^R_{\textup{adv}} = k_{\textup{adv}} \ \hbox{ and } \ \lim_{R \to \infty} \kappa^R_{\textup{rec}} = k_{\textup{rec}}.\]
\item\label{part.monotone-family2-p2} For every $k \in (k_{\textup{rec}},k_{\textup{adv}})$ there exists a {proper} solution of \eref{single-site} with capacity $k$. 
\end{enumerate}
\end{theorem}
This result immediately implies \tref{main-2d-capacities} and also has some additional, and very useful information, about the behavior of the finite radius pinning problem \eref{2d-k-R-setup} as $R \to \infty$.

\begin{remark}
    Note that we have proved in \cite{FKexterior}*{Theorem 4.1} that $k_{\textup{rec}} > -\infty$. However, we only know that $k_{\textup{adv}} < + \infty$ under the assumption that $\max q$ is sufficiently small.  We do not know whether $k_{\textup{adv}} = +\infty$ is actually possible, but the limit theorem above holds regardless.
\end{remark}

\subsection{Exterior solutions with a prescribed capacity}We first prove the following Lemma, which we will use for a sliding barrier argument.

\begin{lemma}\label{l.exterior-barrier}
    There are $k_0>0$ and $C_0 \geq 1$ universal so that for all $|k| \leq k_0$ there exists $\phi_{k}$ with $\phi_{k} = (x_d)_+$ in $B_{1}$ and solving the homogeneous Bernoulli problem \eqref{e.bernoulli} in $U=\R^2\setminus B_{1}$, 
    with the property
\begin{equation}\label{e.sliding-psi-info}
    |\phi_k(x) -( x_d + k \log |x|)| \leq C_0 \ \hbox{ for } \ x \in \{\phi_k>0\} \cap U.
\end{equation}   \end{lemma}
\begin{remark}
Note that for $|k| \geq k_0$ if we define $R_1(k) := |k|/k_0 \geq 1$ then
\[\phi_k(x) := R_1\phi_{\textup{sgn}(k)k_0}(x/R_1)\]
satisfies $\phi_k = (x_d)_+$ in $B_{R_1}$, solves \eqref{e.bernoulli} in $\R^2 \setminus B_{R_1}$, and has the expansion
\[\phi_k(x) = x_d + k\log |x|+O(R_1\log(2+R_1))  \ \hbox{ for } \ x \in \{\phi_k>0\} \setminus B_{R_1}.\]
\end{remark}
\begin{proof}[Sketch]
The solutions can be constructed by solving the approximate problems
\[\begin{cases}
    \Delta \phi^R_k = 0 & \hbox{in } \{\phi^R_k>0\} \cap B_R \setminus B_1\\
    |\grad \phi^R_k| = 1 & \hbox{on } \partial \{\phi^R_k>0\} \cap B_R \setminus B_1\\
    \phi^R_k = (x_d)_+ & \hbox{on } \partial B_1\\
    \phi^R_k = (x_d + k \log |x|)_+ & \hbox{on } \partial B_R
\end{cases}\]
and sending $R \to \infty$. The barriers from \cite{FKexterior}*{Lemma 3.4} can be used via a comparison argument to guarantee that the limit is nontrivial and has the desired capacity $k$.  The argument is similar to what we do below in the proof of \tref{monotone-family2}.
\end{proof}

\subsection{Bounds on the finite radius pinned solutions}
Preliminary to the proof of \tref{monotone-family2} we establish several bounds on the solutions and extremal solutions of the finite radius pinning problem \eref{2d-k-R-setup}. We will need upper and lower bounds on the capacitory growth. The lower bounds use one-sided flatness and the barriers from \lref{exterior-barrier} along with a sliding argument. The upper bound is much more subtle, and it is especially essential for the advancing case where it is quite challenging to establish that solutions are {proper} from one-sided flatness from below due to the possibility of two-plane-like solutions.

First we record some (almost) trivial maximum principle upper and lower bounds on $u^{k,R}_{\textup{adv}}$ and $u^{k,R}_{\textup{rec}}$.  At the same time we show, closely related, one sided bounds on {proper} global solutions of \eref{single-site}.
\begin{lemma}[Maximum principle bounds]\label{l.ukR-trivial-bounds}
    Let $\log R>\max\{1,\frac{1}{|k|}\}$.
    \begin{enumerate}[label = (\alph*)]
        \item\label{part.ukR-trivial-bounds-p1} Assume that $\kappa^R_{\textup{adv}}>0$. For any $0 < k \leq \kappa^R_{\textup{adv}}$ there is $t_* \in [-1,1]$ so that
    \[(x_d+t_*)_+ \leq u^{k,R}_{\textup{adv}}(x) \leq (x_d+k \log R)_+ \ \hbox{ in } \ B_R\]
    and $(x_d+t_*)_+$ touches $u^{k,R}_{\textup{adv}}$ from below at a point in $\partial\{u^{k,R}_{\textup{adv}}>0\} \cap \overline{B_1}$.
    \item\label{part.ukR-trivial-bounds-p2} Assume that $\kappa^R_{\textup{rec}}<0$. For any $0 > k \geq  \kappa^R_{\textup{rec}}$ there is $t_* \in [-1,1]$ so that
    \[ (x_d+k \log R)_+ \leq u^{k,R}_{\textup{rec}}(x) \leq (x_d+t_*)_+ \ \hbox{ in } \ B_R\]
    and $(x_d+t_*)_+$ touches $u^{k,R}_{\textup{rec}}$ from above at a point in $\partial\{u^{k,R}_{\textup{rec}}>0\} \cap \overline{B_1}$.
    \item\label{part.ukR-trivial-bounds-p3} If $u$ is a {proper} solution of \eref{single-site} with capacity $k>0$ (resp. $k<0$) there is $t_* \in [-1,1]$ so that
    \[(x_d+t_*)_+ \leq u(x) \ \hbox{ in } \ \R^d, \ \hbox{(resp. $u(x) \leq (x_d+t_*)_+$)}\]
    and $(x_d+t_*)_+$ touches $u$ from below (resp. above) at a point in $\partial\{u>0\} \cap \overline{B_1}$. In particular $\partial\{u>0\} \cap \overline{B_1} \neq \emptyset$.
    \end{enumerate}
\end{lemma}
\begin{proof}
    The arguments in advancing and receding cases are almost symmetrical, so we just consider the advancing case. The upper bound of part \partref{ukR-trivial-bounds-p1} follows since $(x_d+k\log R)_+$ solves \eref{2d-k-R-setup} and $u^{k,R}_{\textup{adv}}$ is the Perron's method minimal supersolution of \eref{2d-k-R-setup}.

    The touching from below in parts \partref{ukR-trivial-bounds-p2} and \partref{ukR-trivial-bounds-p3} both follow from a very similar sliding comparison. We just present the global case \partref{ukR-trivial-bounds-p3}.  Slide $(x_d+t)_+$ from below until it touches $u$ from below:
    \[t_* := \sup\{t : (x_d+t)_+ \leq u(x)\}.\]
    Since $u$ is {proper} and $k>0$ the above supremum is on a nonempty set, and there is a point $x^* \in \partial \{u>0\}$ with $x^*_d + t_* = 0$. If $x^* \not \in \overline{B_1}$ then \lref{strong-max}\ref{part.strongmax-super} then $u(x) \equiv x_d + t_*$ in $x_d +t_* \geq 0$ which contradicts $u$ being {proper} with capacity $k>0$.  Thus $t_* \in [-1,1]$ and $x^* \in \overline{B_1}$.
\end{proof}
Next we use a sliding argument with the exterior capacitory barriers from \lref{exterior-barrier} to establish that $u^{k,R}_{\textup{adv}}$ and $u^{k,R}_{\textup{rec}}$ have the expected capacitory behavior from ``below", i.e. a lower bound on the growth away from $x_d$.
\begin{lemma}[Capacitory ``lower" bound]\label{l.capacitory-lower-bound}
    For $k\geq0$
    \begin{equation}\label{e.vk-adv-bound}
        u^{k,R}_{\textup{adv}}(x) \geq x_d+k \log |x|-C(k) \ \hbox{ in } \ \{u^{k,R}_{\textup{adv}}>0\},
    \end{equation}
and for $k\leq 0$
\begin{equation}\label{e.vk-rec-bound}
    u^{k,R}_{\textup{rec}}(x) \leq x_d+k \log |x|+C(k) \ \hbox{ in } \ \{u^{k,R}_{\textup{rec}}>0\}.
\end{equation}
\end{lemma}
\begin{proof}    
    
    The case $k=0$ is already implied by \lref{ukR-trivial-bounds}. We only consider the case $k> 0$ and prove \eqref{e.vk-adv-bound}: the argument to prove \eqref{e.vk-rec-bound}  is parallel.  
    
   Let $k>0$ and let $\phi_k$ be from \lref{exterior-barrier} and the remark following it, and set $R_0(k) := \max\{1,|k|/k_0\}$ so that $\phi_k=(x_d)_+$ in $B_{R_0}$ and \eref{sliding-psi-info} holds outside $B_{R_0}$ (up to adjusting $C_0$ depending on $R_0$). If $R \leq 2\max\{R_0,k\}$ then \eqref{e.vk-adv-bound} follows from $u^{k,R}_{\textup{adv}} \geq (x_d-1)_+$ after enlarging $C(k)$, so we assume $R \geq 2\max\{R_0,k\}$. 
   
   Let $\psi_t(x) := \phi_k(x-(1+t)e_d)$. We will perform a sliding comparison argument in the domain $B_R(0) \setminus B_{R_0}((1+t)e_d)$ to conclude  
\begin{equation}\label{e.sliding_k_R}
 \psi_t \leq u^{k,R}_{\textup{adv}} \ \hbox{ in } \ B_R(0) \ \hbox{ for } \ t = 2C_0 \quad \hbox{ for } R \geq 2\max\{R_0, k\}.
\end{equation}
    For all $t>0$ we have, due to \lref{ukR-trivial-bounds},
    \[u^{k,R}_{\textup{adv}}(x) \geq (x_d-1)_+ > (x_d-1-t)_+ = \psi_t(x) \ \hbox{ on } \  B_{R_0}((1+t)e_d) \cap \overline{\{\psi_t>0\}}.\] 
 Due to \eref{sliding-psi-info}, we have
    \[\psi_t(x) \leq x_d-(1+t)+k\log |x-(1+t)e_d| + C_0 \ \hbox{ for } \ x \in \{\psi_t>0\} \setminus B_{R_0}((1+t)e_d)\]
    For $x \in \partial B_R(0) \cap \{\psi_t>0\}$ and $t\geq 2C_0$, using that $R\geq 2k\geq 0$ 
    \begin{align*}
        \psi_t(x) &\leq x_d +k\log (R+1+t)+C_0-(1+t) \\
        & \leq x_d + k\log R + k\frac{1}{R}(1+t)+C_0-(1+t)\\
        & \leq x_d + k\log R + C_0-\frac{1}{2}(1+t)\leq x_d+k\log R = u^{k,R}_{\textup{adv}}(x). 
    \end{align*}
   Hence sliding comparison in $B_R(0) \setminus B_{R_0}((1+t)e_d)$ yields \eqref{e.sliding_k_R}.

    Now \eqref{e.vk-adv-bound} follows by bounds of $\psi_{2C_0}$.  For $|x| \geq 2(1+2C_0)$, \eref{sliding-psi-info} yields    \begin{align*}
        u^{k,R}_{\textup{adv}}(x) &\geq \psi_{2C_0}(x) \geq x_d+k\log |x-(1+2C_0)e_d|-1-3C_0\\
        &\geq x_d+k\log (|x|-(1+2C_0))-1-3C_0\\
        &\geq x_d + k\log |x| - [k \log 2 + 1 +3C_0]
    \end{align*}
   Lastly recall that for $1 \leq |x| \leq 2(1+2C_0)$
    \[u^{k,R}_{\textup{adv}}(x) \geq x_d- 1 \geq x_d + k \log |x| -[1 + k\log 2(1+2C_0)].\]
Hence we have shown \eref{vk-adv-bound}. 
\end{proof}

 Last we show that the finite radius pinned solutions can be extended to global supersolutions at the cost of a slight decrease in the capacitory growth. This is a crucial step which allows us to show that the $u^{k,R}_{\textup{adv}}(x)$ behave like a {proper} solution at all intermediate scales $1 \ll r \ll R$.

\begin{proposition}\label{p.ukR-supersolution-extension}
 Let $K >1$. For any $\delta>0$ there is ${R_0}(\delta,K) \gg 1$ such that for $R \geq {R_0}$ and $0 < k \leq \min\{\kappa^R_{\textup{adv}},K\}$ then
\begin{equation}\label{e.extended-barrier}
    w^{k,R}_{\textup{adv}}(x) :=  \begin{cases}
     u^{k,R}_{\textup{adv}}(x) & \hbox{in } B_R\\
     (x_d+(k-\delta)\log |x|+\delta\log R)_+ & \hbox{in } \R^2 \setminus B_R
 \end{cases}
\end{equation}
 is a supersolution of \eref{single-site} in $\R^2$.
\end{proposition}

\begin{proof}
    By \lref{ukR-trivial-bounds} there is $t_* \in [-1,1]$ such that
    \[x_d+t_* \leq u^{k,R}_{\textup{adv}}(x) \leq (x_d + k\log R)_+  \ \hbox{ in } \ B_R,\]
    and $x_d+t_*$ touches $u^{k,R}_{\textup{adv}}$ from below on its free boundary in $\overline{B_1}$. 

    Let $\eta_0$ from \tref{defect-strong-max} define $r := \eta_0^{-1}(k\log R+1)$. By \tref{defect-strong-max} part \ref{part.main-super}, applied to  $\tilde{u}(x,t):=r^{-1}u(r(x+t_*e_d))$, we find that for any $\delta>0$ there exists $R_0(\delta,K)$ sufficiently large so that for $R \geq R_0$
    \[u^{k,R}_{\textup{adv}}(x) \leq (x_d + \tfrac{\eta_0}{8K}\delta r)_+ \ \hbox{ in } \ B_{r/2}.\]
     Let
    \[\psi^s(x) := x_d+ (k-\delta)\log |x| + \delta \log R-s. \]
    By \cite{FKexterior}*{Lemma 4.4} $\psi^s(x)_+$ is a supersolution of \eref{bernoulli} in $\R^2 \setminus B_3$ for $ 0 \leq s \leq \delta \log R$. Note that 
    \[\psi^s(x)_+ = (x_d + k\log R-s)_+ \leq u^{k,R}_{\textup{adv}}(x) \ \hbox{ on } \ \partial B_R, \hbox{ and } \]
    
    \[\psi^s(x) = x_d + (k-\delta)\log r + \delta \log R-s \quad \hbox{ on } \partial B_r.\]
    Choose $R_0=R_0(K)$ so that $\frac{\log r}{\log R} \leq \frac{1}{2}$ and $\eta_0 r \geq 2$. Then
    \begin{align*}
        (k-\delta)\log r + \delta \log R & \geq  (1-\frac{\log r}{\log R})\delta \log R \\
        &\geq \frac{1}{2}\delta \log R= \frac{1}{2k} \delta (\eta_0r-1)\geq \frac{\eta_0}{4K} \delta r.
    \end{align*}
    Hence for $0<s < \frac{\eta_0}{4K} \delta r$ we have  $ u^{k,R}_{\textup{adv}}(x) \prec \psi^s(x)_+$ on $\partial B_r$, and thus 

    \[w_s(x):= \begin{cases}
        u^{k,R}_{\textup{adv}}(x) & x \in B_r \\
        \min \{\psi^s(x)_+, u^{k,R}_{\textup{adv}}(x)\} & x \in B_R \setminus B_r \\
        \psi^s(x)_+ & x \in \R^2 \setminus B_R
    \end{cases}\]
    is a supersolution of \eref{single-site} in $\R^2$.  By sending $s \to 0$ and using uniform stability (\tref{stability}  \ref{part.stability-super}), $w_0(x)$ is also a supersolution of \eref{single-site} in $\R^2$.
    
    Since $\psi^0_+ \geq u^{k,R}_{\textup{adv}}$ on $\partial (B_R \setminus B_r)$,  and since $u^{k,R}_{\textup{adv}}$ is a minimal supersolution of \eref{bernoulli} in $B_R \setminus B_r$, we conclude that 
     \[u^{k,R}_{\textup{adv}}(x) \leq \psi^0(x)_+ \ \hbox{ in } \ B_{R} \setminus B_r.\]
     Thus $w^{k,R}_{\textup{adv}}(x) = w_0(x)$ in $\R^2$ and so $w^{k,R}_{\textup{adv}}$ is a supersolution of \eref{single-site} in $\R^2$.
\end{proof}

An analogous result holds in the receding case as well. The explicit supersolution $(x_d+k\log |x|)_+$ does not work as a subsolution, and we use $\phi_k$ instead.  
\begin{proposition}
    \label{p.ukR-subsolution-extension}
 Let $K >1$. For any $\delta>0$ there are ${R_0}(\delta,K) \gg 1$ and $C(K)>0$ such that for $R \geq {R_0}$ and $0 > k \geq \max\{\kappa^R_{\textup{rec}},-K\}$ then, with $C = C(K)$,
\begin{equation}\label{e.extended-barrier-phi}
    w^{k,R}_{\textup{rec}}(x) :=  \begin{cases}
    u^{k,R}_{\textup{rec}}(x) & x \in B_{R/2}\\
        \max\{u^{k,R}_{\textup{rec}}(x),\phi_{k+\delta}(x-(\delta \log R+C)e_d)\} & x \in B_R \setminus B_{R/2}\\
        \phi_{k+\delta}(x-(\delta \log R+C)e_d) & x \in \R^2 \setminus B_R
    \end{cases}
\end{equation}
 is a continuous subsolution of \eref{single-site} in $\R^2$.
\end{proposition}
The proof relies on \tref{defect-strong-max} part \ref{part.main-sub} in an analogous way to the application of \tref{defect-strong-max} part \ref{part.main-super} in the proof of \pref{ukR-supersolution-extension}.

\subsection{Limit of the finite radius pinning thresholds}
Now we have the necessary set-up to return to the proof of \tref{monotone-family2}. Both the advancing and receding case have their own important challenges. In the advancing case the solutions are only flat from below, which potentially allows two-plane behavior to appear at intermediate scales. In the receding case we need to construct a monotone increasing sequence of finite $R$ solutions in order to avoid degeneracy issues in the limit. In each case \pref{ukR-supersolution-extension} or \pref{ukR-subsolution-extension} are used, respectively, to rule out these behaviors. Recall that the defect strong maximum principle \tref{defect-strong-max} plays a crucial role in the proof of those Propositions.  Since we regard the issue of ruling out two-plane solutions as more fundamental we focus on the advancing case, but we also carefully present the key details of the monotonicity argument in the receding case.
\begin{proof}[Proof of \tref{monotone-family2}]

We will establish the equation for the limit by proving an upper bound for the limit supremum, and a lower bound for the limit infimum.  We will focus on the advancing case which has a major additional challenge since it is trickier to establish that global solutions are {proper} in the flat from below case.  The differences in the receding case will be pointed out at the end of each part of the proof.

{\bf Step 1.} ($\liminf_{R \to \infty} \kappa^R_{\textup{adv}} \geq k_{\textup{adv}}$) The case $k_{\textup{adv}} = 0$ is immediate so we can assume that $k_{\textup{adv}}>0$. It will suffice to show that:
 \begin{equation}\label{e.zero_set}
 \begin{array}{c}
 \hbox{ for any $0<k<k_{\textup{adv}}$ there is $R_1(k)$ large so that for $R \geq R_1$, } \\
 \{u^{k,R}_{\textup{adv}}=0\} \cap \overline{B_1}\neq \emptyset \ \hbox{ and therefore } \ \kappa^R_{\textup{adv}} \geq k.
 \end{array}
 \end{equation}
      Notice that \eref{zero_set} implies that $\liminf_{R \to \infty} \kappa^R_{\textup{adv}} \geq k$ for all $0<k < k_{\textup{adv}}$ and so $\liminf_{R \to \infty} \kappa^R_{\textup{adv}} \geq k_{\textup{adv}}$.

     Now we proceed to establish the claim \eref{zero_set}. Let $0<k<k_{\textup{adv}}$. By \dref{cap-2d} there is $k<k_0\leq k_{\textup{adv}}$,  and a {proper} solution $u$ of \eref{single-site} with capacity $k_0$. By \lref{ukR-trivial-bounds} part \partref{ukR-trivial-bounds-p3} we know the free boundary of $u$ must touch the support of $q$, i.e. $\{u=0\} \cap \overline{B_1} \neq \emptyset$.

Now we compare $u$ with $u^{k,R}_{\textup{adv}}$. Since $u(x) \geq (x_d-1)_+$,  $\psi(x):=\min\{u^{k,R}_{\textup{adv}}(x),u(x)\}$ is above $(x_d-1)_+$. Moreover, since $k<k_0$ we can choose $R_1$ large such that for $R \geq R_1$
    \[u(x) \geq (x_d+ k\log R)_+ \ \hbox{ on }  \partial B_R,\]
and thus $\psi$ is a supersolution of \eqref{e.2d-k-R-setup}.  Since $u^{k,R}_{\textup{adv}}(x)$ is the minimal supersolution of \eqref{e.2d-k-R-setup} above $(x_d-1)_+$ so we conclude that $u^{k,R}_{\textup{adv}} \leq \psi \leq u$. But this means that 
    \[\{u^{k,R}_{\textup{adv}}=0\} \cap \overline{B_1}\supset \{u=0\} \cap \overline{B_1} \neq \emptyset.\]
    Thus we proved \eref{zero_set} and therefore established $\liminf_{R \to \infty} \kappa^R_{\textup{adv}} \geq k_{\textup{adv}}$.

    For the receding case the proof of the inequality $\limsup_{R \to \infty} \kappa^R_{\textup{rec}} \leq k_{\textup{rec}}$ is closely analogous to the above.

    {\bf Step 2.} ($ \limsup_{R \to \infty} \kappa^R_{\textup{adv}} \leq k_{\textup{adv}}$ and existence of solutions with prescribed capacity)  Define $ \kappa_\infty:=\limsup_{R \to \infty} \kappa^R_{\textup{adv}}$.  It will suffice to show that for every $0<k < \kappa_\infty$ there is a {proper} solution of \eref{single-site} with capacity $k$. We allow the possibility that $\kappa_\infty$ and/or $k_{\textup{adv}}$ are equal to $+\infty$.
    
     Let $0<k < \kappa_\infty$. Let  $R_j \to \infty$ along which $\kappa_{j}:=\kappa^{R_j}_{\textup{adv}} \to \kappa_\infty$. Up to a further subsequence we can assume that $0<k  < \kappa_j$ for all $j$.  Then define $u_j := u^{k,R_j}_{\textup{adv}}$. We will show that, up to taking a further subsequence, the $u_j$ converge to a global {proper} solution $u_\infty$ with capacity $k$.
    
    For compactness, note that, since $0 < k < \kappa^{R_j}_{\textup{adv}}$  there is at least one point $x_j \in \{u_j=0\} \cap \overline{B_1}$. So the Lipschitz estimate \lref{lipschitz-estimate} implies that $u_j(x) \leq C|x-x_j|$ in $B_{R_j}$. By uniform Lipschitz continuity, \lref{lipschitz-estimate} again, along a subsequence (not relabeled) the $u_j$ converge locally uniformly in $\R^d$ to some $u_\infty$ and $\{u_\infty = 0\} \cap \overline{B_1} \neq \emptyset$. Also recall that $u_j(x) \geq (x_d-1)_+$ so $u_\infty(x) \geq (x_d - 1)_+$ and so $u_\infty$ is one-sided flat and nontrivial.  By \lref{extremal-stability} $u_\infty$ is a solution of \eref{single-site}.
    
    By \lref{capacitory-lower-bound}, 
    \[u_j(x) \geq (x_d +k\log |x| - C)_+ \ \hbox{ in } \ B_R.\]
     Taking the limit of the lower bound
  \begin{equation}\label{e.ukR-limit-LB}
  u_\infty(x) \geq (x_d +k\log |x| - C)_+ \ \hbox{ in } \ \R^d.
  \end{equation}
    
    Next we aim to establish a similar upper bound, in the advancing case -- as we currently consider -- we need some kind of upper bound to obtain that $u_\infty$ is {proper}. But, moreover, we want to show that $u_\infty$ has capacity exactly $k$ so we need a matching upper bound in any case.  This is where we make use of \pref{ukR-supersolution-extension}. Let $\delta := \min\{\frac{1}{2}(\kappa_\infty - k),1\}$. Then the \pref{ukR-supersolution-extension} and the definition of $\kappa_\infty$ as the limit supremum imply that there is some $R_* > 1$ sufficiently large, depending on $\delta>0$, $k$, and on the convergence rate to the limit supremum, so that $k+\delta < \kappa^{R_*}_{\textup{adv}}$ and
    \[w (x) := \begin{cases}
    u^{k+\delta,R_*}_{\textup{adv}}(x) & \hbox{in } B_{R_*} \\
    (x_d + k\log |x| + \delta \log R_*)_+ & \hbox{in } \R^d \setminus B_{R_*}
     \end{cases}\]
     is a supersolution of \eref{single-site} in $\R^d$.  Then for any $j$ large enough that $R_j > R_*$ we have that $w$ is a supersolution of \eref{2d-k-R-setup} in $B_{R_j}$ and, since $u_j$ is the minimal supersolution of that problem, we conclude that
     \[ u_j(x) \leq w(x) \ \hbox{ in } \ B_{R_j} \ \hbox{ for all $j$ large enough that $R_j > R_*$.}\]
     Taking the limit we find
     \[u_\infty(x) \leq w(x) \ \hbox{ in } \ \R^d.\]
     Combined with the lower bound \eref{ukR-limit-LB} this implies that $u_\infty$ is {proper} and has capacity exactly $k$.
     
     Thus we have proved that $\limsup_{R \to \infty} \kappa^R_{\textup{adv}} \leq k_{\textup{adv}}$ and also that, for every $k \in (0,k_{\textup{adv}})$, there exists a {proper} solution of \eref{single-site} with capacity $k$, namely $u_\infty$ constructed above (note that the existence for $k=0$ is trivial e.g. $(x_d+2)_+$ works).

In the receding case we follow an analogous strategy to prove $\kappa_\infty:=\liminf_{R \to \infty} \kappa^R_{\textup{rec}} \geq k_{\textup{rec}}$, but more care is required to establish that $u_\infty$ is a solution due to the weaker stability property for maximal subsolutions. Instead of using \lref{nondegen-stability},
we construct a monotone increasing sequence $u^{k_j,R_j}_{\textup{rec}}$ so that we can apply \cref{monotone-stability} instead.  For $\kappa_{\infty} < k < 0$, let $k_j := (1-\sum_{n=j}^\infty 2^{-n})k$ and  $\delta_j := 2^{-j}k$. Note that $k_j+\delta_j = k_{j+1}$. 

Given a good choice of $R_j$ we explain how to make a good choice $R_{j+1}$. Assume that we have chosen $R_j \geq R_0(\delta_{j})$ so that \pref{ukR-subsolution-extension} applies and $w^{k_j,R_j}_{\textup{rec}}$, as defined in \eref{extended-barrier-phi} with $\delta = \delta_{j}$, is a subsolution of \eref{single-site}. By the definition of $w^{k_j,R_j}_{\textup{rec}}$ and \eref{sliding-psi-info}
\[ w^{k_j,R_j}_{\textup{rec}}(x) \leq  (x_d+(k_j+\delta_j)\log |x|)_+ \ \hbox{ for } \ |x| \geq R_j.\]
Thus the maximal subsolution property of $u^{k_{j+1},R}_{\textup{rec}}$ gives that $u^{k_{j+1},R}_{\textup{rec}}\geq w^{k_j,R_j}_{\textup{rec}}$ in $B_R(0)$ for any $R \geq R_j$ and hence $u^{k_{j+1},R}_{\textup{rec}} \geq u^{k_{j},R_j}_{\textup{rec}}$ in $B_{R_j}(0)$.  Now choose $R_{j+1} \geq R_j$ sufficiently large so that $R_{j+1} \geq R_0(\delta_{j+1})$ and also so that $\kappa^{R_{j+1}}_{\textup{rec}} < k < k_{j+1}$.  With this monotone choice of sequence the rest of the proof is parallel to the advancing case already handled in detail.

\end{proof}

\section{\texorpdfstring{Pinned solutions on a single site in $d\geq3$}{Pinned solutions on a single site in higher dimensions}}\label{s.monotone-family}

In this section we will analyze the solution pinned on a single site in $d \geq 3$ and prove \tref{main-3d-capacities}. Namely we will show that there are a pair of laminations of extremal solutions to \eref{single-site}. These two families bound all the {proper} solutions of \eref{single-site} and, up to closure, contain all the solutions with extremal values of the capacity. 

For $d=2$, since the capacity appears at higher order in the asymptotic expansion, the solution with the extremal pinning force is the last one visible as free boundary is pulled past the defect. For $d\geq 3$ this is not the case and
and different pinning forces may be viewed at different heights. 
This leads us to use a different approximation procedure to construct the global pinned solutions, based on solving half-space problems indexed by the height variable.

 \begin{theorem}\label{t.monotone-family}
Let $d \geq 3$. There exist families of  maximal subsolutions $(u_{\textup{rec}}(x;s))_{s\in\R}$ and minimal supersolutions $(u_{\textup{adv}}(x;s))_{s\in\R}$ of \eref{single-site}  such that the following holds: 
    \begin{enumerate}[label = (\roman*)]
\item\label{part.minimal-super-family}
     Both $u_{\textup{adv}}$ and $u_{\textup{rec}}$ are non-decreasing with respect to $s$.
     Furthermore    
        \begin{equation}\label{expansion} u_{\textup{adv}/\textup{rec}}(x;s) = x_d +  s - \varkappa_{\textup{adv}/\textup{rec}}(s)\frac{1}{|x|^{d-2}}+O(\frac{1}{|x|^{d-1}}) \ \hbox{ as } \ \{u>0\}\ni x \to \infty.
        \end{equation}

        \item\label{part.ordering} Let $u$ be a {proper} solution of \eref{single-site} with height $s$ and capacity $k$.  Then 
        \[u_{\textup{adv}}(x;s-) \leq u(x) \leq u_{\textup{rec}}(x;s+) \ \hbox{ and } \ \varkappa_{\textup{rec}}(s+) \leq k \leq \varkappa_{\textup{adv}}(s-).\]
        If $u$ is a minimal supersolution, then 
        \[ u(x) \leq u_{\textup{adv}}(x;s+) \ \hbox{ and } \ k \geq \varkappa_{\textup{adv}}(s+).\]
        If $u$ is a maximal subsolution, then
        \[u_{\textup{rec}}(x;s-) \leq u(x) \ \hbox{ and } \ k \leq \varkappa_{\textup{rec}}(s-) .\]
        \end{enumerate}
        
        \end{theorem}

The following theorem studies the range of capacities as $s$ varies.

\begin{theorem}\label{t.capacity-properties} 
The capacity functions satisfy the following properties:
        \begin{enumerate}[label =(\roman*)]
        \item\label{part.capacity-properties-0} The following capacity functions are well-defined
\begin{equation}\label{e.kappa-relation}
    \kappa_{\textup{rec}}(s): = \varkappa_{\textup{rec}}(s+) \ \hbox{ and } \ \kappa_{\textup{adv}}(s) := \varkappa_{\textup{adv}}(s-),
\end{equation}
in particular the limits implicit on the right hand sides exist.
            \item\label{part.capacity-properties-1}  $\kappa_{\textup{adv}/\textup{rec}}$ are compactly supported.
            \item\label{part.capacity-properties-2}  $\kappa_{\textup{adv}}$ is left continuous and USC, and $\kappa_{\textup{rec}}$ is right continuous and LSC.
            \item\label{part.capacity-properties-3} $\kappa_{\textup{adv}/\textup{rec}}$ have at most countably many jump discontinuities. 
            \item\label{part.capacity-properties-4}
            For every $s \in \R$ \begin{equation}\label{e.kappa-3d-def}
    \begin{split}
        \kappa_{\textup{adv}}(s)=\sup\{k: \exists \hbox{ a {proper} solution of \eqref{e.single-site} with capacity $k$ and height $s$}\}\\
        \kappa_{\textup{rec}}(s)\leq\inf\{k: \exists \hbox{ a {proper} solution of \eqref{e.single-site} with capacity $k$ and height $s$}\},
    \end{split}
\end{equation}
equality also holds in the latter at all but countably many $s$.

        \end{enumerate}
       
\end{theorem}

As usual we are just considering the {proper} solutions in the normalized setting with  direction $e_d$. The continuity properties of $\kappa_{\textup{adv}/\textup{rec}}(s;q)$ with respect to $q$ can also be studied in this framework, we skip the analysis there since it is not used for the main results of the paper.

\subsection{Half-space approximation problem}

In order to construct the laminations of global solutions in \tref{monotone-family}, we need an approximation procedure by some boundary value problems in subdomains of $\R^d$.

Fix $s \in \R$, let $R > -s$ and consider the half-space problem
\begin{equation}\label{e.finite-strip-1}
    \begin{cases}
    \Delta w = 0 & \hbox{ in } \{w>0\} \cap \{x_d < R\},\\
    |\grad w| = 1+q(x) & \hbox{ on } \partial \{w>0\} \cap \{x_d < R\},\\
    w = x_d+s  & \hbox{ on }  \{x_d \geq R\},
\end{cases}
\end{equation}
with the additional boundedness requirement 
\begin{equation}\label{e.fb-upper-bound}
    \sup_{\{w>0\}} |w(x) - (x_d+s)| < + \infty.
\end{equation}
In this section we show that $w$ converges to $(x_d+s)_+$ as $|x'| \to \infty$ uniformly in $R$. This is preparation for taking the limit $R \to \infty$ in the next section.

We begin by pointing out a simple a-priori maximum principle bound. For any $w$ solving \eref{finite-strip-1}-\eref{fb-upper-bound}, one can verify by sliding comparison \lref{sliding-comparison}, see e.g.  the proof of \lref{ukR-trivial-bounds}, that 
\begin{equation}\label{e.strip-solution-bounds}
 (x_d+s \wedge (-1))_+ \leq w(x) \leq (x_d + s \vee 1)_+.
\end{equation}

The proof of this quantitative asymptotics will be split into two lemmas. First, in \lref{transverse-limit-strip}, we prove the qualitative version.

\begin{lemma}\label{l.transverse-limit-strip}
    If $w$ solves \eref{finite-strip-1}- \eref{fb-upper-bound}, then for $d\geq 2$ we have
    \[\lim_{r \to \infty} \sup_{\{w>0\} \setminus B_r(0)} |w(x) - (x_d+s)| = 0.\]
\end{lemma}

\begin{proof}
   
         Take $s=0$, the general $s \in \R$ case is the same. We will show that 
         \[  \liminf_{\{w>0\}\ni x \to \infty} (w(x) - x_d) = 0,\,\, \limsup_{\{w>0\} \ni x \to \infty} (w(x)-x_d) = 0. \quad \]
Note that          both limits are finite since $w(x) - x_d$ is bounded.
Suppose that
    \[\liminf_{\{w>0\} \ni |x| \to \infty} (w(x)-x_d) = - \ep < 0.\]
    Then there is a sequence $x_k\in\{w>0\} \to \infty$ so that $w(x_k) -(x_k)_d \to - \ep$. Let $x_k'$ be the projection of $x_k$ onto $x \cdot e_d = 0$. Taking a further subsequence we can assume that $(x_k)_d \to a \in [\ep,R]$, since $w \geq 0$ and so
    \[a := \lim_{k \to \infty} (x_k)_d = \lim_{k \to \infty} (w(x_k) + \ep) \geq \ep.\]  Define
    \[w_k(x) := w(x_k'+x).\]
    Since the $w_k$ are uniformly Lipschitz, along a subsequence $w_k$ converge locally uniformly to a Lipschitz continuous $w_\infty$. By \tref{stability}, $w_\infty$ solves 
 \begin{equation}\label{e.finite-strip-2}
    \begin{cases}
    \Delta w_\infty = 0 & \hbox{ in } \{w_\infty>0\} \cap \{x_d < R\};\\
    |\grad w_\infty| \leq 1 & \hbox{ on } \partial \{w_\infty>0\} \cap \{x_d < R\};\\
w_\infty(x) = x_d  & \hbox{ on }  \{x_d \geq R\}.
\end{cases}
\end{equation}
Furthermore, 
    \[w_\infty(0,a) - a +\ep = 0 \ \hbox{ and } \ w_\infty(x) - x_d +\ep\geq  0\ \hbox{ in } \ \{w_\infty>0\}.\]
   Thus $x_d-\e$ touches $w_{\infty}$ from below at $(0,a)$. This leads to a contradiction due to strong maximum principle and \lref{strong-max}\ref{part.strongmax-super}.

    The $\limsup$  side of the argument is parallel, using \tref{stability}\ref{part.stability-pair} and \lref{strong-max}\ref{part.strongmax-sub}.
   
\end{proof}

Next we show a stronger, quantitative estimate for \lref{transverse-limit-strip}.

\begin{lemma}\label{l.finite-strip-uniform-bdry-layer}
   
     Let $d \geq 3$. There exists $C(\min q, \max q)\geq 1$ such that, if $w$ solves \eref{finite-strip-1} - \eref{fb-upper-bound} with $R+s>0$, then 
    \[ |w(x) - (x_d+s)| \leq C(1+|x+se_d|)^{2-d} \ \hbox{ in } \ \{w>0\}.\]
\end{lemma}

\begin{proof}
By \cite{FKexterior}*{Proposition 4.2} there exists $\phi$ a subsolution of 
\[\Delta \phi = 0 \ \hbox{in } \{\phi>0\} \ \hbox{ and } \
         |\grad \phi| = 1+ (\max q){\bf 1}_{B_1' \times \R} \ \hbox{on } \ \partial \{\phi>0\},\]
         with
\[\phi(x) \geq x_d -C(d,\max q) (1+|x|)^{2-d} \ \hbox{ in } \ \{\phi >0\}.\]
By  \lref{sliding-comparison}, we have  $\phi(x+te_d) \leq w(x)$ for $t \leq s$. Therefore
\[w(x) \geq \phi(x+se_d) \geq (x_d+s) - C(1+|x+se_d|)^{2-d}.\]
The upper bound can be shown similarly, using the supersolution barrier constructed in \cite{FKexterior}*{Proposition 4.2}.

\end{proof}

\subsection{Construction of the lamination in the whole space}

Now we consider the minimal supersolutions $u^R_{\textup{adv}}(x;s)$ and maximal subsolutions $u^R_{\textup{rec}}(x;s)$ of \eref{finite-strip-1}.
More precisely
\[u^R_{\textup{adv}}(x;s):=\inf_{w \in \mathcal{S}^R_{\textup{super}}(s)} w(x) \ \hbox{ and } \ u^R_{\textup{rec}}(x;s):=\sup_{w \in \mathcal{S}^R_{\textup{sub}}(s)} w(x)\]
where
\[\mathcal{S}^R_{\textup{super/sub}}(s):=\left\{w \in C(\R^d): \begin{array}{cc}\hbox{$w$ satisfies \eref{fb-upper-bound} and}\\\hbox{ is a super/subsolution of \eref{finite-strip-1}}\end{array}\right\}.\]
(see \sref{perron-method} for Perron's method.)

\begin{lemma}\label{l.monotone}
  For fixed $x,R$ the maps $s \mapsto u^R_{\textup{rec}}(x;s)$ and $s \mapsto u^R_{\textup{adv}}(x;s)$ are monotone increasing.  
\end{lemma}
  \begin{proof} 
  Let $-R \leq s < s'$ then $u^R_{\textup{adv}}(x;s') \wedge u^R_{\textup{adv}}(x;s)$ is a supersolution of \eref{finite-strip-1} and so $u^R_{\textup{adv}}(x;s') \wedge u^R_{\textup{adv}}(x;s) \geq u^R_{\textup{adv}}(x;s)$.  Similar arguments show the monotonicity of $s \mapsto u^R_{\textup{rec}}(x;s)$.
  \end{proof}

We will now send $R\to\infty$ to construct extremal solutions in $\R^d$.

\begin{lemma}\label{l.limit-R-defn-uadv}
Let $s\in \R$. For any $R_k \to \infty$, along a subsequence $u^{R_k}_{\textup{adv}}(\cdot;s)$ converges to a minimal supersolution $u_{\textup{adv}}(\cdot;s)$ of \eref{single-site}, locally uniformly in $\R^d$. Moreover $u_{\textup{adv}}(\cdot;s)$   is {proper} and has height $s$. 
\end{lemma}
The proof of Lemma~\ref{l.limit-R-defn-uadv} is parallel to and simpler than that of the next Lemma. The receding case is more delicate due to the lack of nondegeneracy in maximal subsolutions. We will thus only present the proof in this case.  Similar to the idea in the proof of \tref{monotone-family2} we construct a monotonic approximating sequence, the details are a bit different from the $d=2$ case since the monotonicity relies on taking advantage of the far-field behavior.

\begin{lemma}\label{l.limit-R-defn-urec}
    Let $s\in \R$. For any $R_k \nearrow \infty,$ there is a subsequence $R_{k_j}$ with $s_j \nearrow s$ so that $u^{R_{k_j}}_{\textup{rec}}(\cdot;s_j)$ converges  to a maximal subsolution $u_{\textup{rec}}(\cdot;s)$ of \eref{single-site} locally uniformly in $\R^d$. Moreover, $u_{\textup{rec}}(\cdot;s)$ is {proper} with height $s$. 
\end{lemma}

\begin{proof} 

 By \lref{finite-strip-uniform-bdry-layer} for any $R \geq R'$, $s \geq s'$, and $R' \geq 2|s|$:
\begin{equation}\label{e.increasing-R-s}
    u^R_{\textup{rec}}(x;s)-u^{R'}_{\textup{rec}}(x;s') \geq s-s' - CR'^{2-d} \ \hbox{ on } \ (\R^d \setminus B_{R'}(0)) \cap \overline{\{u^{R'}_{\textup{rec}}(\cdot;s') >0\}}.
\end{equation}
We may assume that $R_k \geq \max\{1,2|s|\}$ for all $k$. Let $s_j := s-\frac{1}{j}$ so that $s_{j+1}-s_j>0$ and $s_{j}\nearrow s$. Now choose a subsequence $R_{k_j}$ so that $CR_{k_j}^{2-d}<s_{j+1}-s_j$. Call $u_j :=u^{R_{k_j}}_{\textup{rec}}$.  With this choice both $u_{j+1}$ and $u_j$ are maximal subsolutions of \eref{single-site} in $B_{R_{k_j}}(0)$ and, by \eref{increasing-R-s} and choice of $R_{k_j}$ we have $u_{j+1}>u_j$ on $(\R^d \setminus B_{R_{k_j}}(0)) \cap \overline{\{u_j>0\}}$.  Thus the maximal subsolution property of $u_{j+1}$ in $B_{R_{k_j}}(0)$ gives $u_{j+1} \geq u_j$ in  $\R^d$.

By \eqref{e.strip-solution-bounds} and \lref{lipschitz-estimate}, $u_j$ are uniformly bounded and uniformly Lipschitz continuous. Thus, along a subsequence, the local uniform convergence to some $u$ follows. Since the sequence was monotone we can apply \cref{monotone-stability} and \lref{extremal-stability} to conclude that $u$ is a maximal subsolution. Taking the limit in \lref{finite-strip-uniform-bdry-layer} gives that $u(x) = x_d + s +o(1)$ in $\{u>0\}$ as $|x| \to \infty$.

\end{proof}

\subsection{ Proof of \tref{monotone-family}.}
  Let $u_{\textup{adv}}$ and $u_{\textup{rec}}$ from \lref{limit-R-defn-uadv} and \lref{limit-R-defn-urec}. Their monotonicity in $s$ follows from Lemma~\ref{l.monotone}.  The asymptotic expansion of these functions, as well as  existence and finiteness of the capacity functions $\varkappa_{\textup{adv/rec}}(s)$,  follows from \cite{FKexterior}*{Theorem 4.1}. This concludes part \ref{part.minimal-super-family} of the theorem.

To show part \ref{part.ordering}, let $u$ be a {proper} solution of \eref{single-site} with height $s$ and capacity $k$. Let us choose $s'<s$. We first observe that
  $u_{\textup{adv}}(x;s') \leq u(x)$.
 To see this, let $r$ be sufficiently large so that \eqref{expansion} yields
  \[u_{\textup{adv}}(x;s')<u(x) \ \hbox{ on } \ \overline{\{u_{\textup{adv}}(\cdot;s') > 0 \}} \setminus B_r.\]
  Then minimality of $u_{\textup{adv}}(x;s')$ implies that $u_{\textup{adv}}(x;s') \leq u(x) $.

  By \cite{FKexterior}*{Theorem 4.1} there is $C_0=C_0(d, \min q, \max q)$ so that for any $v$ solving \eref{single-site} in $\R^d$ with capacity $\kappa$ and height $\sigma$, we have
  \[|\kappa-R^{d-2}[ (R+\sigma)-v(Re_d)]| \leq C_0 R^{-1} \quad\hbox{ for any } R\geq 1. \]
  Applying this to $v=u_{\textup{adv}}(\cdot;s')$ and $v=u$,  as well as $u_{\textup{adv}}(\cdot;s') \leq u $,  we find:
  \begin{align*}
     \varkappa_{\textup{adv}}(s') &\geq  R^{d-2}[ (R+s')-u_{\textup{adv}}(Re_d;s')] - C_0R^{-1} \\
     & \geq R^{d-2}[ (R+s)-u(Re_d)]-R^{d-2}(s-s') - C_0R^{-1}\\
     &\geq k -R^{d-2}(s-s') - 2C_0R^{-1}.
  \end{align*}
  Taking the limit $s'\nearrow s$ first and then $R \to \infty$ yields $\varkappa_{\textup{adv}}(s-) \geq k$.  Note that the uniform constant $C_0$ from \cite{FKexterior}*{Theorem 4.1} is crucial here.

If $u(x)$ is a minimal supersolution, then arguing as above yields $u(x) \leq  u_{\textup{adv}}(x;s')$ for any $s'>s$ and then
\[k \leq \varkappa_{\textup{adv}}(s+).\]

The arguments are parallel for the receding case.
  
\hfill$\Box$

\subsection{Proof of \tref{capacity-properties}}

Arguing as in the proof of \tref{monotone-family} using \cite{FKexterior}*{Theorem 4.1} we find that
\[|\varkappa_{\textup{adv}}(s)-R^{d-2}[ (R+s')-u_{\textup{adv}}(Re_d;s')]| \leq C_0R^{-1}\]
and so
\[
|\varkappa_{\textup{adv}}(s)-\varkappa_{\textup{adv}}(s')| \leq R^{d-2}|u_{\textup{adv}}(Re_d;s)-u_{\textup{adv}}(Re_d;s')|+R^{d-2}(s-s')+2C_0R^{-1}.
\]
Thus if $s\mapsto u_{\textup{adv}}(\cdot;s)$ is continuous at $s$ then taking the limit above $s' \to s$ first and then $R \to \infty$ yields that $\varkappa_{\textup{adv}}$ is continuous at $s$. Arguing similarly, using that $u_{\textup{adv}}(Re_d;s\pm)$ left and right limits exist at every $s$ shows that left and right limits of $\varkappa_{\textup{adv}}$ exist at every $s \in \R$. The receding case is argued the same. This shows \ref{part.capacity-properties-0} and \ref{part.capacity-properties-3}. Also \ref{part.capacity-properties-4} follows now from \tref{monotone-family}\ref{part.ordering}.

Next, we explain the semi-continuity properties of $\kappa_{\textup{adv}}$ and $\kappa_{\textup{rec}}$. We check the receding case since advancing can be addressed similarly.  Due to the right continuity of $\kappa_{\textup{rec}}$ lower semicontinuity is equivalent to $\kappa_{\textup{rec}}(s-) \geq \kappa_{\textup{rec}}(s)$.  Note that $u_{\textup{rec}}(\cdot;s-) := \lim_{s' \nearrow s} u_{\textup{rec}}(\cdot;s')$  is a proper solution of \eref{single-site} with height $s$ by \cref{monotone-stability}. By the same arguments as in the previous paragraph, applying \cite{FKexterior}*{Theorem 4.1} again, $u_{\textup{rec}}(\cdot;s-)$ has capacity $\kappa_{\textup{rec}}(s-)$. Since $u_{\textup{rec}}(\cdot;s-)$ is another {proper} solution of \eref{single-site} with height $s$,  \tref{monotone-family}\ref{part.ordering} implies that $\kappa_{\textup{rec}}(s-) \geq \varkappa_{\textup{rec}}(s+)=\kappa_{\textup{rec}}(s)$.

 Finally we check the compact support \ref{part.capacity-properties-1}. This follows from \cite{FKexterior}*{Theorem 4.1}, which yields the universal bound $|s| \leq C_0(d,\max |q|)$ for any proper solution with $k \neq 0$. This shows that $\varkappa_{\textup{adv/rec}}(s) = 0$ for all $|s| > C_0(d,q)$ and we conclude \ref{part.capacity-properties-1}.

\hfill$\Box$

\section{Example of nonzero pinning force via linearization}\label{s.nonzero-example}

In this section we will show that the capacities are nontrivial for a generic family of $Q$.  As usual we will work on the case $e = e_d$, results for general normal direction $e \in S^{d-1}$ can be obtained by rotation.  Let $\Phi$ denote the, un-normalized, upward pointing fundamental solution
\begin{equation}\label{e.fundamental}
    \Phi(x) := \begin{cases}
    -\log |x| & d =2\\
    |x|^{2-d} & d \geq 3
\end{cases} 
\end{equation}
and the normalizing constant $\gamma_d$ is defined
\begin{equation}\label{e.normalizing}
    \hbox{$\gamma_2 = \pi$ and $\gamma_d = \frac{1}{2}(d-2)|\partial B_1|$ in $d \geq 3$.}
\end{equation}
Our main tool is an asymptotic expansion of the capacity for  defects of the form $Q = 1+ q$ with small $q$. 

\begin{theorem}\label{t.nonzero}
   Let $d \geq 2$. Let $Q$ be of the form $Q=1+q$ with $q=\sigma \tilde{q}$,  $\tilde{q}\in C^{0,1}_c(B_1(0))$, and $ 0 < \sigma \leq \sigma_0$ sufficiently small depending on $\|\tilde{q}\|_{C^{0,1}}$. For every $s \in \R$ there is a solution $u$ of \eref{single-site} such that 
   \[ |u(x) - (x_d-s)| \leq C\sigma |\Phi(2+|x|)| \ \hbox{ in } \ \{u>0\}\]
   with capacity 
    \begin{equation}\label{e.linear-response}
    k = \frac{\sigma}{\gamma_d}\int_{\{x_d = s\}} \tilde{q}\ dS+O(\sigma^2).
    \end{equation}
       As a consequence 
    \[[k_{\textup{rec}},k_{\textup{adv}}] \supset  \left[\frac{\sigma}{\gamma_d}\min_{s}\int_{\{x_d = s\}} \tilde{q}\ dS+O(\sigma^2),\,\,\frac{\sigma}{\gamma_d}\max_s \int_{\{x_d = s\}} \tilde{q}\ dS+O(\sigma^2)\right].\]
 \end{theorem}

\tref{nonzero} is for small $q$. However, we can still use this information for large $q$ by exploiting monotonicity and locality. Note that, for $0<a<1$, $\kappa_{\textup{adv}}(s)$ decreases if we replace $q(x)$ with $aq(x)_+ - q(x)_-$. Similarly $\kappa_{\textup{rec}}(s)$ increases if we replace $q(x)$ with $q(x)_+ - aq(x)_-$. It follows easily from the definitions of $\kappa_{\textup{adv/rec}}$ that these quantities are monotone increasing with respect to $q$. Combining this fact with \tref{nonzero} we can deduce some information in the nonlinear regime:

\begin{corollary}\label{c.nonzero-ks}
  Suppose that $Q = 1+q$ with $q\in C_c^1(B_1(0))$.
  \begin{itemize}
\item If there is some interval $I \subset \R$ such that $q|_{\{x_d \in  I\}}$ is non-negative and not identically zero, then $k_{\textup{adv}}>0$. 
\item If there is some interval $I \subset \R$ such that $q|_{\{x_d \in  I\}}$ is non-positive and not identically zero, then $k_{\textup{rec}}<0$.
       \end{itemize}
 \end{corollary}
 We emphasize that the difference with \tref{nonzero} is that there is no requirement here that the defect $q$ be small in magnitude.  Based on the criteria in \cref{nonzero-ks} it is not so hard to arrange scenarios where $k_{\textup{adv}}(e) > 0 > k_{\textup{rec}}(e)$ for \emph{all} directions $e \in S^{d-1}$ see \fref{arrangement}. 

 \begin{figure}
     \centering
     \input{figures/ellipses-pic.tikz}
     \caption{Defect $q$ so that for every direction $e$ there exist two planes $P_\pm = \{x \cdot e = s_\pm\}$ such that $\pm q \geq 0$ on $P_\pm$ and are nontrivial on $P_\pm$.}
     \label{f.arrangement}
 \end{figure}

We proceed to prove \tref{nonzero}. The desired solution will be constructed in hodograph coordinates (see \cite{FKexterior}*{Section 2.4}), for linearization purposes.
The corresponding PDE  is
\begin{equation}\label{e.hod_nonlinear}
\begin{cases}
    \textup{tr}(A(\grad_y v)D^2_yv) = 0 & \hbox{ in } \ \{y_d>0\},\\
    (1+\sigma \tilde{q}(y',v))(1+\partial_{y_d}v) = \sqrt{1+|\grad_y 'v|^2} & \hbox{ on } \partial \{y_d>0\}.
\end{cases}
\end{equation}
Here $A(p)$ is the coefficient matrix of the hodograph PDE, see \cite{FKexterior}*{Section 2.4},
\[A(p) = \left[\begin{array}{cc}
    \textup{I}_{d-1} & (1+p_d)^{-1}p' \\
    (1+p_d)^{-1}(p')^T & \frac{1+|p'|^2}{(1+p_d)^2}
\end{array}\right],\]
which is uniformly elliptic with $\|A(p)-I\| \leq C|p|$ for $|p| \leq 1/2$.

\begin{proposition}\label{p.nonzero-hodograph}
  Let $0<\sigma\ll 1$. Then for every $s \in \R$ there is a solution $v_s$ of \eref{hod_nonlinear} of the form
  \[
v_{s}(y) = s+\sigma w_s(y) + \sigma^2 r(y),
\]
 where $|r(y)| \leq C|\Phi(2+|y|)|$ and 
 \[w_s(y) := \frac{1}{\gamma_d}\int_{\partial \R^d_+} \Phi(y-z') \tilde{q}(z',s) \ dS(z') \ \hbox{ for } \ y \in \R^d_+,\]
 i.e. $w_s$ solves the linearized problem with decay as $x_d\to\infty$:
\begin{equation}\label{e.hodograph-PDE}
\begin{cases}
    \Delta w_s = 0 & \hbox{ in } \ \R^d_+, \\
    \partial_{y_d} w_s = -\tilde{q}(y',s) & \hbox{ on } \partial \R^d_+.
\end{cases}
\end{equation}
In particular the capacity $k_s$ of the solution $v_s$, as in \cite{FKexterior}*{Theorem 3.2}, satisfies
\[k_s = \frac{\sigma}{\gamma_d} \int_{\partial \R^d_+} \tilde{q}(z',s) \ dS(z')+O(\sigma^2).\]
\end{proposition}

\tref{nonzero} is an immediate corollary of \pref{nonzero-hodograph} by taking the inverse hodograph transform of the solution $v_s$ which, as per the discussion in \cite{FKexterior}*{Section 2.4}, produces a solution $u_s$ of \eref{single-site} with the same capacity.

\begin{proof}
  Take $s = 0$ and we will omit the $s$ subscripts, the argument is the same for other values of $s \in \R$. We will write the nonlinear PDE solved by $r =\sigma^{-2}(v(y)  - s-\sigma w(y))$, prove that it has a solution with the desired bounds, and then define $v(y) := s + \sigma w(y)  + \sigma^2 r(y)$ to attain the desired solution of \eref{hod_nonlinear}. 

    The nonlinear Neumann problem for $r$ has the form
    \begin{equation}\label{e.hod_nonlinear_error-PDE}
\begin{cases}
    \textup{tr}((I+\sigma B_\sigma(y,\grad r))D^2r) +\textup{tr}(B_\sigma(y,\grad r)D^2w)=0 & \hbox{ in } \ \{y_d>0\}\\
\partial_{y_d} r = N_\sigma(y,r,\grad'r) & \hbox{ on } \partial \{y_d>0\}.
\end{cases}
\end{equation}
The coefficients are defined so that if $r$ solves \eref{hod_nonlinear_error-PDE} then $v(y) := s + \sigma w(y)  + \sigma^2 r(y)$ solves \eref{hod_nonlinear}. More specifically
      \[B_\sigma(y,p):= \sigma^{-1}(A(\sigma\grad w+\sigma^2p)-I)\]
      and
      \[N_\sigma(y,r,p') := \frac{1}{\sigma^2}\left[\frac{\sqrt{1+|\sigma\grad_y' w + \sigma^2 p'|^2}}{1+\sigma \tilde{q}(y',\sigma w + \sigma^2 r)}-(1+\sigma \tilde{q}(y',0))\right].\]
      From now on we drop the $\sigma$ subscript since that is fixed through the proof. Note that $B(y,r,p)$ and $N(y,r,p')$ satisfy the following bounds
\begin{equation}\label{e.B-N-estimates}
\begin{split}
        |B| &\leq C((2+|y|)^{1-d}+\sigma |p|)\\
         |\grad_p B| &\leq C\sigma \\
      |N| &\leq C((2+|y|)^{2-2d}+\sigma^2|p'|^2),\\
      |\grad_{p'} N| &\leq C( \sigma(2+|y|)^{1-d}+\sigma^2|p'|), \\  |\partial_r N| &\leq C(1+|p'|)\sigma {\bf 1}_{B_1}(y)
\end{split}
\end{equation}

We will do a fixed point argument to deal with the non-proper $r$ dependence in the nonlinear Neumann condition in \eref{hod_nonlinear_error-PDE}. We define the functional space setting. First define $C^1$ norm weighted by the fundamental solution scaling
\[\|r\|_\Phi := \sup_{y \in \R^d_+}\left[\frac{|r(y)|}{ |\Phi(2+|y|)|}+\frac{|\grad r(y)|}{|\grad \Phi(2+|y|)|}\right].\]
Then define
\[X := \left\{r \in C^1(\R^d_+ \cup \partial \R^d_+) : \|r\|_\Phi\leq M\right\}.\]
This is a complete metric space under the $d_{\Phi}(r,r'):=\| r-r'\|_\Phi$ metric. 
For each $r \in X$ we will show there exists a solution $\rho =: \mathbf{\Psi}[r]$ of
\begin{equation}\label{e.hod_nonlinear_error-PDE-FP}
\begin{cases}
    \textup{tr}((I+\sigma B(y,\grad r))D^2\rho) +\textup{tr}(B(y,\grad r)D^2w)=0 & \hbox{ in } \ \{y_d>0\}\\
\partial_{y_d} \rho = N(y,r,\grad'r) & \hbox{ on } \partial \{y_d>0\}
\end{cases}
\end{equation}
given by the Neumann Green's kernel for the elliptic operator $\textup{tr}((I+\sigma B(y,\grad r)) \cdot)$ in $\R^d_+$. For $M$ sufficiently large and $\sigma>0$ sufficiently small we will show that $\mathbf{\Psi} : X \to X$ and is a contraction. In particular assume that $\sigma M \leq 1$. We do not show all details from here, but just give a sketch.

To show that $\Psi(X) \subset X$ note that, applying the bounds \eref{B-N-estimates}, $\rho = \Psi[r]$ solves
\begin{equation}\label{e.hod_nonlinear_error-PDE-Mbound-FP}
    \begin{cases}
    |\textup{tr}((I+\sigma B(y,\grad r))D^2\rho)|\leq  C(1+\sigma M)(2+|y|)^{1-2d} & \hbox{ in } \ \{y_d>0\}\\
|\partial_{y_d} \rho| \leq C(1+\sigma^2M^2)(2+|y|)^{2-2d} & \hbox{ on } \partial \{y_d>0\}.
\end{cases}
\end{equation}
Recall that, given $M$, we will choose $\sigma$ small so that $\sigma M \leq 1$. By standard Neumann kernel bounds for variable coefficient operators applied in \eref{hod_nonlinear_error-PDE-Mbound-FP}
\[|\rho| \leq C|\Phi(2+|y|)| \ \hbox{ and } \ |\grad \rho | \leq C|\grad \Phi(2+|y|)|\]
so we choose $M$ sufficiently large depending on the dimensional constants in the previous line so that $\rho = \Psi[r] \in X$. 

To show the contraction property: let $r,r' \in X$ and define $\nu = \mathbf{\Psi}[r]-\mathbf{\Psi}[r']$. Then $\nu$ solves
\begin{equation}\label{e.hod_nonlinear_error-PDE-lipschitz-FP}
\begin{cases}
    |\textup{tr}((I+\sigma B(y,\grad r))D^2\nu)| \leq C\sigma(2+|y|)^{1-2d}\|r-r'\|_{\Phi} & \hbox{ in } \ \{y_d>0\}\\
|\partial_{y_d} \nu| \leq C\sigma(2+|y|)^{2-2d}\|r-r'\|_{\Phi} & \hbox{ on } \partial \{y_d>0\}.
\end{cases}
\end{equation}
Applying Neumann kernel bounds again in \eref{hod_nonlinear_error-PDE-lipschitz-FP}
\[|\nu| \leq C\sigma \|r-r'\|_{\Phi}|\Phi(2+|y|)| \ \hbox{ and } \ |\grad \nu| \leq C\sigma \|r-r'\|_{\Phi}|\grad \Phi(2+|y|)|.\]
So $\|\nu\|_{\Phi} \leq C\sigma \|r-r'\|_{\Phi}$ and for $\sigma$ sufficiently small this makes $\mathbf{\Psi}$ a contraction mapping of $X$.

Finally the computation of the expansion for the capacity just follows by computing $k =\lim_{|y| \to \infty} \Phi(y)^{-1}v(y)$, using the bound $|r(y)| \leq C|\Phi(2+|y|)|$, and observing that $\Phi(y)^{-1}\Phi(y-z') \to 1$ as $|y| \to \infty$ uniformly on $z'$ in the support of $\tilde{q}(z',0)$.

 \end{proof}

\section{Pinned slopes in a periodic medium}\label{s.periodic-prelim}
In this section we outline the proof of \tref{asymptotics},  and provide necessary ingredients to prepare for \sref{homogenization-lower-bound} and \sref{homogenization-upper-bound},
where the proof will be carried out.

\subsection{Outline of the proof of \tref{asymptotics}} 
The proof of \tref{asymptotics} is divided into two parts, the lower bound away from $1$ and the upper bound away from $1$, respectively: the lower bounds
\begin{equation}\label{e.lower_bound}
    \liminf_{\delta \to 0} \frac{Q^\delta_{\textup{adv}}-1} {\delta^{d-1}} \geq \gamma_d|\xi|^{-1}k_{\textup{adv}}, \quad \liminf_{\delta \to 0}\frac{1-Q^{\delta}_{\textup{rec}}}{\delta^{d-1}} \geq -\gamma_d|\xi|^{-1}k_{\textup{rec}}\end{equation}
and the upper bounds
\begin{equation}\label{e.upper_bound}
   \limsup_{\delta \to 0} \frac{Q^\delta_{\textup{adv}}-1}{\delta^{d-1}} \leq \gamma_d|\xi|^{-1}k_{\textup{adv}}, \quad  \limsup_{\delta \to 0} \frac{1-Q^\delta_{\textup{rec}}}{\delta^{d-1}} \leq -\gamma_d|\xi|^{-1}k_{\textup{rec}}.
\end{equation}

Both of these bounds rely on a certain cell problem \eqref{e.CM-cell}, to be introduced in \sref{semipermeable-cell}.

 In \sref{homogenization-lower-bound}, we establish \eref{lower_bound}. For this it suffices to construct a plane-like supersolution with the desired slope. Such a barrier will be constructed in \lref{patched-barrier}, by patching together an outer solution, based on $\omega$ from \eref{CM-cell}, with an inner solution near the defect based on one of the single site profiles constructed in Sections \ref{s.pinned-solutions-2d} and \ref{s.monotone-family}.

In \sref{homogenization-upper-bound} we proceed to \eref{upper_bound}: Here we show that the barrier constructed above is essentially optimal.  This is naturally a more difficult direction, and will be the content of \sref{homogenization-initial-flatness} - \ref{s.homogenization-capacity}.   Due to the characterization of the pinning interval endpoints in \tref{cell-solutions-reference} from \cite{feldman2021limit}, there is a plane-like solution $u_\delta$ of \eref{general-Q-bernoulli} achieving the maximal macroscopic slope $Q_{\textup{adv}}$ and satisfying the strong Birkhoff property (see \dref{birkhoff-prop} below).  We will show that $u_{\delta}$ matches the profile of our barrier away from the defects, which amounts to \tref{asymptotics}. An important ingredient in achieving this is the asymptotic decomposition of $u_{\delta}$ via the semipermeable membrane cell problem \eqref{e.CM-cell-ptwise}, based on refined compactness properties of $u_{\delta}$ using its periodicity and  near-defect behavior understood from earlier sections. The uniqueness / rigidity properties of \eqref{e.CM-cell-ptwise} then allow us to compute the asymptotic upper bound.

\subsection{Rotation to $e=e_d$ case} As usual in this paper, we will normalize the problem to reduce to the direction $e=e_d$, paying attention to the fact that the periodicity lattice is now rotated as well. We use the fact that $e = \xi/|\xi|$ is rational and $\xi \in \Z^d$ is irreducible to explain the structure of the rotated lattice, especially the intersection of the lattice with the plane $\partial \R^d_+$.

For $v\in \R^d,$ denote $v^\perp$ to be the linear subspace orthogonal to $v$. For $\xi \in \Z^d$, irreducible, the lattice $\Z^d \cap \xi^\perp$ is $(d-1)$-dimensional. The area of the unit period cell of this lattice is $|\xi|$. We will now perform a rotation $O$ sending $e$ to $e_d$.  This also rotates the underlying periodicity lattice to $O\Z^d$. We define
$$
\mathcal{Z}:=O(\Z^d \cap \xi^\perp) \subset \{x_d=0\} = \partial \R^d_+.
$$
Note that $\mathcal{Z}$ is a $(d-1)$-dimensional integer lattice.  Fix a basis $\zeta_1,\dots,\zeta_{d-1}$ of $\mathcal{Z}$ (for instance a reduced basis). Let $\Box_{\mathcal{Z}} \subset \{x_d=0\}$ be the fundamental domain / parallelepiped of the lattice $\mathcal{Z}$ centered at $0 \in \mathcal{Z}$ corresponding to this basis. The volume of the fundamental domain is (independent of the choice of basis actually)
\begin{equation}\label{e.box-size}
|\Box_{\mathcal{Z}}|= |\xi|.
\end{equation} 
 We also denote  
 \begin{equation}\label{e.cube-notation-def}
 \cube^{\mathcal{Z}}:= \Box_{\mathcal{Z}} \times [-1,1] \ \hbox{ and } \ \cube^{\mathcal{Z}}_r:=r\cube^{\mathcal{Z}}.
 \end{equation}
 These fattened cells are convenient since the lattice translations (almost) disjointly cover the channel region $\{|x_d| \leq 1\}$. Finally we note that there is some $\rho_0(\xi)>0$ so that $B_{\rho_0}(0) \subset \cube^{\mathcal{Z}}$. In $d=2$ the lattice $\mathcal{Z}$ is generated by $\xi^\perp = (\xi_2,-\xi_1)$ so we can take $\rho_0 = \min\{|\xi|/2,1\} \geq 1/2$ in that case.

\subsection{Pinning interval in periodic homogenization}\label{s.pinning-interval-recall} In this section we collect some known results on pinning interval and plane-like solutions for the Bernoulli free boundary problem in periodic media:
    \begin{equation}\label{e.general-Q-bernoulli}
        \begin{cases}
            \Delta u = 0 & \hbox{in } \{u>0\}\\
            |\grad u | = Q(x) & \hbox{on } \partial \{u>0\},
        \end{cases}
    \end{equation}
  where $Q(x)>0$ is $\Z^d$-periodic and continuous.

    \begin{definition}\label{d.pinning-interval}
     Given $e \in S^{d-1}$ we define the \emph{advancing slope}
     \[Q_{\textup{adv}}(e) := \sup\{ \alpha : \exists \hbox{ a supersolution $u \in C(\R^d;\R_+)$ of \eref{general-Q-bernoulli} with $u(x) \geq \alpha (e \cdot x)_+$}\}\]
     and the \emph{receding slope}
     \[Q_{\textup{rec}}(e) := \inf\{ \alpha : \exists \hbox{ a nontrivial subsolution $u\in C(\R^d;\R_+)$ of \eref{general-Q-bernoulli} with $u(x) \leq \alpha (e \cdot x)_+$}\}.\]
     The \emph{pinning interval} is the interval
     \[I(e) := [Q_{\textup{rec}}(e),Q_{\textup{adv}}(e)].\]
     Note that $Q_{\textup{rec}}(e;Q)$, $Q_{\textup{adv}}(e;Q)$, and $I(e;Q)$ depend implicitly on the periodic field $Q(x)$, but we usually hide the dependence. When $Q = Q_\delta$ is the dilute defect medium of \eref{model} we write, as in the introduction, $Q^\delta_{\textup{adv}}(e) := Q_{\textup{adv}}(e;Q_\delta)$ and $Q^\delta_{\textup{rec}}(e) := Q_{\textup{rec}}(e;Q_\delta)$.
    \end{definition}
    It is straightforward to check that
    \[ 0<\min Q \leq Q_{\textup{rec}}(e) \leq Q_{\textup{adv}}(e) \leq \max Q<+\infty.\]
    Next we define the Birkhoff property which is an analogue of periodicity in this context.
    \begin{definition}\label{d.birkhoff-prop}
     $u$ satisfies the \emph{Birkhoff property} in direction $e$ if
    \[u(x+k) \geq u(x) \ \hbox{ for all } \ k \in \Z^d \hbox{ s.t. } k \cdot e >0.\]
     $u$ satisfies the \emph{strong Birkhoff property} if the same holds for all $k \in \Z^d$ with $k \cdot e \geq 0$. 
     \end{definition}
     
    In particular any $u$ with the strong Birkhoff property are periodic with respect to tangential lattice translations (if there are any such).
    
    We quote a result now summarizing several of the results in \cite{feldman2021limit} on the pinning interval and the existence of corresponding plane-like solutions of \eref{general-Q-bernoulli}. 
    
    \begin{theorem}[\cite{feldman2021limit}]\label{t.cell-solutions-reference}
    The following hold:
        \begin{enumerate}[label = (\alph*)]
            \item\label{part.cell-solutions-reference-1} (Energy minimizing slope) For all $e \in S^{d-1}$ we have $\langle Q^2 \rangle^{1/2} \in I(e)$.
            \item\label{part.cell-solutions-reference-2} (Existence of plane-like solutions in the pinning interval) For any $\alpha \in I(e)$ there exists a solution $u$ of \eref{general-Q-bernoulli} satisfying the strong Birkhoff property in direction $e$, and satisfying the flatness condition
            \begin{equation}\label{e.plane-like-defn}
                \alpha (x \cdot e - C)_+ \leq u(x) \leq  \alpha (x \cdot e + C)_+
            \end{equation}
             and $C \geq 0$ depends only on $d$, $\min Q$ and $\max Q$. Such solutions are called \emph{plane-like solutions}. Furthermore if $\alpha = Q_{\textup{adv}}(e)$ then $u$ can be taken to be a minimal supersolution, if $\alpha = Q_{\textup{rec}}(e)$ then $u$ can be taken to be a maximal subsolution, and if $\alpha = \langle Q^2 \rangle^{1/2}$ then $u$ can be taken to be an energy minimizer.
        \end{enumerate}
    \end{theorem}

\subsection{A cell problem of semi-permeable boundaries}\label{s.semipermeable-cell}

In this section we describe a certain cell problem, related to the homogenization theory of perforated domains, for example see \cite{PapanicolaouVaradhan1980} and many others. The solutions of this problem describe the asymptotic behavior of \eref{model} between the defects.  The main issue we want to clarify is the relation between integral solutions and pointwise solutions. The former characterization is preferable for computing the asymptotic slope, while the latter is easier to relate to the viscosity solutions of \eref{model}.  

Recall that we will work in a normalized setting, with the inward normal vector $e = e_d$ and a $(d-1)$-dimensional integer lattice $\mathcal{Z} \subset\partial \R^d_+$.
  Consider the semipermeable boundary cell problem
\begin{equation}\label{e.CM-cell}
    \left\{\begin{array}{ll}
    \Delta \omega = 0 & \hbox{in } \R^d_+ ,\\
    \partial_{x_d}\omega= 1-\sum_{z \in \mathcal{Z}} c_*\delta_z & \hbox{on } \partial\R^d_+,\\
   
    \end{array}
\right.
\end{equation} with
\begin{equation}\label{e.CM-cell-1}\langle \partial_{x_d}\omega(\cdot,x_d)\rangle' = 0 \hbox{ for all } x_d>0.
\end{equation}
Here $\delta_z$ denotes the Dirac mass at $z$. 
For now this PDE boundary value problem should be interpreted in the weak divergence form sense. Note that $c_*:=|\Box_{\mathcal{Z}}| $ is the unique constant  for which \eqref{e.CM-cell}-\eqref{e.CM-cell-1} are solvable: since $\langle \partial_{x_d}\omega(\cdot,x_d)\rangle'$ is independent of $x_d$ by \lref{periodic-bdry-layer} (applied to $\omega$ on $\{x_d>\ep\}$ for each $\ep>0$, where it is bounded and periodic), formally the existence of a solution to \eref{CM-cell}-\eqref{e.CM-cell-1} requires that
\[0 = \langle 1-\sum_{z \in \mathcal{Z}} c_* \delta_z \rangle' = 1-\frac{1}{|\Box_{\mathcal{Z}}|}c_* \ \hbox{ i.e. } \ c_* = |\Box_\mathcal{Z}| = |\xi|.\]
The formal argument can be made rigorous for variational weak solutions of \eref{CM-cell} using a standard test function choice.  Note that $\omega$ is also just a modified Green's function for a periodic problem for the half-Laplacian.

The pointwise formulation of \eref{CM-cell} and \eref{CM-cell-1} is 
\begin{equation}\label{e.CM-cell-ptwise}
    \left\{\begin{array}{ll}
    \Delta \omega = 0 & \hbox{in } \R^d_+ ,\\
    \partial_{x_d}\omega= 1 & \hbox{on } \partial\R^d_+ \setminus \mathcal{Z}\\
    \lim_{x \to z}\Phi(x-z)^{-1}\omega(x) = \frac{1}{\gamma_d}c_* & \hbox{for } z \in\mathcal{Z}, \\
    \langle \partial_{x_d}\omega(\cdot,x_d)\rangle' = 0 & \hbox{for } x_d>0.
    \end{array}
\right.
\end{equation}
where $\Phi$ denotes the, un-normalized, upward pointing fundamental solution and $\gamma_d$ is the normalizing constant, see \eref{fundamental} and \eref{normalizing} above.

The following lemma explains the relation between distributional and pointwise notions and has some useful information about removing singularities.

\begin{lemma}\label{l.divergence-non-divergence-hole-info}
    Suppose that $w \in C(\overline{B_1^+} \setminus \{0\})$ solves, for some $A \leq B$,
    \[\begin{cases}
        \Delta w = 0 & \hbox{in } \ B_1^+\\
        \partial_{x_d} w = 0 & \hbox{on } \ (\partial \R^d_+ \cap B_1) \setminus \{0\}\\
        \displaystyle\limsup_{x \to 0} w(x)/\Phi(x) \leq B \\
        \displaystyle\liminf_{x \to 0} w(x)/\Phi(x) \geq A.
    \end{cases}\]
    Then, calling $h$ to be the harmonic function in $B_1$ which is equal to  $w$ on $\partial B_1$ (extended to the lower half sphere by even reflection),
    \begin{equation}
        \label{e.removable-sing-0}
    w(x) = c (\Phi(x)- \Phi(e_d)) + h(x) \ \hbox{ for some } A \leq c \leq B.
    \end{equation}
    Furthermore $w$ solves, in the standard distributional weak sense,
      \[\begin{cases}
        \Delta w = 0 & \hbox{in } \ B_1^+\\
        \partial_{x_d} w = -\gamma_dc\delta_0 & \hbox{on } \ \partial \R^d_+ \cap B_1.
    \end{cases}\]
\end{lemma}
\begin{proof}
    We give a quick sketch of this classical result. By even reflection we can suppose that $w$ is harmonic in $B_1 \setminus \{0\}$ with the same bounds. Then there is $A\leq k\leq B$ such that $\lim_{x\to 0} w(x)/\Phi(x) = k$. This can be proved via a Kelvin transform argument similar to the proof of \cite{FKexterior}*{Theorem 3.2}. Then let $h(x) := w(x) - k (\Phi(x)-\Phi(e_d))$. Then $h(x) = w(x)$ on $\partial B_1$.  Note that  $\lim_{|x| \to 0} h(x)/\Phi(x) = 0$, thus standard removable singularity arguments show that $h$ is harmonic in $B_1$. Then we conclude by using the distributional derivatives of $\Phi$, noting that $\Delta \Phi = -2 \gamma_d \delta_0$ and the factor of $\frac{1}{2}$ arises from integrating over an upper half sphere instead of the full sphere.
 
\end{proof}

One useful consequence of the previous Lemma is the local decomposition / extension for solutions of \eref{CM-cell-ptwise}:
\begin{equation}\label{e.removable-sing-omega}
  \omega(x) = \tfrac{1}{\gamma_d}c_*\Phi(x) + x_d + h(x) \ \hbox{ in } \ \cube^{\mathcal{Z}}_{3/2} \cap \R^d_+
  \end{equation}
  with $h$ a bounded harmonic function in $\cube^{\mathcal{Z}}_{3/2}$ which is even symmetric with respect to $x_d\mapsto -x_d$. Since the right hand side naturally extends to the lower half cube $\cube^{\mathcal{Z}}_{3/2} \cap \R^d_-$ we can also view that as an extension of the left hand side.

\begin{theorem}\label{t.membrane-cell-ptwise}
    There exists a $\mathcal{Z}$-periodic solution $\omega$ of \eref{CM-cell-ptwise} which is bounded on $\{x_d \geq 1\}$ if and only if $c_* = |\Box_\mathcal{Z}|$.  Furthermore any two such solutions differ by a constant.
\end{theorem}

\begin{proof}
We omit the proof of existence: the existence of modified Green's functions for elliptic problems on compact manifolds is classical, see for example \cite{AubinThierry1998Snpi}*{Section 2.3}.

If $\omega$ solves \eref{CM-cell-ptwise} then \lref{divergence-non-divergence-hole-info} implies that $\omega$ is a variational weak solution of \eref{CM-cell}-\eqref{e.CM-cell-1} and then \lref{periodic-bdry-layer} below implies that $c_* = |\Box_{\mathcal{Z}}| = |\xi|$.

   Finally we discuss uniqueness, suppose that $\omega_0$ and $\omega_1$ both solve \eref{CM-cell-ptwise}. Let $h(y):= \omega_1(y) - \omega_0(y)$. Since $h$ is harmonic in $\R^d_+$ with zero Neumann data on $\partial \R^d_+ \setminus \mathcal{Z}$ we can extend $h$ by even reflection to be harmonic on $\R^d \setminus \mathcal{Z}$. 
   
    Due to  the asymptotic condition at the holes 
    \[\lim_{y \to z} \frac{h(y)}{\Phi(y-z)} = 0 \ \hbox{ for each } \ z \in \mathcal{Z},\]
    a standard removable singularity result for harmonic functions implies that $h$ is an entire harmonic function. Since $h$ is $\mathcal{Z}$-periodic and bounded on $\{x_d \geq 1\}$ by assumption, hence on $|x_d| \geq 1$ by the even reflection, it is also bounded on $|x_d| \leq 1$ by maximum principle, and so by Liouville it is constant.
\end{proof}
\begin{lemma}\label{l.periodic-bdry-layer}
    Suppose that $u$ solves
    \[\Delta u = 0 \ \hbox{ in } \ \{x_d >0\}\]
    and $u$ is bounded, and $\mathcal{Z}$-periodic for some integer lattice $\mathcal{Z}$ with $\textup{span}(\mathcal{Z}) = \R^{d-1} \times \{0\}$. Then $s:=\lim_{x_d \to +\infty} u(x)$ exists,
    \[|u(x) - s| \leq Ce^{-cx_d} \ \hbox{ and } \ |D^ku| \leq C_ke^{-cx_d}.\]
    Furthermore
    \[ \langle \partial_{x_d}u(\cdot,x_d) \rangle' = 0 \ \hbox{ for all } \ x_d>0\]
    where $\langle \cdot \rangle'$ denotes the average over the fundamental domain of the lattice $\mathcal{Z}$.
\end{lemma}
\begin{proof}
The exponential convergence follows from periodicity, maximum principle (here is where boundedness is used), and interior elliptic estimates. See, for example, \cite{FeldmanZhang2019}*{Lemma 5.2}.  The formula for the average of the normal derivative follows from the integration by parts identities:
    \[\frac{d}{dx_d}\langle \partial_{x_d}u(\cdot,x_d) \rangle' = \langle \partial_{x_d}^2u(\cdot,x_d) \rangle' = -\langle \Delta' u(\cdot,x_d) \rangle' = 0.\]
    The limit $\partial_{x_d}u(\cdot,x_d) \to 0$ as $x_d \to \infty$ fixes the constant of integration.
    
\end{proof}

\section{Asymptotic lower bound}\label{s.homogenization-lower-bound}
In this section we create periodic sub and supersolution barriers for the full problem \eqref{e.model}, by taking a solution of the single-site problem and patching it with a solution of \eqref{e.CM-cell}. As a corollary we will obtain the asymptotic lower bounds:

\begin{proposition}
    \label{p.asymptotic-exp-lower-bound}
Let $\xi \in \Z^d$ irreducible and $e = \xi/|\xi|$ then 
    \[
\liminf_{\delta \to 0} \frac{Q^\delta_{\textup{adv}}(e)-1}{\delta^{d-1}} \geq \gamma_d|\xi|^{-1}k_{\textup{adv}}(e) \ \hbox{ and } \ \liminf_{\delta \to 0} \frac{1-Q^\delta_{\textup{rec}}(e)}{\delta^{d-1}} \geq -\gamma_d|\xi|^{-1}k_{\textup{rec}}(e).
\]
\end{proposition}

To construct our patched barrier, first we discuss the {\it outer solution}, namely the part of the barrier away from the defects, based on \eqref{e.CM-cell}. This will be based on a translation of the profile
\begin{equation}\label{e.model-outer-soln}
    \tilde{w}_{out}(x) := (1+\sigma\alpha)x_d +\gamma_d\delta^{d-1}|\xi|^{-1}k[x_d - \omega(x)].
\end{equation}
Here $\sigma \in \{-1,+1\}$, $\alpha>0$, and $k \in \R$ are parameters to be adjusted, and $\omega$ is the unique $\mathcal{Z}$-periodic solution of \eref{CM-cell-ptwise} normalized by $\min_{\R^d_+} \omega = 0$, see \tref{membrane-cell-ptwise}.

\begin{lemma}\label{l.outer-solution-alpha-info}
Let $\tilde{w}_{out}$ be given by \eqref{e.model-outer-soln} with $0<\alpha\leq 1/2$, and suppose that for some $r\in [\delta, 1/2]$  and for a sufficiently large $C_0 = C_0(d, \mathcal{Z})$ we have \begin{equation}\label{e.alpha}C_0k^2(\delta/r)^{2(d-1)}r^{d-2}\Phi(r) \leq \alpha.
\end{equation}

\begin{enumerate}[label = (\roman*)]
    \item If $\sigma=-1$  then $\tilde{w}_{out}$ is a supersolution of \eref{bernoulli} in $\R^d\setminus \cup_{z \in \mathcal{Z}} B_r(z) $, i.e.,
\[|\grad \tilde{w}_{out}(x)|^2 \leq 1 \ \hbox{ on } \ \partial \{\tilde{w}_{out}>0\} \setminus \cup_{z \in \mathcal{Z}} B_r(z).\]
\item
If $\sigma=1$ then $\tilde{w}_{out}$ is a subsolution of \eref{bernoulli} in $\R^d\setminus \cup_{z \in \mathcal{Z}} B_r(z)$, i.e.,
\[|\grad \tilde{w}_{out}(x)|^2 \geq 1 \ \hbox{ on } \ \partial \{\tilde{w}_{out}>0\}\setminus \cup_{z \in \mathcal{Z}} B_r(z).\]
\end{enumerate}
\end{lemma}

\begin{proof}
By periodicity it suffices to show the result in $\cube^{\mathcal{Z}} $. First we estimate the location of the free boundary. Denoting $\beta:= (1+\sigma\alpha)$, for $x \in \partial \{\tilde{w}_{out}>0\}\cap\cube^{\mathcal{Z}} $ we have
\[ 0 = \tilde{w}_{out}(x) = \beta x_d + \gamma_d\delta^{d-1}|\xi|^{-1}k[x_d - \omega(x)].\]
  \eref{CM-cell-ptwise} and \lref{divergence-non-divergence-hole-info}  yields that
\[0 = \beta x_d + \gamma_d\delta^{d-1}|\xi|^{-1}k[-\gamma_d^{-1}|\xi|\Phi(x)+O(1)] = \beta x_d - \delta^{d-1}k[\Phi(x)+O(1)]. \]
Since we are in the unit cell we can bound the $O(1)$ term by $|x|^{2-d}$ up to a constant multiple depending on the unit cell diameter. Rearranging the previous displayed equation, and using that $|\beta| \geq 1/2$, then gives
\begin{equation}\label{e.wout-FB-location}
    |x_d| \leq C\delta^{d-1} |k|(C+\Phi(x)) \ \hbox{ for } \ x_d \in \partial \{\tilde{w}_{out}>0\} \cap (\cube^{\mathcal{Z}}  \setminus  B_r(0)).
\end{equation}

 On the same set we can use the expansion from \eref{removable-sing-omega} to write
    \[\grad \tilde{w}_{out}(x) = \beta e_d - \delta^{d-1}k\left[\grad \Phi(x) + \grad h(x)\right] \ \hbox{ in } \ \cube^{\mathcal{Z}}.\]
    Here $h$ is a harmonic function in $\cube^{\mathcal{Z}}_{3/2}$ which is even with respect to $x_d$ and so
    \begin{equation}\label{e.h-grad-bound}
        |\grad h(x) \cdot e_d| \leq C|x_d| \ \hbox{ in } \ \cube^{\mathcal{Z}}.
    \end{equation}
Now we compute the gradient along $\partial \{\tilde{w}_{out}>0\}$ via
\begin{equation}\label{e.wout-gradsquared}
    |\grad \tilde{w}_{out}(x)|^2 = \beta^2-2\beta\delta^{d-1}k\left[\grad \Phi(x) + \grad h(x)\right] \cdot e_d+O(k^2\delta^{2(d-1)}r^{2(1-d)}).
\end{equation}
    Using \eref{h-grad-bound} we can bound
    \[|\left[\grad \Phi(x) + \grad h(x)\right] \cdot e_d| \leq C|x_d|(\frac{1}{|x|^d}+1) \leq C\frac{|x_d|}{|x|^d} \ \hbox{ in } \ \cube^{\mathcal{Z}}.\]
    Due to \eref{wout-FB-location}, on $x \in \partial \{\tilde{w}_{out}>0\} \cap \cube^{\mathcal{Z}} \setminus B_r(0)$,
    \begin{align*}
        \left|2\beta\delta^{d-1}k\left[\grad \Phi(x) + \grad h(x)\right] \cdot e_d\right|  
        &\leq Ck^2\delta^{2(d-1)}r^{-d}\Phi(r).
    \end{align*}
    This is of the same order as the quadratic term in \eref{wout-gradsquared}, except in $d=2$ where it is larger by a factor of $\log (1/r)$.
    
    So, in the case $\sigma = -1$ and thus $\beta=(1-\alpha)$, using that $\beta^2\leq \beta$ for $\alpha \leq 1/2$,
    \[|\grad \tilde{w}_{out}(x)|^2 \leq 1-\alpha +Ck^2\delta^{2(d-1)}r^{-d}\Phi(r) \leq 1 \ \hbox{ for } \ x \in \partial \{\tilde{w}_{out}>0\} \cap (\cube^{\mathcal{Z}} \setminus B_r(0))\]
    where the right-most inequality is from \eqref{e.alpha}. By $\mathcal{Z}$-periodicity we derive the supersolution property on the entire $\R^d \setminus \cup_{z\in\mathcal{Z}} B_r(z)$. The argument in the case $\sigma = +1$ is even easier since the quadratic error term is positive and $(1+\alpha)^2 \geq 1+2\alpha$ for all $\alpha >0$.

\end{proof}

 Next we show that any {proper} solution of \eref{single-site} in $\R^d \setminus B_1$ with capacity $k$ can be patched together with an appropriate outer solution  to create a global periodic sub- or supersolution with the desired asymptotic slope.

\begin{lemma}\label{l.patched-barrier}
Suppose that $w_{in}$ is a solution of \eref{single-site} in $\R^d$, for example as constructed in \sref{pinned-solutions-2d} or \sref{monotone-family}, satisfying the expansion, in $d \geq 3$,
\[w_{in}(x) = x_d+s-k|x|^{2-d}+O(|x|^{1-d}) \ \hbox{ as } \ \overline{\{w_{in}>0\}}\ni x \to \infty,\]
or in $d=2$
\[w_{in}(x) = x_d+k\log |x|+O(1) \ \hbox{ as } \ \overline{\{w_{in}>0\}}\ni x \to \infty.\]
Then for any $\ep>0$ and $\delta>0$ sufficiently small depending on $\ep$ there is a supersolution (resp. subsolution) $w_{per}(x)$ of \eqref{e.model} which is $\mathcal{Z}$-periodic with asymptotic slope at least $1+\gamma_d|\xi|^{-1}(k-\ep)\delta^{d-1}$ (resp. at most $1+\gamma_d|\xi|^{-1}(k+\ep)\delta^{d-1}$).
\end{lemma}

\begin{proof}
 We will only present the construction of the supersolution $w_{per}$, the subsolution construction is parallel. 
 
 Let $0 < \ep \leq 1$. For the outer solution, we use a translated version of $\tilde{w}_{out}$ from \eqref{e.model-outer-soln} with $\sigma=-1$:
\begin{equation}
\begin{split}
    w_{out}(x):= &(1- \alpha )(x_d + \delta (s+s_{\Delta}))\\
    &+\gamma_d|\xi|^{-1}\delta^{d-1}(k-\ep)[x_d+\delta(s+s_\Delta)-\omega(x+\delta (s+s_\Delta) e_d)].
\end{split}
\end{equation}
Here the constants $0<\alpha<1$, $s_{\Delta} \in \R$, and also $s$ (only free in $d=2$), are to be chosen later in the proof.

Next for $\Lambda \geq 2$, another constant to be chosen later, define 
    \[w_{per}(x):= \begin{cases}
        \delta w_{in}(\delta^{-1} (x-z)) & |x-z| \leq \frac{1}{\Lambda}r \hbox{ for some $z \in \mathcal{Z}$,}\\
        \min\{\delta w_{in}(\delta^{-1} (x-z)),w_{out}(x)\} & \frac{1}{\Lambda}r \leq |x-z| \leq \Lambda r \hbox{ for some $z \in \mathcal{Z}$},\\
        w_{out}(x) & \textup{dist}(x,\mathcal{Z}) \geq \Lambda r.
    \end{cases}\]
   Here $r=r(\Lambda)$ is chosen sufficiently small so that $\delta \ll r/\Lambda $ and $\Lambda r \ll 1$ so that $B_{\Lambda r} \subset \cube^{\mathcal{Z}}_{1/2}$. We now choose 
\begin{equation}\label{e.barrier-alpha-choice}
\alpha := C_0k^2\Lambda^{2(d-1)}(\delta/r)^{2(d-1)}(r/\Lambda)^{d-2}\Phi(r/\Lambda),
\end{equation}
i.e. \eqref{e.alpha} with $r/\Lambda$ in place of $r$, to enable Lemma~\ref{l.outer-solution-alpha-info} in the region $\textup{dist}(x,\mathcal{Z}) \geq \Lambda^{-1}r$. Note that $\alpha \leq C(k)\Lambda^{4}\delta^2\log\Lambda$ in $d=2$ (with the choice $r = c/\Lambda$ made below) and $\alpha \leq C(k)\Lambda^{3d-4}\delta^{2(d-1)}r^{2-d}$ in $d \geq 3$. In either case $\alpha = o(\delta^{d-1})$ as $\delta \to 0$ with $\Lambda, r$ fixed.

To show  $w_{per}$ is a supersolution, the goal is to obtain
    \begin{equation}\label{e.outer-crossing-ineq}
        w_{out}(x)<\delta w_{in}(\delta^{-1} x) \ \hbox{ on } \ x \in \partial B_{\Lambda r}(0),
    \end{equation}
    and
    \begin{equation}\label{e.inner-crossing-ineq}
        w_{out}(x)>\delta w_{in}(\delta^{-1} x) \ \hbox{ on } \ x \in \partial B_{\frac{1}{\Lambda}r}(0).
    \end{equation}
\eqref{e.outer-crossing-ineq}-\eqref{e.inner-crossing-ineq} yields that  $w_{per}$
    is a  supersolution in $\Box_{\mathcal{Z}} \times \R$, and we can conclude using the $\mathcal{Z}$-periodicity of $w_{out}$.

   For $d=2$, for patching purposes, we compare the expansions of inner and outer solutions on $|x| = \rho$:
    \[\delta w_{in}(\delta^{-1} x) = x_d +k\delta \log \frac{\rho}{\delta}+O(\delta)\]
    and
  \[w_{out}(x) = (1-\alpha)\left[x_d+\delta(s+s_{\Delta})\right] +\delta(k-\ep)[\log \rho+O(1)].\]
  So, choosing $s := k\log\frac{1}{\delta}$ in the second expansion,
     \begin{align*}
         w_{out}(x)-\delta w_{in}(\delta^{-1} x) &= \delta \left[-\ep\log \rho+s_{\Delta}+s-k\log\frac{1}{\delta}\right]+O(\alpha \rho+\alpha\delta|s+s_\Delta|+\delta)\\
         &=\delta \left[-\ep\log \rho+s_{\Delta}\right]+O(\alpha \rho+\delta)\\
         &=: A(\rho) + E(\rho).
    \end{align*}
        So if we can choose parameters so that \begin{equation}\label{conditions}
    A(\Lambda r)<0, \,A(r/\Lambda)>0, \ \hbox{ and }\ |E(\rho)| < |A(\rho)| \ \hbox{ for } \ \rho \in \{ r/\Lambda,\Lambda r\},
    \end{equation} then we will achieve \eref{outer-crossing-ineq} and \eref{inner-crossing-ineq}.

    Let us now choose 
    \begin{equation}\label{e.parameter_1-2d}
   s_{\Delta} := \ep \log r \ \hbox{ and } \  r := c\Lambda^{-1}.
   \end{equation}
   With these choices
    \[A(\rho) = \delta \ep \log \frac{r}{\rho}, \  A(r/\Lambda) = \delta \ep\log \Lambda>0 \ \hbox{ and } \ A(\Lambda r) = - \delta \ep\log \Lambda<0.\]

   Moreover, $\alpha  \leq C\Lambda^{4}\delta^{2}\log \Lambda$ and, since $\rho \leq \Lambda r = c$,  \[ |E(\rho)| \leq C( \alpha \rho + \delta) \leq C_1(\Lambda^4\delta^{2}\log\Lambda+\delta).
   \]
   So for any $\ep>0$ we can choose $\Lambda$ sufficiently large so that $\ep \log \Lambda \geq 2C_1$ and then for $\delta$ sufficiently small $|E(\rho)| \leq |A(\rho)|$ for $\rho \in \{r/\Lambda,\Lambda r\}$. 
   
   Lastly, in both $d=2$ and the case $d \geq 3$ treated next, recall from \eqref{e.removable-sing-omega} that the asymptotic slope of this supersolution will be
    \begin{equation}\label{e.barrier-asymptotic-slope-err}
        \partial_{x_d}w_{per}(x) \to 1 -\alpha+\gamma_d|\xi|^{-1}(k-\ep)\delta^{d-1} \ \hbox{ as } \ x_d \to +\infty,
    \end{equation}
    and since $\alpha = o(\delta^{d-1})$ as $\delta \to 0$ by the choice of $\alpha$ above in \eqref{e.barrier-alpha-choice}, we can choose $\delta$ smaller if necessary so that the asymptotic slope is at least $1 +\gamma_d|\xi|^{-1}(k-2\ep)\delta^{d-1}$.

For $d \geq 3$, we take $\Lambda = 2$, and as before compare the expansions on $|x|=\rho$: 
    \[\delta w_{in}(\delta^{-1} x) = x_d + \delta s-k\delta^{d-1}\rho^{2-d} + O(\delta^{d} \rho^{1-d})\]
    and
    \[w_{out}(x) = (1-\alpha)\left[x_d+\delta(s+s_{\Delta})\right] +\delta^{d-1}(k-\ep)[-\rho^{2-d}+O(1+\delta \rho^{1-d})].\]
    So, still evaluating on $|x| = \rho$,
    \begin{align*}
         w_{out}(x)-\delta w_{in}(\delta^{-1} x) &= \delta \left[\ep(\delta/\rho)^{d-2}+s_{\Delta}\right]+O(\alpha \rho + \delta^{d-1}+\delta^d \rho^{1-d})\\
         &=: A(\rho) + E(\rho).
    \end{align*}

    As in the $d=2$ case it will suffice to show \begin{equation}\label{e.conditions-barrier-3d}
    A(2r)<0, \,A(r/2)>0, \ \hbox{ and }\ |E(\rho) | < |A(\rho)| \ \hbox{ for } \ \rho \in \{r/2,2r\},
    \end{equation} then we will achieve \eref{outer-crossing-ineq} and \eref{inner-crossing-ineq}.

  Let us now choose $s_{\Delta} := -\ep (\delta/r)^{d-2}$ so that
    \[A(\rho) = \delta^{d-1} \ep \left[\tfrac{1}{\rho^{d-2}}-\tfrac{1}{r^{d-2}}\right] \]
    and 
    \[A(r/2) = (2^{d-2}-1)\ep\delta^{d-1}r^{2-d}>0 \ \hbox{ and } \ A(2r) = -(1-2^{2-d})\ep\delta^{d-1}r^{2-d}<0.\]
 Indeed we have
    \[\min\{|A(r/2)| , |A(2r)|\} \geq c_1 \ep\delta^{d-1}r^{2-d}\]
   and the error term bound
  \begin{equation}\label{e.barrier-3d-E-errors}
        \begin{split}
       |E| &\leq C( \delta^{d-1}+\alpha r + \delta^dr^{1-d}) \\
       &\leq C(\delta^{d-1}+(\delta/r)^{2(d-1)}r+\delta^dr^{1-d})\\
       & = C_2\delta^{d-1}(1+\delta^{d-1}r^{3-2d}+\delta r^{1-d}).
     \end{split}
  \end{equation}
   Now choose $r$ sufficiently small so that $c_1\ep r^{2-d} \geq 2C_2 $, and then choose $\delta$ sufficiently small so that $C_2(\delta^{d-1}r^{3-2d}+\delta r^{1-d}) \leq C_2$. We then achieve \eref{conditions-barrier-3d} to conclude.

\end{proof}

\begin{proof}[Proof of \pref{asymptotic-exp-lower-bound}]
By \tref{monotone-family2} (in $d=2$) or \tref{monotone-family} (in $d \geq 3$), for any $0 < K <k_{\textup{adv}}(e)$ there is a proper solution $w_{in}$ of \eref{single-site} with capacity $k \geq K$.  \lref{patched-barrier} in turn delivers, for all sufficiently small $\delta>0$, a $\mathcal{Z}$-periodic supersolution $w_{per}(x)$ of \eqref{e.model} with asymptotic slope at least $1+\gamma_d|\xi|^{-1}K\delta^{d-1}$. Then, by the definition of $Q^\delta_{\textup{adv}}(e)$ in \dref{pinning-interval}, we have
\[Q^\delta_{\textup{adv}}(e) \geq 1+\gamma_d|\xi|^{-1}K\delta^{d-1}.\]
   Rearranging, taking first the limit $\delta \to 0$, and then $K\nearrow k_{\textup{adv}}(e)$ gives the desired bound for $Q^{\delta}_{\textup{adv}}(e)$. The bound of  $Q^\delta_{\textup{rec}}(e)$ is proved symmetrically with the subsolution counterpart of $\omega_{per}$.    
\end{proof}

\section{Asymptotic upper bound}\label{s.homogenization-upper-bound}
Now we complete the proof of \tref{asymptotics} by showing the upper bound part \eref{upper_bound} of the asymptotic expansion.

\begin{proposition}
    \label{p.asymptotic-exp-upper-bound}
Let $\xi \in \Z^d$ irreducible and $e = \xi/|\xi|$ then 
    \[
\limsup_{\delta \to 0} \frac{Q^\delta_{\textup{adv}}(e)-1}{\delta^{d-1}} \leq \gamma_d|\xi|^{-1}k_{\textup{adv}}(e) \ \hbox{ and } \ \limsup_{\delta \to 0} \frac{1-Q^\delta_{\textup{rec}}(e)}{\delta^{d-1}} \leq -\gamma_d|\xi|^{-1}k_{\textup{rec}}(e).
\]
\end{proposition}

Along the route to prove this Proposition we will also prove \tref{udelta-expansion}.

\subsection{Set-up and normalization}\label{s.UB-section-hypotheses} Throughout this section  $u_\delta$ will denote a strong Birkhoff plane-like solution of \eref{model} with asymptotic slope $\alpha=\alpha_\delta=1+\mu_\delta \in [Q^\delta_{\textup{rec}},Q^\delta_{\textup{adv}}]$. The existence of such solutions is guaranteed by \tref{cell-solutions-reference} part \ref{part.cell-solutions-reference-2} with $Q(x)=Q_{\delta}(x) = 1+\sum_{z \in \mathcal{Z}} q(\frac{x-z}{\delta})$.  In particular, due to \eref{plane-like-defn}, there is a constant $C(d,q)$ independent of $\delta$ and $\xi$ such that \begin{equation}\label{e.udelta-C-flat}
    (1+\mu_\delta)(x_d - C)_+ \leq u_\delta(x) \leq (1+\mu_\delta)(x_d+C)_+.
\end{equation}
We will denote
\[
 \Omega_{\delta}:= \{u_{\delta}>0\}.
 \]

    \begin{remark}
        For most of this section we do \emph{not} need to assume that $1+\mu_{\delta}$ is one of the extremal slopes, however  we will assume that $\mu_\delta \neq 0$. As we will see, this is sufficient to prove \pref{asymptotic-exp-upper-bound} and \tref{udelta-expansion}.  
    \end{remark}

Define
\begin{equation}\label{e.sdelta-pm-def}
\begin{split}
      &s_\delta^+ := \inf \{ s: u_\delta(x) \leq (1+\mu_\delta)(x_d+s)_+\} \ \hbox{ and } \\
      &s_\delta^- := \sup \{ s: (1+\mu_\delta)(x_d+s)_+ \leq u_\delta(x)\}.  
\end{split}
\end{equation}
Since $u_\delta$ are plane-like in the sense of \eref{plane-like-defn}, $s^\pm_{\delta}$ are well-defined finite values, bounded independent of $\delta$. One can also check that $(1+\mu_\delta)(x_d+s_\delta^\pm)_+$ must touch $u_\delta$ from above/below at the free boundary $\partial\Omega_{\delta}$: by maximum principle and periodicity it remains to rule out infimum occurring only at $x_d=\infty$, which would contradict the definition of $s^{\pm}_{\delta}$ by Harnack inequality. 

\begin{lemma}\label{l.normalization-sdelta-deltabd}
    Let $u_{\delta}, \mu_{\delta}$ and $s^{\pm}_{\delta}$ as given above. If $\mu_\delta> 0$ (resp. $\mu_\delta < 0$) then, up to a lattice translation, 
    \begin{equation}\label{e.udelta-translation-normalization}
        -\delta\leq s_\delta^- \ \hbox{(resp. $s_\delta^+ \leq \delta$) and } \  \partial \Omega_{\delta}\cap \overline{B_\delta(0)} \neq \emptyset. 
    \end{equation}
\end{lemma}

\begin{proof}
   We will prove for $\mu_\delta>0$: the other is symmetric. Since $\mu_\delta>0$ we have $(x_d+s_\delta^-)_+ \leq (1+\mu_\delta)(x_d+s_\delta^-)_+$, so $(x_d+s_\delta^-)_+$ touches $u_\delta$ from below at the same free boundary point where $(1+\mu_\delta)(x_d+s_\delta^-)_+$ touches from below. By a lattice translation we can assume that this touching point is in $[-1/2,1/2)^d$. If the touching point is not in $\overline{B_\delta(0)}$ then the strong maximum principle \lref{strong-max}\ref{part.strongmax-super} implies that $u_\delta \equiv (x_d+s_\delta^-)_+$ , contradicting \eref{udelta-C-flat}. 
\end{proof}
 \begin{remark}
   Since we have assumed that $\mu_\delta \neq 0$ the hypothesis of \lref{normalization-sdelta-deltabd} is satisfied. So, for the remainder of the section we will assume that $u_\delta$ is properly translated so that \eref{udelta-translation-normalization} holds.
   \end{remark}
\subsection{Initial $o_\delta(1)$-flatness}\label{s.homogenization-initial-flatness}    Our first result constrains the defect sites intersecting $\partial \{u_\delta >0\}$ to only those  $o_{\delta}(1)$-close to  $\{x_d = 0\}$ -- recalling that we have normalized by a lattice translation so that \eref{udelta-translation-normalization} holds.

\begin{lemma}\label{l.initial-flatness}
    Let $u_\delta$ be a sequence of strong Birkhoff solutions satisfying the hypotheses described in \sref{UB-section-hypotheses}. Then $\mu_{\delta} \to 0$ and $s_\delta^\pm \to 0$ as $\delta \to 0$.
\end{lemma}

\begin{proof}
   Due to \lref{lipschitz-estimate}, $u_\delta$ are uniformly Lipschitz continuous.   Every subsequence of $\delta \to 0$ has a further subsequence (not re-labeled) so that $u_\delta$ converges locally uniformly to some $\mathcal{Z}$-periodic $w$, $\mu_\delta \to \mu_0$ as $\delta\to 0$, and also $\mu_\delta$ has a fixed sign along this subsequence. To conclude, it is enough to show that $\mu_{\delta}, s^{\pm}_{\delta}\to 0$ along every such subsequence. Of course $w$ is harmonic in its positivity set and we also claim that:  
\begin{equation}\label{e.bernoulli-initial-flatness}
 w \hbox{ is a supersolution of the homogeneous problem }\eqref{e.bernoulli}.
\end{equation} 
  Let us postpone the proof of \eqref{e.bernoulli-initial-flatness} for now and complete the proof of the Lemma. 

We will argue for the case $\mu_\delta>0$ and $\mu_0 \geq 0$: the other case is similar. Taking the limit of \eref{udelta-C-flat} and \eref{udelta-translation-normalization} we find $w(0) = 0$ and 
\begin{equation}\label{e.w-flat-C}
      (1+\mu_0)(x_d)_+ \leq w(x) \leq (1+\mu_0)(x_d+C)_+ \ \hbox{ and } \ 0 \in \partial \{w>0\}.
\end{equation}
If $\mu_0>0$ then, by \eref{w-flat-C}, the strict subsolution $(1+\mu_0)(x_d)_+$ touches $w$ from below at $0$, contradicting \eref{bernoulli-initial-flatness}.  Thus $\mu_0 = 0$. Then applying the strong maximum principle \lref{strong-max}\ref{part.strongmax-super} we conclude that $w(x) \equiv (x_d)_+$. Since $\mu_\delta>0$ we already know that $s_\delta^- \to 0$, thus it remains to show that $s_{\delta}^+\to 0$.  Recall that $(1+\mu_{\delta})(x_d+s_\delta^+)_+$ touches $u_\delta$ from above at some $x_\delta \in \partial \Omega_{\delta} \cap B_{C}(0)$.  Then $x_\delta$ is an outer-regular point, so non-degeneracy at outer regular free boundary points (\lref{linear-growth-outer-reg}), and local uniform convergence of $u_\delta(x) \to (x_d)_+$ implies that $s_\delta^+ \to 0$.  

We conclude that $\mu_\delta \to 0$ and $s_\delta^\pm \to 0$ along the full sequence.

Now we proceed to prove \eqref{e.bernoulli-initial-flatness}: we only need to check the free boundary condition. Consider a test function $\varphi$ touching $w$ strictly from below at $x_0 \in \partial \{w>0\}$ with $\Delta \varphi(x_0) > 0$. The only interesting case is when $x_0 \in \mathcal{Z}$, so we may assume that $x_0=0$. Then there is a sequence $x_\delta \to 0$ and $c_\delta \to 0$ so that $\varphi(x) + c_\delta$ touches $u_\delta$ from below at $x_\delta$. Since $\Delta \varphi(x_\delta) >0$, $x_\delta\in \partial \Omega_{\delta}$ and  $c_\delta = - \varphi(x_\delta)$. We can conclude now except when $x_\delta \in B_\delta(0)$ for all $\delta>0$.  In this case consider the blow-up sequences
    \[\varphi_\delta(y) := \frac{1}{\delta}(\varphi(x_\delta + \delta y) - \varphi(x_\delta)) \ \hbox{ and } \ \tilde{w}_\delta(y) := \frac{1}{\delta} u_\delta (x_\delta + \delta y)\]
    with $\varphi_\delta$ touching $\tilde{w}_\delta$ from below at $0 \in \partial \{\tilde{w}_\delta>0\}$.  As $\delta \to 0$, since $\varphi$ is $C^1$,
    \[\varphi_\delta(y) \to \grad \varphi(0) \cdot y.\]
   Along a further subsequence, we may assume that $x_\delta/\delta \to \tau \in \overline{B_1(0)}$ and $\tilde{w}_\delta$ converges locally uniformly to $\tilde{w}$ which, by \tref{stability}, is a relaxed solution of the single-site problem
    \[ \Delta \tilde{w} = 0 \ \hbox{ in } \ \{\tilde{w}>0\} \ \hbox{ with } \ |\grad \tilde{w}| = Q(\tau + y).\]
    Finally, since $\varphi_\delta \leq \tilde{w}_\delta$, passing to the limit we find that
    \[\tilde{w}(y) \geq \grad \varphi(0) \cdot y .\]
    If $|\grad \varphi(x_0)| \leq 1$ we are done, while if $|\grad \varphi(x_0)| \geq 1$ then $\tilde{w}$ is one-sided flat in the sense of \eref{below-bound}:
    \[\tilde{w}(y) \geq (\grad \varphi(x_0) \cdot y)_+ \geq (e \cdot y)_+ \ \hbox{ with } \ e :=\frac{\grad \varphi(x_0)}{|\grad \varphi(x_0)|}.\]
    Hence \pref{blow-down} applies and we know the blow-down is of the form
    \[ \lim_{r \to \infty}\frac{1}{r}\tilde{w}(ry) =(e \cdot y)_+ \ \hbox{ locally uniformly in } \ y \cdot e \geq 0.\]
    but
    \[(\grad \varphi(x_0) \cdot y)_+ \leq \lim_{r \to \infty}\frac{1}{r}\tilde{w}(ry)  = (e \cdot y)_+  \ \hbox{ in } \ y \cdot e \geq 0\]
    and so $|\grad \varphi(x_0)| \leq 1$.
    
    The proof in the case $\mu_\delta <0$ is similar. Instead of \eref{bernoulli-initial-flatness} we argue that $(w,E)$ with $E:={\limsup}^*\{w_\delta>0\}$ is a {relaxed subsolution}. Taking the limit of \eref{udelta-C-flat} and \eref{udelta-translation-normalization} gives 
    \begin{equation*}
      (1+\mu_0)(x_d-C)_+ \leq w(x) \leq (1+\mu_0)(x_d)_+, \ E \subset \{x_d \geq 0\}, \ \hbox{ and } \ 0 \in \partial E.
\end{equation*}
and we can apply \lref{strong-max}\ref{part.strongmax-sub}.  Note that $0 \in \partial E$ since it is a limit of points of $\partial E$ with an exterior touching half-space so \lref{linear-growth-outer-reg} applies at those points.
  
\end{proof}
\begin{remark}
From \tref{cell-solutions-reference} part \ref{part.cell-solutions-reference-2} we can derive at this point:
    \[\lim_{\delta \to 0} Q^\delta_{\textup{adv}} = \lim_{\delta \to 0}Q^\delta_{\textup{rec}} =1. \]
\end{remark}

\subsection{Refined flatness and fundamental solution bounds in $d =2$}Next we improve the flatness in Lemma~\ref{l.initial-flatness} to the optimal scale $O(|\mu_\delta|+\delta)$ and establish the asymptotic behavior of $u_{\delta}$ near the holes. The proof for $d= 2$ will be presented here, and the next section will address $d \geq 3$.
\begin{proposition}\label{p.2d-capacitory-bounds}
    Let $d=2$ and let $u_\delta$ and $\mu_{\delta}$ be as given in \sref{UB-section-hypotheses}. Then let $s_\delta = s_\delta^\pm$ when $\pm \mu_\delta>0$ and if we define
    \[\omega_\delta(x) := \frac{1}{\delta}[(1+\mu_\delta)(x_d+s_\delta)-u_\delta(x)],\] then in the case $\mu_\delta > 0$
    \[ 0\leq {\omega}_{\delta}(x) \leq \kappa^{\frac{1}{2\delta}}_{\textup{adv}}  (\log \tfrac{1}{|x|})_+ +C(1+\tfrac{\mu_\delta}{\delta}) \quad\hbox{ for } x \in (\Box_{\mathcal{Z}} \times \R) \cap \Omega_{\delta},\]
    and, in the case $\mu_\delta < 0$ ,
    \[0 \geq {\omega}_{\delta}(x) \geq \kappa^{\frac{1}{2\delta}}_{\textup{rec}} (\log \tfrac{1}{|x|})_+-C(1+\tfrac{|\mu_\delta|}{\delta}) \quad\hbox{ for } x \in (\Box_{\mathcal{Z}} \times \R) \cap \Omega_{\delta}.\]
    Recall  that $\kappa^R_{\textup{adv}/\textup{rec}}$ satisfies, by \tref{monotone-family2}, $\lim_{R \to \infty} \kappa^R_{\textup{adv}/\textup{rec}} = k_{\textup{adv}/\textup{rec}}$.
\end{proposition}

     We start with a Lemma to establish $O(|\mu_\delta| + \delta)$ flatness away from the defects.
\begin{lemma}\label{l.2d-harnack-flatness}
Under the hypotheses of \pref{2d-capacitory-bounds}, if $\mu_\delta > 0$ then
\begin{equation}\label{e.wdelta-harnack-bounds-2d}
 0\leq    (1+\mu_\delta)(x_d + s_\delta^+)-u_\delta(x) \leq C(\delta+\mu_\delta) \ \hbox{ on } \ \Omega_\delta  \setminus \cup_{z \in \mathcal{Z}} B_{1/2}(z),
\end{equation}
where $C>1$ is some constant. A symmetrical result holds if $\mu_\delta<0$ with $s_\delta^-$.
\end{lemma}
 
\begin{proof}
 Let's assume that $\mu_\delta>0$, the case $\mu_\delta < 0$ is analogous. By \lref{normalization-sdelta-deltabd},
     \begin{equation}\label{e.wdelta-initial-flatness-bounds-2d}
          (1+\mu_\delta)(x_d-\delta) \leq u_\delta(x) \leq (1+\mu_\delta)(x_d+s_\delta^+),
     \end{equation}
 which, in particular, yields the lower bound in \eqref{e.wdelta-harnack-bounds-2d}. By the discussion from \sref{UB-section-hypotheses} (see the paragraph above \lref{normalization-sdelta-deltabd}) there is $x^0 \in \partial \Omega_\delta \cap \cube^{\mathcal{Z}}$ (see notation defined in \eqref{e.cube-notation-def}) with $x^0_d+s_\delta^+ = 0$. For the remainder of the proof we denote $s_\delta = s_\delta^+$.

  Recall $\cube^{\mathcal{Z}}_r:=r\cube^{\mathcal{Z}}$ and define
    \begin{equation}\label{e.sdelta-defn-2d}
        S := \max_{\overline{\cube^{\mathcal{Z}}\cap \Omega_\delta}} (u_\delta(x) - x_d).
    \end{equation}
     Since $x^0 \in \partial \Omega_\delta \cap \cube^\mathcal{Z}$ and $u_\delta(x^0) = 0$, 
    \begin{equation}\label{e.S-LB-sdelta}
        S \geq  u_\delta(x^0)- x^0_d = - x^0_d = s_\delta.
    \end{equation}

        Note that $(x_d + S)_+$ touches $u_\delta$ from above in $\cube^{\mathcal{Z}}\cap \Omega_\delta$ at $x^1$, where the maximum in \eqref{e.sdelta-defn-2d} is achieved.
            Since $u_\delta(x) - x_d$ is harmonic in $\Omega_\delta \cap \cube^{\mathcal{Z}}$,  either (i) $x^1 \in \cube^{\mathcal{Z}} \cap \partial \Omega_\delta$ or (ii) $x^1 \in \Omega_{\delta} \cap \partial \cube^{\mathcal{Z}}$.

               First suppose (i). Then $(x_d+S)_+$ touches $u_\delta$ from above at $x^1$ in an open neighborhood. If $x^1 \not\in \overline{B_{\delta}(0)}$ then \lref{strong-max}\ref{part.strongmax-sub} implies $u_\delta(x) \equiv (x_d+S)_+$ which contradicts $\mu_\delta>0$.  If $x^1 \in B_\delta(0) \cap \partial \Omega_\delta$ then $S \leq \delta$. Due to \eref{S-LB-sdelta} $s_{\delta} = -x^0_d \leq \delta$. Since \eqref{e.wdelta-initial-flatness-bounds-2d} yields $ (1+\mu_{\delta})x_d -u_{\delta}(x) \leq C\delta$ in $\Omega_{\delta}$, we conclude  \eqref{e.wdelta-harnack-bounds-2d} in case of (i).

    It remains now to consider the case (ii). We will aim to apply Harnack inequality, \cite{FKexterior}*{Corollary 2.7}, in an annular domain $\cube^{\mathcal{Z}}_{3/2} \setminus B_{1/4}$. Note  
    \[ \max_{\overline{\cube^{\mathcal{Z}}_{3/2}\cap \Omega_\delta}}(u_\delta(x) - x_d) \leq \max_{\overline{\cube^{\mathcal{Z}}_{3/2}\cap \Omega_\delta}} (u_\delta(x) - (1+\mu_\delta)x_d)+\mu_\delta \max_{\overline{\cube^{\mathcal{Z}}_{3/2}\cap \Omega_\delta}} x_d \leq s_\delta+C_0\mu_\delta,\]
     So
    \[(x_d - \delta)_+ \leq u_\delta(x) \leq (x_d + s_\delta + C_0\mu_\delta)_+ \ \hbox{ in } \ \Omega_\delta \cap \cube^{\mathcal{Z}}_{3/2}.\]
     Due to \lref{initial-flatness}, $s_\delta + C\mu_\delta \to 0$ as $\delta \to 0$, and thus the Harnack inequality \cite{FKexterior}*{Corollary 2.7} applies for sufficiently small $\delta$. Call $A := \cube^{\mathcal{Z}}_{1} \setminus B_{1/2}$. Since $\cube^{\mathcal{Z}}_{3/2} \setminus B_{1/4}\supset \supset A$, \cite{FKexterior}*{Corollary 2.7}   implies that
   \begin{equation}\label{e.2d-bounds-harnack-1}
       \max_{\Omega_\delta \cap A} (x_d+s_\delta + C_0\mu_\delta - u_\delta(x)) \leq C\min_{\Omega_\delta \cap A}(x_d+s_\delta + C_0\mu_\delta - u_\delta(x)).
   \end{equation}
    Utilizing (ii) to plug in $x^1\in\partial \cube^{\mathcal{Z}}_1 \subset A$  on the right hand side yields
    \[\min_{\Omega_\delta \cap A}(x_d+s_\delta + C_0\mu_\delta - u_\delta(x)) \leq s_\delta + C_0\mu_\delta - S \leq C_0\mu_\delta,\]
    using \eref{S-LB-sdelta} for the last inequality.  Plugging this back into \eref{2d-bounds-harnack-1} we get
    \begin{equation}\label{e.2d-A-lower-bound-result}
        \max_{\Omega_\delta \cap A}(x_d+s_\delta - u_\delta(x)) \leq C\mu_\delta.
    \end{equation}
   Lastly we look for the upper bound in $\{x_d \geq 1\}$. Note that \lref{initial-flatness} ensures $\{x_d \geq 1\} \subset \Omega_\delta$ for $\delta\ll 1$. Since $\Box_{\mathcal{Z}} \times \{1\} \subset A$, using $\mathcal{Z}$-periodicity of $u_\delta$,
    \begin{align*}
        \max_{\{x_d = 1\}}((1+\mu_\delta)(x_d+s_\delta)-u_\delta(x) ) &= \max_{\Box_{\mathcal{Z}} \times \{1\}}((1+\mu_\delta)(x_d+s_\delta)-u_\delta(x)) \\
        &\leq \max_{\Omega_\delta \cap A}((1+\mu_\delta)(x_d+s_\delta)-u_\delta(x)) \\
        &\leq  C\mu_\delta,
    \end{align*}
 where the last line is due to \eqref{e.2d-A-lower-bound-result}.   Since $(1+\mu_\delta)(x_d+s_\delta)-u_\delta(x)$ is bounded harmonic function in $\{x_d >1\}$ , the maximum principle in this domain implies that
    \[\max_{\{x_d \geq 1\}} ((1+\mu_\delta)(x_d+s_\delta)-u_\delta(x)) \leq  C\mu_\delta. \]
    Combining this with \eref{2d-A-lower-bound-result} gives the upper bound in $\Omega_\delta \setminus \cup_{z \in \mathcal{Z}} B_{1/2}(z)$.
\end{proof}
    \begin{proof}[Proof of \pref{2d-capacitory-bounds}]
    We prove the case $\mu_\delta > 0$: the other case is symmetrical. By \lref{2d-harnack-flatness},  $s_\delta = s_\delta^+$ satisfies      
      \[u_\delta(x) \geq x_d + s_\delta - C_0(\mu_\delta+\delta) \ \hbox{ on } \ \partial B_{1/2}(0) \cap \Omega_\delta.\]
    Define
 \begin{equation}\label{e.Rdelta-kdelta-def}
     R_\delta := \frac{1}{2\delta} \ \hbox{ and } \ k_\delta := \frac{1}{\delta \log \frac{1}{2\delta}}[s_\delta -  C_0(\mu_\delta+\delta)],
 \end{equation}
   and
   \[\varphi(x):= \delta u^{k_\delta,R_\delta}_{\textup{adv}}(\tfrac{x}{\delta}),\]
   where $u^{k,R}_{\textup{adv}}$ is from \dref{adv-2d-R}. Due to the definition of $u^{k,R}_{adv}$ and $k_{\delta}$,  
   \[  \varphi(x) = (x_d + \delta k_\delta \log \tfrac{|x|}{\delta})_+ = (x_d+s_\delta - C_0(\mu_\delta+\delta))_+ \ \hbox{ on } \ \partial B_{1/2}(0).\]
   Since $\varphi$ is the minimal supersolution of \eqref{e.general-Q-bernoulli} in $B_{1/2}(0)$ with this boundary data, $\varphi(x) \leq u_\delta(x)$ in $B_{1/2}(0)$. Since $\{u_\delta =0\} \cap \overline{B_\delta(0)} \neq \emptyset$ and so, since $\varphi \leq u_\delta$, also $\{\varphi = 0\} \cap \overline{B_\delta(0)} \neq \emptyset$.  Therefore, by the definition of $\kappa^{R_\delta}_{\textup{adv}}$ in \dref{adv-2d-R},
\begin{equation}\label{e.kdelta-UB}
    k_\delta \leq \kappa_{\textup{adv}}^{R_\delta}.
\end{equation}
   By \lref{capacitory-lower-bound}
   \[ \varphi(x) \geq x_d + k_\delta\delta\log \tfrac{|x|}{\delta} - C_1\delta \ \hbox{ in } \ \{\varphi>0\}.\]
   So we obtain the desired upper bound in $B_{1/2}(0) \cap \{\varphi >0\}$:
   \begin{align*}
       \omega_\delta(x)  &\leq \frac{1}{\delta}[(1+\mu_\delta)(x_d+s_\delta)-\varphi(x)]\\
       &\leq \frac{1}{\delta}[(1+\mu_\delta) (x_d+s_\delta)-(x_d + k_\delta\delta\log |x|+k_\delta\delta \log \tfrac{1}{\delta} - C_1\delta)]\\
       & =\frac{1}{\delta}[ k_\delta\delta\log \tfrac{1}{|x|}+(s_\delta-k_\delta\delta\log \frac{1}{\delta}) + C_1\delta+\mu_\delta (x_d+s_\delta)]\\
       & =  \kappa_{\textup{adv}}^{R_\delta} \log \tfrac{1}{|x|}+C_0(1+\mu_\delta/\delta)+C_1+\delta^{-1}\mu_\delta(x_d+s_\delta).
   \end{align*}
   For the last equation we used \eref{kdelta-UB} to bound $k_\delta$ in the first term, and the definition of $k_\delta$ in \eref{Rdelta-kdelta-def} for the middle term in parenthesis.

   Combining this upper bound inside of $B_{1/2}(0)$ with periodicity and the bounds from \lref{2d-harnack-flatness} outside of $\cup_{z \in \mathcal{Z}}B_{1/2}(z)$, we obtain the result.
\end{proof}

\subsection{Refined flatness and fundamental solution bounds in $d \geq 3$}\label{section:compactness}

\begin{proposition}\label{p.refined-flatness-3d}
   Let $d\geq 3$, and let $\mu_{\delta}$ and $u_{\delta}$ as given in \sref{UB-section-hypotheses}. Then there is $s_\delta \to 0$ so that if we call
    \[\omega_\delta(x) := \frac{1}{\delta^{d-1}}[(1+\mu_\delta)(x_d+s_\delta)-u_\delta(x)],\] 
    then  for $x \in (\Box_{\mathcal{Z}} \times \R) \cap \Omega_{\delta}$,
    \[ k_{\textup{rec}} |x|^{2-d} -C(1+\delta|x|^{1-d}+\tfrac{|\mu_\delta|}{\delta^{d-1}})\leq \omega_{\delta}(x)\leq k_{\textup{adv}} |x|^{2-d}+C(1+\delta|x|^{1-d}+\tfrac{|\mu_\delta|}{\delta^{d-1}}).\]
\end{proposition}

We expect that the maximum and minimum heights $s_\delta^\pm$ defined in \eref{sdelta-pm-def} are at least $O(\delta)$, since this level of oscillation can be achieved just within the defect.  However for $d\geq 3$ we expect the oscillations of $u_\delta$ away from the defect to be of lower order, $\delta^{d-1}$.  This requires more precise choice of height $s_\delta$, in contrast to the $d=2$ case where it was sufficient to choose the height as one of $s^{\pm}_{\delta}$.

\begin{proof}
 
  The argument will proceed in three steps. First we obtain the bounds in $\cube^{\mathcal{Z}}$ (defined in \eref{cube-notation-def}), using the single-site barriers constructed in \tref{monotone-family}.  Next we extend these to global bounds by using the asymptotic growth of  $u_\delta(x) = (1+\mu_\delta)x_d + O(1)$ as $x_d \to +\infty$ and using the $\mathcal{Z}$-periodicity. Finally we use Harnack inequality and Lemma~\ref{l.initial-flatness} to 
  conclude.

   Let $u_{\textup{adv}/\textup{rec}}(x;s)$ be from \tref{monotone-family}. By \tref{monotone-family} and \tref{main-3d-capacities}, for all but countably many $s \in \R$,
   \begin{equation}\label{e.sdelta-bounds-expansion}
       |u_{\textup{adv}/\textup{rec}}(x;s) - (x_d + s - \kappa_{\textup{adv}/\textup{rec}}(s)|x|^{2-d})| \leq C_1|x|^{1-d}
   \end{equation}
   with $\sup_s|\kappa_{\textup{adv}/\textup{rec}}(s)|$ and $C_1$ having universal bounds.  Therefore we can choose a large universal constant $C_0$ so that
   \[ 2 \max_{x\in \partial \cube^{\mathcal{Z}}}\left[\max\{k_{\textup{adv}},|k_{\textup{rec}}|\}|x|^{2-d}+C_1\delta|x|^{1-d}\right] \leq C_0 \ \hbox{ for all } \ 0 < \delta < 1/2. \]
Define the candidates for $s_{\delta}$ as
\begin{equation}\label{e.sdelta-pm-defn}
\begin{split}
    t_\delta^+ &:= \max_{(\cube^{\mathcal{Z}} \setminus B_{1/2}) \cap \Omega_\delta}\{u_\delta(x) - x_d\}+C_0\delta^{d-1}, \ \hbox{ and } \\
    t_\delta^- &:= \min_{(\cube^{\mathcal{Z}} \setminus B_{1/2}) \cap \Omega_\delta}\{u_\delta(x) - x_d\}-C_0\delta^{d-1}.
\end{split}
\end{equation}

Now we focus on the bounds in $\cube^{\mathcal{Z}}$. We will argue only for the lower bound of $\omega_{\delta}$: the proof of the other bound is symmetric. Define the barrier 
    \[\varphi(x):=\delta u_{\textup{rec}}(\tfrac{x}{\delta};\tfrac{t_\delta^+}{\delta}) \ \hbox{ and } \ k_\delta:=\kappa_{\textup{rec}}(\tfrac{t_\delta^+}{\delta}).\]

    Evaluating on $\overline{\Omega_\delta} \cap \partial \cube^{\mathcal{Z}}$ and using  \eref{sdelta-bounds-expansion},
\begin{align*}
    \varphi(x) = \delta u_{\textup{rec}}(\tfrac{x}{\delta};\tfrac{t_\delta^+}{\delta})  &\geq  x_d+t_\delta^+ -k_\delta\delta^{d-1}|x|^{2-d}-C_1\delta^d|x|^{1-d}\\
    &\geq x_d +\max_{\partial \cube^{\mathcal{Z}} \cap \Omega_\delta}\{u_\delta(x) - x_d\}+C_0\delta^{d-1}-\frac{1}{2}C_0\delta^{d-1}\\
    &> x_d +\max_{\partial \cube^{\mathcal{Z}} \cap \Omega_\delta}\{u_\delta(x) - x_d\}\\
    &\geq u_\delta(x).
\end{align*}
    Thus $\varphi > u_\delta$ on $\overline{\Omega_\delta} \cap \partial \cube^{\mathcal{Z}}$. Since $\varphi$ is a maximal subsolution of \eref{model} in $\Box^\mathcal{Z}$, we conclude that  $u_\delta \leq \varphi$ in $\cube^{\mathcal{Z}}$ or
 \begin{equation}\label{e.wdelta-upper-inproof}
     u_\delta(x) \leq x_d+t_\delta^++\delta^{d-1}[-k_\delta|x|^{2-d}+C_0+C_1\delta|x|^{1-d}] \ \hbox{ in } \cube^{\mathcal{Z}}.
 \end{equation}
 Also by definition $k_\delta \geq \min_s \kappa_{\textup{rec}}(s) = k_{\textup{rec}}$.

Now we proceed to obtain global bounds. By the $\mathcal{Z}$-periodicity of $u_{\delta}$ and  \lref{periodic-bdry-layer},
   \begin{align*}
       \max_{\{x_d \geq 1\}} (u_\delta(x) - (1+\mu_\delta)x_d) &= \max_{\{x_d=1\}}(u_\delta(x) - (1+\mu_\delta)x_d) \\
       &= \sup_{\Box_{\mathcal{Z}} \times \{1\}}(u_\delta(x) - (1+\mu_\delta)x_d)\leq t_\delta^++(-\mu_\delta)_+,
   \end{align*}
 Symmetrical arguments for the minimum imply that
   \begin{equation}\label{e.wdelta-outer-region-bound}
       t_\delta^--(\mu_\delta)_+ \leq u_\delta(x) - (1+\mu_\delta)x_d \leq t_\delta^++(-\mu_\delta)_+ \ \hbox{ in } \{\ x_d \geq 1\}.
   \end{equation}

Thus, combining \eref{wdelta-upper-inproof} (in $x_d \leq 1$) and \eref{wdelta-outer-region-bound} (in $x_d \geq 1$) yields the global bounds
\begin{equation}\label{e.UB-wdelta-inproof}
    u_\delta(x) \leq (1+\mu_\delta)(x_d+t_\delta^+)+\delta^{d-1}[-k_{\textup{rec}}|x|^{2-d}+C+C\delta|x|^{1-d}]+C|\mu_\delta| \ \hbox{ in } \ \Omega_\delta \cap [\Box_{\mathcal{Z}} \times \R]
\end{equation}
and 
\begin{equation}\label{e.LB-wdelta-inproof}
    u_\delta(x) \geq (1+\mu_\delta)(x_d+t_\delta^-)+\delta^{d-1}[-k_{\textup{adv}}|x|^{2-d}-C-C\delta|x|^{1-d}]-C|\mu_\delta| \ \hbox{ in } \ \Omega_\delta \cap [\Box_{\mathcal{Z}} \times \R].
\end{equation}

   Finally, we aim to apply the Harnack inequality \cite{FKexterior}*{Corollary 2.7} to $u_\delta$ in the ``annulus" domain $A:= \cube^{\mathcal{Z}}_{3/2} \setminus \overline{B_{1/4}}$. Note that $u_{\delta}$ solves the homogeneous Bernoulli problem \eqref{e.bernoulli} in $A$ and satisfies, from \eref{LB-wdelta-inproof} above,
   \[u_{\delta}(x) \geq x_d+t_\delta^--C(\delta^{d-1}+|\mu_\delta|) \ \hbox{ in } \ \Omega_\delta \cap A.\]
    By \lref{initial-flatness}, \cite{FKexterior}*{Corollary 2.7} applies to $u_\delta$ in $A$ and we have
\begin{multline}
\sup_{\Omega_\delta \cap(\cube^{\mathcal{Z}} \setminus B_{1/2})} (u_{\delta}(x)-(x_d+t_\delta^--C(\delta^{d-1}+|\mu_\delta|))) \leq\\ C\inf_{\Omega_\delta \cap(\cube^{\mathcal{Z}} \setminus B_{1/2})}  (u_{\delta}(x)-(x_d+t_\delta^--C(\delta^{d-1}+|\mu_\delta|))).
\end{multline}
Plugging in the point on $\Omega_\delta \cap(\cube^{\mathcal{Z}} \setminus B_{1/2})$ where $u_{\delta}(x)-x_d+C_0\delta^{d-1} = t_\delta^+$ is achieved on the left, and the point where $u_{\delta}(x)-x_d-C_0\delta^{d-1} = t_\delta^-$ is achieved on the right gives
\[t_\delta^+-t_\delta^- \leq C(\delta^{d-1}+|\mu_\delta|). \]
Taking, for example, $s_\delta = t_\delta^-$ we conclude by the above estimate, \eref{UB-wdelta-inproof}, and \eref{LB-wdelta-inproof}.

Finally, since $(1+\mu_\delta)(x_d+s_\delta^-)_+ \le u_\delta \le (1+\mu_\delta)(x_d+s_\delta^+)_+$, on $(\Box^{\mathcal Z}\setminus B_{1/2})\cap\Omega_\delta$ we have $|u_\delta(x)-x_d| \le C(|s_\delta^+|+|s_\delta^-|+|\mu_\delta|)$, so $|s_\delta| = |t_\delta^-| \le C(|s_\delta^+|+|s_\delta^-|+|\mu_\delta|+\delta^{d-1}) \to 0$ by Lemma~\ref{l.initial-flatness}.
   
\end{proof}

\subsection{A rescaling limit leading to the perforated domain cell problem}\label{s.homogenization-capacity}

We are finally ready to prove \pref{asymptotic-exp-upper-bound} and  \tref{udelta-expansion} in this section. 

Let $\omega_{\delta}$ and $s_\delta$ be given from Propositions \ref{p.2d-capacitory-bounds} and \ref{p.refined-flatness-3d}. At the moment we do not yet have an upper bound  on $\mu_\delta/\delta^{d-1}$, thus we will first consider the limit of
\begin{equation}\label{e.tilde-omega-def}
    \tilde{\omega}_\delta(x): = \frac{\delta^{d-1}}{\mu_\delta} \omega_\delta(x) =  \frac{1}{\mu_\delta}[(1+\mu_\delta)(x_d+s_\delta)-u_\delta(x)].
\end{equation}
The main technical result of this section is on the precompactness of the $\tilde{\omega}_\delta$ and the limiting PDE boundary value problem.

\begin{lemma}\label{l.omega-neumann}
    Consider a sequence $\delta\to 0$ such that $\mu_\delta \neq 0$ and  $\frac{\mu_\delta}{\delta^{d-1}}\to \lambda \in [-\infty,\infty] \setminus \{0\}$. Then the $\tilde{\omega}_\delta$ are locally uniformly bounded and $C^{1,\alpha}$ in $\R^d_+ \setminus \mathcal{Z}$, and any subsequential limit of $\tilde{\omega}_\delta$,  $\omega$, is a $\mathcal{Z}$-periodic solution of
    \begin{equation}\label{e.tildeomega-limit-PDE}
        \begin{cases}
        \Delta \omega = 0 & \hbox{in } \R^d_+\\
        \partial_{y_d}\omega = 1  &\hbox{on } \partial \R^d_+ \setminus \mathcal{Z},\\
        \langle \partial_{y_d} \omega(\cdot,y_d) \rangle' = 0 &\hbox{for } y_d>0.
     \end{cases}
    \end{equation}
 Moreover when $\lambda>0$, with the notation $(\infty)^{-1}=0$,   
     \begin{equation}\label{e.tildeomega-limit-holes}
         \begin{cases}
             \displaystyle\limsup_{x \to z} \dfrac{\omega(x)}{\Phi(x-z)} \leq k_{\textup{adv}}\lambda^{-1} & \hbox{for } z \in \mathcal{Z},\\
             \displaystyle\liminf_{x \to z} \dfrac{\omega(x)}{\Phi(x-z)} \geq k_{\textup{rec}}\lambda^{-1} & \hbox{for } z \in \mathcal{Z}.
         \end{cases}
     \end{equation}
   When $\lambda<0$, \eqref{e.tildeomega-limit-holes} holds with $k_{\textup{adv}}$ and $k_{\textup{rec}}$ interchanged.
\end{lemma}

For clarity of presentation, we postpone the proof of the Lemma until  the end of the section to finish the proofs of \pref{asymptotic-exp-upper-bound} and \tref{udelta-expansion}.

\begin{proof}[Proof of \pref{asymptotic-exp-upper-bound} and \tref{udelta-expansion}]

We will only argue for $Q^{\delta}_{\textup{adv}}$: the other case is similar.
    
 First we prove \pref{asymptotic-exp-upper-bound}. Let $u_{\delta_j}$ be a sequence of strong Birkhoff solutions with slopes $1+\mu_{\delta_j}$ such that $\frac{\mu_{\delta_j}}{\delta_j^{d-1}} \to \lambda \in [0,+\infty]$. To prove \pref{asymptotic-exp-upper-bound} it is enough to show that \begin{equation}\label{e.upper_2}\lambda \leq \gamma_d|\xi|^{-1}k_{\textup{adv}}.\end{equation} 
 If $\lambda = 0$ we conclude since $k_{\textup{adv}} \geq 0$. If $k_{\textup{adv}} = +\infty$ the claim is also trivial, so we can assume that $\lambda>0$ and $k_{\textup{adv}}<+\infty$, and apply \lref{omega-neumann}. Let $\omega$ be as given in \lref{omega-neumann} and define
    \[A := k_{\textup{rec}}\lambda^{-1} \ \hbox{ and } \ B:= k_{\textup{adv}}\lambda^{-1}.\]
    By \eref{tildeomega-limit-holes}, \lref{divergence-non-divergence-hole-info}, and $\mathcal{Z}$-periodicity  there is some
$c_*\in [\gamma_dA,\gamma_dB]$ such that $\omega$ solves \eref{CM-cell-ptwise} and \eref{CM-cell-1} and then \tref{membrane-cell-ptwise} implies that
\[c_* = |\Box_{\mathcal{Z}}| = |\xi|.\] 
Therefore
\begin{equation}\label{e.xi-upper-bound}
|\xi| \leq \gamma_d B = \gamma_d k_{\textup{adv}}\lambda^{-1}
\end{equation}
yielding \eqref{e.upper_2}. We have now proved \pref{asymptotic-exp-upper-bound}.

Note that, in particular, \eqref{e.xi-upper-bound} implies that $\lambda$ is finite, which allows us to consider the original error function $\omega_{\delta}$. Since the solution of \eref{CM-cell-ptwise} is unique modulo constants, it follows that the full sequence $\tilde{\omega}_{\delta_j}$ converges locally uniformly in $\overline{\R^d_+} \setminus \mathcal{Z}$ to $\omega$ (modulo constants). Thus 
\begin{equation}\label{e.convergence}\omega_{\delta_j} = \frac{\mu_{\delta_j}}{\delta_j^{d-1}} \tilde{\omega}_{\delta_j} \to \lambda \omega \ \hbox{ locally uniformly (modulo constants) in } \ \R^d_+ \setminus \mathcal{Z}.
\end{equation}

 We now take as assumption that $0 <k_{\textup{adv}} < +\infty$. Let us choose $u_{\delta}$ to be the minimal supersolutions with $\mu_\delta = Q^\delta_{\textup{adv}}-1$. Then, combining \eqref{e.upper_2} with \pref{asymptotic-exp-lower-bound} for the limit infimum, we conclude that
\[\lim_{\delta \to 0} \frac{\mu_\delta}{\delta^{d-1}} = \gamma_d|\xi|^{-1}k_{\textup{adv}}.\]
  \eqref{e.convergence} now yields that
\[\omega_\delta \to \gamma_d|\xi|^{-1}k_{\textup{adv}} \omega \ \hbox{ locally uniformly (modulo constants) in } \ \R^d_+ \setminus \mathcal{Z},\]
which concludes the proof of Theorem~\ref{t.udelta-expansion}.
\end{proof}

\begin{remark}
    The above proof is given for general $\mu_{\delta}$ to also show the convergence of $\omega_\delta$ for other sequences of strong Birkhoff plane-like solutions (not just with extremal slopes), which exist for every slope $\alpha_\delta$ in the pinning interval $[Q^{\delta}_{\textup{rec}},Q^\delta_{\textup{adv}}]$ by \tref{cell-solutions-reference}, as long as $\frac{\alpha_\delta-1}{\delta^{d-1}}$ converges to a nonzero value. 
\end{remark}
Finally we return to prove \lref{omega-neumann}.
\begin{proof}[Proof of \lref{omega-neumann}]
First we show compactness and the asymptotic bounds near the holes. By hypothesis $\lambda_\delta := \frac{\mu_\delta}{\delta^{d-1}}\to \lambda\neq 0$ (but possibly equal to $\pm\infty$), so  we may assume that $\mu_\delta$ has a fixed sign and $\lambda_\delta^{-1}$ are bounded. We will continue in the case $\mu_\delta>0$: the other case is similar. 

By multiplying $\lambda_\delta^{-1} = \frac{\delta^{d-1}}{\mu_\delta}>0$ on the bounds from \pref{2d-capacitory-bounds} (in $d=2$) or \pref{refined-flatness-3d} (in $d \geq 3$), we have
\[|\tilde{\omega}_\delta| \leq C(1+\lambda_\delta^{-1}) \ \hbox{ in }  \ \Omega_\delta \setminus \cup_{z \in \mathcal{Z}} B_{1/2}(z),\]
and the bounds near the defects, in $B_{1/2}(0) \cap \Omega_\delta$:
\begin{equation}\label{e.tildeomega-bound-holes}
    \begin{split}
    &\tilde{\omega}_{\delta}(x)\leq \lambda_\delta^{-1}k_{\textup{adv}}^\delta \Phi(x)+C+C\lambda_\delta^{-1}(1+\delta|x|^{1-d}),\\
    &\tilde{\omega}_{\delta}(x) \geq \lambda_\delta^{-1} k_{\textup{rec}}^\delta \Phi(x) -C-C\lambda_\delta^{-1}(1+\delta|x|^{1-d}),
\end{split}
\end{equation}
where $\lim_{\delta \to 0} k_{\textup{adv}/\textup{rec}}^\delta = k_{\textup{adv}/\textup{rec}}$. By $\mathcal{Z}$-periodicity  translated versions of the bounds hold in $B_{1/2}(z)$ for $z \in \mathcal{Z}$.   

  Let $K$ be a compact subset of $\R^d \setminus \mathcal{Z}$. The uniform bounds above imply that the sequence $\tilde{\omega}_\delta$ is uniformly bounded on $K \cap \overline{\Omega}_\delta$, and \cite{FKexterior}*{Lemma 2.3} implies that $\tilde{\omega}_\delta$ are uniformly $C^{1,\alpha}$ on $K \cap \overline{\Omega}_\delta$. 
Note that, by Lemma~\ref{l.initial-flatness},
\begin{equation}\label{e.Omega-delta-convergence-est}
    \{x_d + s_\delta^- > 0\} \subset \Omega_\delta \subset \{x_d +s_\delta^+> 0\} \ \hbox{ with } \ s_\delta^\pm \to 0
\end{equation}
so $\Omega_\delta$ is converging to $\R^d_+$ as $\delta\to 0$.  Putting these together,  \eref{tildeomega-limit-holes} follows now easily from \eref{tildeomega-bound-holes} by first fixing $z\in \mathcal{Z}$, taking the limit $\delta \to 0$ and then dividing by $\Phi(x-z)$ and taking the limit $x \to z$.

Next we establish the PDE conditions for $\omega$. Combined with $\omega_\delta$ being harmonic in $\Omega_\delta$ and the convergence of $\Omega_{\delta}$ from \eref{Omega-delta-convergence-est}, we conclude that $\omega$ is harmonic in $\R^d_+$. 

It remains to establish the Neumann condition away from the holes. The proof follows arguments from De Silva \cite{DeSilva}*{Lemma 4.1}.  Consider a smooth test function $\varphi$ touching $\omega$ from above at some $x_0 \in \partial \R^d_+ \setminus \mathcal{Z}$ with $\Delta \varphi(x_0)<0$. Then there is $\overline{\Omega_\delta}  \ni x_\delta \to x_0$ and $c_\delta \to 0$ so that $\varphi+c_{\delta}$ touches $\tilde{\omega}_\delta$ from above at $x_\delta$. Since $\tilde{\omega}_\delta$ is harmonic in $\Omega_\delta$, it must be that $x_\delta \in \partial \Omega_\delta$. Recalling the definition of $\tilde{\omega}_\delta$,
    \[\varphi(x)+c_\delta - \frac{1}{\mu_\delta}((1+\mu_\delta) (x_d+s_{\delta})_+-u_\delta(x)) \ \hbox{ has a local maximum $0$ at $x_\delta$,}\]
    and so if we define
    \[\tilde{\varphi}(x) := (1+\mu_\delta) (x_d+s_\delta)_+-\mu_\delta(\varphi(x)+c_\delta),\]
     then $u_\delta - \tilde{\varphi}$ has a local minimum $0$ at $x_\delta$. Since $x_0 \not\in \mathcal{Z}$, $x_\delta \not\in \cup_{z \in \mathcal{Z}}B_\delta(z)$ for $\delta\ll 1$. Since $Q=1$ in a local neighborhood of $x_{\delta}$, from the viscosity supersolution property of $u_{\delta}$ it follows that 
\begin{align*}
1 &\geq |\grad \tilde\varphi(x_\delta)|^2 = |(1+\mu_\delta)e_d-\mu_\delta\grad \varphi(x_\delta)|^2\\
&=(1+\mu_\delta)^2-2\mu_\delta(1+\mu_\delta) \partial_{d}\varphi(x_\delta)+\mu_\delta^2|\grad \varphi(x_\delta)|^2\\
&=1+2\mu_\delta(1-\partial_{d}\varphi(x_\delta)) + \mu_\delta^2|e_d-\grad \varphi(x_\delta)|^2.
\end{align*}
Rearranging this gives
\[\partial_{d}\varphi(x_\delta) \geq 1 + \frac{1}{2}\mu_\delta|e_d-\grad \varphi(x_\delta)|^2\]
and taking the limit as $\delta \to 0$
\[\partial_d \varphi(x_0) \geq 1.\]
Thus $\omega$ satisfies viscosity subsolution property for the Neumann condition $\partial_d \omega = 1$ on $\partial \R^d_+ \setminus \mathcal{Z}$. The viscosity supersolution property is proved similarly.

Lastly, note that  $(1+\mu_\delta)x_d - u_\delta(x)$ is a bounded harmonic function in $\Omega_\delta$, which is periodic with respect to $\mathcal{Z}$ translations. Thus, by \lref{periodic-bdry-layer}, $\langle \partial_d u_\delta(\cdot ,y_d) \rangle' = 1+\mu_\delta$ for any $y_d>0$ large enough that $\{x \cdot e_d \geq y_d\} \subset \Omega_\delta$, where $\langle \cdot \rangle'$ is the tangential period cell average.  Due to the Hausdorff convergence of $\overline{\Omega_\delta} \to \overline{\R^d_+}$, this implies that $\langle \partial_d\omega(\cdot,y_d) \rangle' = 0$ for all $y_d>0$.
\end{proof}

  \bibliographystyle{amsplain}
\bibliography{single-site-articles}

\end{document}

%% file: figures/ellipses-pic.tikz
\begin{tikzpicture}
\draw (2,1) ellipse (.5 and 1) node {$+$};
\draw (-2,-1) ellipse (.5 and 1) node {$-$};
\draw (2,-1) ellipse (1 and .5) node {$-$};
\draw (-2,1) ellipse (1 and .5) node {$+$};
\end{tikzpicture}

%% file: single-site-dilute-defects.bbl
\begin{bibdiv}
\begin{biblist}

\bib{AbedinFeldman}{article}{
      author={Abedin, Farhan},
      author={Feldman, William~M.},
       title={Quantitative convergence of the ``bulk'' free boundary in an
  oscillatory obstacle problem},
        date={2024},
        ISSN={1463-9963,1463-9971},
     journal={Interfaces Free Bound.},
      volume={26},
      number={1},
       pages={31\ndash 44},
         url={https://doi.org/10.4171/ifb/501},
      review={\MR{4705640}},
}

\bib{AlbertiDeSimone}{article}{
      author={Alberti, Giovanni},
      author={DeSimone, Antonio},
       title={Quasistatic evolution of sessile drops and contact angle
  hysteresis},
        date={2011},
        ISSN={0003-9527,1432-0673},
     journal={Arch. Ration. Mech. Anal.},
      volume={202},
      number={1},
       pages={295\ndash 348},
         url={https://doi.org/10.1007/s00205-011-0427-x},
      review={\MR{2835870}},
}

\bib{AubinThierry1998Snpi}{book}{
      author={Aubin, Thierry},
       title={Some nonlinear problems in riemannian geometry},
    language={eng},
      series={Springer monographs in mathematics},
        date={1998},
        ISBN={3540607528},
}

\bib{ConstantinLeDoussalWiese}{article}{
      author={Bachas, Constantin},
      author={Le~Doussal, Pierre},
      author={Wiese, Kay~J\"org},
       title={Wetting and minimal surfaces},
        date={2007Mar},
     journal={Phys. Rev. E},
      volume={75},
       pages={031601},
         url={https://link.aps.org/doi/10.1103/PhysRevE.75.031601},
}

\bib{Bico_2001}{article}{
      author={Bico, J},
      author={Tordeux, C},
      author={Quéré, D},
       title={Rough wetting},
        date={2001-07},
        ISSN={1286-4854},
     journal={Europhysics Letters (EPL)},
      volume={55},
      number={2},
       pages={214–220},
         url={http://dx.doi.org/10.1209/epl/i2001-00402-x},
}

\bib{BraidesScilla}{article}{
      author={Braides, Andrea},
      author={Scilla, Giovanni},
       title={Motion of discrete interfaces in periodic media},
        date={2013},
        ISSN={1463-9963,1463-9971},
     journal={Interfaces Free Bound.},
      volume={15},
      number={4},
       pages={451\ndash 476},
         url={https://doi.org/10.4171/IFB/310},
      review={\MR{3148630}},
}

\bib{CaffarelliMellet-StickingDrop}{article}{
      author={Caffarelli, L.~A.},
      author={Mellet, A.},
       title={Capillary drops: contact angle hysteresis and sticking drops},
        date={2007},
        ISSN={0944-2669,1432-0835},
     journal={Calc. Var. Partial Differential Equations},
      volume={29},
      number={2},
       pages={141\ndash 160},
         url={https://doi.org/10.1007/s00526-006-0036-y},
      review={\MR{2307770}},
}

\bib{CaffarelliMellet-RandomObstacle}{article}{
      author={Caffarelli, L.~A.},
      author={Mellet, A.},
       title={Random homogenization of an obstacle problem},
        date={2009},
        ISSN={0294-1449,1873-1430},
     journal={Ann. Inst. H. Poincar\'e{} C Anal. Non Lin\'eaire},
      volume={26},
      number={2},
       pages={375\ndash 395},
         url={https://doi.org/10.1016/j.anihpc.2007.09.001},
      review={\MR{2504035}},
}

\bib{caffarelli-mellet}{article}{
      author={Caffarelli, Luis~A},
      author={Mellet, Antoine},
       title={Capillary drops on an inhomogeneous surface},
        date={2007},
     journal={Contemporary Mathematics},
      volume={446},
       pages={175\ndash 202},
}

\bib{CS}{book}{
      author={Caffarelli, Luis~A},
      author={Salsa, Sandro},
       title={A geometric approach to free boundary problems},
   publisher={American Mathematical Soc.},
        date={2005},
      volume={68},
}

\bib{Chodosh-Li}{article}{
      author={Chodosh, Otis},
      author={Edelen, Nick},
      author={Li, Chao},
       title={Improved regularity for minimizing capillary hypersurfaces},
        date={2024},
     journal={arXiv preprint arXiv:2401.08028},
}

\bib{Chu_2010}{article}{
      author={Chu, Kuang-Han},
      author={Xiao, Rong},
      author={Wang, Evelyn~N.},
       title={Uni-directional liquid spreading on asymmetric nanostructured
  surfaces},
        date={2010-03},
        ISSN={1476-4660},
     journal={Nature Materials},
      volume={9},
      number={5},
       pages={413–417},
         url={http://dx.doi.org/10.1038/nmat2726},
}

\bib{CioranescuUnfolding}{article}{
      author={Cioranescu, D.},
      author={Damlamian, A.},
      author={Donato, P.},
      author={Griso, G.},
      author={Zaki, R.},
       title={The periodic unfolding method in domains with holes},
        date={2012},
     journal={SIAM Journal on Mathematical Analysis},
      volume={44},
      number={2},
       pages={718\ndash 760},
}

\bib{Cioranescu-Donato}{article}{
      author={Cioranescu, Doina},
      author={Donato, Patrizia},
       title={Homogénéisation du problème de neumann non homogène dans des
  ouverts perfores},
        date={1988},
     journal={Asymptotic Analysis},
      volume={1},
      number={2},
       pages={115\ndash 138},
}

\bib{cioranescu1997strange}{incollection}{
      author={Cioranescu, Doina},
      author={Murat, Fran{\c{c}}ois},
       title={A strange term coming from nowhere},
        date={1997},
   booktitle={Topics in the mathematical modelling of composite materials},
   publisher={Springer},
       pages={45\ndash 93},
}

\bib{Courbin_2007}{article}{
      author={Courbin, Laurent},
      author={Denieul, Etienne},
      author={Dressaire, Emilie},
      author={Roper, Marcus},
      author={Ajdari, Armand},
      author={Stone, Howard~A.},
       title={Imbibition by polygonal spreading on microdecorated surfaces},
        date={2007-08},
        ISSN={1476-4660},
     journal={Nature Materials},
      volume={6},
      number={9},
       pages={661–664},
         url={http://dx.doi.org/10.1038/nmat1978},
}

\bib{CourteBattacharyaDondl}{article}{
      author={Courte, Luca},
      author={Bhattacharya, Kaushik},
      author={Dondl, Patrick},
       title={Bounds on precipitate hardening of line and surface defects in
  solids},
        date={2020},
        ISSN={0044-2275,1420-9039},
     journal={Z. Angew. Math. Phys.},
      volume={71},
      number={3},
       pages={Paper No. 99, 14},
         url={https://doi.org/10.1007/s00033-020-01327-3},
      review={\MR{4108717}},
}

\bib{CourteDondlStefanelli}{article}{
      author={Courte, Luca},
      author={Dondl, Patrick},
      author={Stefanelli, Ulisse},
       title={Pinning of interfaces by localized dry friction},
        date={2020},
        ISSN={0022-0396},
     journal={Journal of Differential Equations},
      volume={269},
      number={9},
       pages={7356\ndash 7381},
  url={https://www.sciencedirect.com/science/article/pii/S0022039620303193},
}

\bib{DeRosa25}{article}{
      author={De~Rosa, Antonio},
      author={Neumayer, Robin},
      author={Resende, Reinaldo},
       title={Rigidity of critical points of hydrophobic capillary
  functionals},
        date={2025},
     journal={arXiv preprint arXiv:2509.22532},
}

\bib{DeSimoneGrunewaldOtto}{article}{
      author={DeSimone, Antonio},
      author={Grunewald, Natalie},
      author={Otto, Felix},
       title={A new model for contact angle hysteresis},
        date={2007},
        ISSN={1556-1801},
     journal={Netw. Heterog. Media},
      volume={2},
      number={2},
       pages={211\ndash 225},
         url={https://doi.org/10.3934/nhm.2007.2.211},
}

\bib{DirrYip}{article}{
      author={Dirr, N.},
      author={Yip, N.~K.},
       title={Pinning and de-pinning phenomena in front propagation in
  heterogeneous media},
        date={2006},
        ISSN={1463-9963,1463-9971},
     journal={Interfaces Free Bound.},
      volume={8},
      number={1},
       pages={79\ndash 109},
         url={https://doi.org/10.4171/IFB/136},
      review={\MR{2231253}},
}

\bib{eral2013contact}{article}{
      author={Eral, Husseyin~B},
      author={Mannetje, DJCM},
      author={Oh, Jung~Min},
       title={Contact angle hysteresis: a review of fundamentals and
  applications},
        date={2013},
     journal={Colloid and polymer science},
      volume={291},
      number={2},
       pages={247\ndash 260},
}

\bib{FeldmanPozar2024}{article}{
      author={Feldman, W.~M.},
      author={Po\v{z}\'ar, N.},
       title={On convex comparison for exterior bernoulli problems with
  discontinuous anisotropy},
        date={2024},
     journal={Interfaces and Free Boundaries},
      volume={26},
      number={4},
       pages={459\ndash 468},
         url={https://doi.org/10.4171/IFB/525},
}

\bib{feldman2021limit}{article}{
      author={Feldman, William~M},
       title={Limit shapes of local minimizers for the {A}lt--{C}affarelli
  energy functional in inhomogeneous media},
        date={2021},
     journal={Archive for Rational Mechanics and Analysis},
      volume={240},
       pages={1255\ndash 1322},
}

\bib{FK}{article}{
      author={Feldman, William~M},
      author={Kim, Inwon~C},
       title={Dynamic stability of equilibrium capillary drops},
        date={2014},
     journal={Archive for Rational Mechanics and Analysis},
      volume={211},
      number={3},
       pages={819\ndash 878},
}

\bib{FKoriginal}{article}{
      author={Feldman, William~M},
      author={Kim, Inwon~C},
       title={The pinning effect of dilute defects},
        date={2025},
      eprint={2512.11152},
         url={https://arxiv.org/abs/2512.11152},
}

\bib{FKexterior}{article}{
      author={Feldman, William~M},
      author={Kim, Inwon~C},
       title={Solutions of the bernoulli one-phase problem with a defect},
        date={2026},
     journal={arXiv preprint},
}

\bib{feldman2024obstacleapproachrateindependent}{misc}{
      author={Feldman, William~M},
      author={Kim, Inwon~C},
      author={Požár, Norbert},
       title={An obstacle approach to rate independent droplet evolution},
        date={2024},
         url={https://arxiv.org/abs/2410.06931},
}

\bib{feldman-kim-pozar}{article}{
      author={Feldman, William~M.},
      author={Kim, Inwon~C.},
      author={Požár, Norbert},
       title={An obstacle approach to rate-independent droplet evolution},
        date={2026},
     journal={Forum of Mathematics, Sigma},
      volume={14},
       pages={e38},
}

\bib{feldman-smart}{article}{
      author={Feldman, William~M},
      author={Smart, Charles~K},
       title={A free boundary problem with facets},
        date={2019},
     journal={Archive for Rational Mechanics and Analysis},
      volume={232},
      number={1},
       pages={389\ndash 435},
}

\bib{FeldmanZhang2019}{article}{
      author={Feldman, William~M.},
      author={Zhang, Yuming~Paul},
       title={Continuity properties for divergence form boundary-data
  homogenization problems},
        date={2019},
     journal={Analysis \& PDE},
      volume={12},
      number={8},
       pages={1963\ndash 2002},
         url={https://doi.org/10.2140/apde.2019.12.1963},
}

\bib{guillen-kim}{article}{
      author={Guillen, Nestor},
      author={Kim, Inwon},
       title={Quasistatic droplets in randomly perforated domains},
        date={2015},
     journal={Archive for Rational Mechanics and Analysis},
      volume={215},
      number={1},
       pages={211\ndash 281},
}

\bib{heida2020}{article}{
      author={Heida, Martin},
       title={Stochastic homogenization on randomly perforated domains},
        date={2020},
     journal={arXiv preprint arXiv:2001.10373},
}

\bib{JerisonKamburov2}{article}{
      author={Jerison, David},
      author={Kamburov, Nikola},
       title={Free boundaries subject to topological constraints},
        date={2019},
        ISSN={1078-0947,1553-5231},
     journal={Discrete Contin. Dyn. Syst.},
      volume={39},
      number={12},
       pages={7213\ndash 7248},
         url={https://doi.org/10.3934/dcds.2019301},
      review={\MR{4026187}},
}

\bib{Joanny-deGennes}{article}{
      author={Joanny, J.~F.},
      author={de~Gennes, P.~G.},
       title={A model for contact angle hysteresis},
        date={198407},
        ISSN={0021-9606},
     journal={The Journal of Chemical Physics},
      volume={81},
      number={1},
       pages={552\ndash 562},
  eprint={https://pubs.aip.org/aip/jcp/article-pdf/81/1/552/18948471/552_1_online.pdf},
         url={https://doi.org/10.1063/1.447337},
}

\bib{joanny1990motion}{article}{
      author={Joanny, JF},
      author={Robbins, Mark~O},
       title={Motion of a contact line on a heterogeneous surface},
        date={1990},
     journal={The Journal of chemical physics},
      volume={92},
      number={5},
       pages={3206\ndash 3212},
}

\bib{KamburovWang}{article}{
      author={Kamburov, Nikola},
      author={Wang, Kelei},
       title={Nondegeneracy for stable solutions to the one-phase free boundary
  problem},
        date={2024},
     journal={Math. Ann.},
      volume={388},
      number={3},
       pages={2705\ndash 2726},
}

\bib{kim-contact}{article}{
      author={Kim, Inwon~C},
       title={Homogenization of a model problem on contact angle dynamics},
        date={2008},
     journal={Communications in Partial Differential Equations},
      volume={33},
      number={7},
       pages={1235\ndash 1271},
}

\bib{kinderlehrer1977regularity}{article}{
      author={Kinderlehrer, David},
      author={Nirenberg, Louis},
       title={Regularity in free boundary problems},
        date={1977},
     journal={Annali della Scuola Normale Superiore di Pisa-Classe di Scienze},
      volume={4},
      number={2},
       pages={373\ndash 391},
}

\bib{KriventsovWeiss}{article}{
      author={Kriventsov, Dennis},
      author={Weiss, Georg~S.},
       title={Rectifiability, finite {H}ausdorff measure, and compactness for
  non-minimizing {B}ernoulli free boundaries},
        date={2025},
        ISSN={0010-3640,1097-0312},
     journal={Comm. Pure Appl. Math.},
      volume={78},
      number={3},
       pages={545\ndash 591},
         url={https://doi.org/10.1002/cpa.22226},
      review={\MR{4850026}},
}

\bib{Lieberman}{article}{
      author={Lieberman, Gary~M.},
       title={Local estimates for subsolutions and supersolutions of oblique
  derivative problems for general second order elliptic equations},
        date={1987},
        ISSN={0002-9947,1088-6850},
     journal={Trans. Amer. Math. Soc.},
      volume={304},
      number={1},
       pages={343\ndash 353},
         url={https://doi.org/10.2307/2000717},
      review={\MR{906819}},
}

\bib{PapanicolaouVaradhan1980}{inproceedings}{
      author={Papanicolaou, G.},
      author={Varadhan, S.~R.~S.},
       title={Diffusion in regions with many small holes},
        date={1980},
   booktitle={Proceedings of the vilnius conference on probability theory and
  mathematical statistics},
      editor={Grigelionis, B.},
      series={Lecture Notes in Control and Information Sciences},
      volume={25},
   publisher={Springer},
       pages={190\ndash 206},
}

\bib{quere2008wetting}{article}{
      author={Qu{\'e}r{\'e}, David},
       title={Wetting and roughness},
        date={2008},
     journal={Annu. Rev. Mater. Res.},
      volume={38},
       pages={71\ndash 99},
}

\bib{Raj_2014}{article}{
      author={Raj, Rishi},
      author={Adera, Solomon},
      author={Enright, Ryan},
      author={Wang, Evelyn~N.},
       title={High-resolution liquid patterns via three-dimensional droplet
  shape control},
        date={2014-09},
        ISSN={2041-1723},
     journal={Nature Communications},
      volume={5},
      number={1},
         url={http://dx.doi.org/10.1038/ncomms5975},
}

\bib{raphael1989dynamics}{article}{
      author={Raphael, Elie},
      author={de~Gennes, Pierre-Gilles},
       title={Dynamics of wetting with nonideal surfaces. the single defect
  problem},
        date={1989},
     journal={The Journal of chemical physics},
      volume={90},
      number={12},
       pages={7577\ndash 7584},
}

\bib{reyssat2009contact}{article}{
      author={Reyssat, Mathilde},
      author={Qu{\'e}r{\'e}, David},
       title={Contact angle hysteresis generated by strong dilute defects},
        date={2009},
     journal={The Journal of Physical Chemistry B},
      volume={113},
      number={12},
       pages={3906\ndash 3909},
}

\bib{schwartz1985contact}{article}{
      author={Schwartz, Leonard~W},
      author={Garoff, Stephen},
       title={Contact angle hysteresis on heterogeneous surfaces},
        date={1985},
     journal={Langmuir},
      volume={1},
      number={2},
       pages={219\ndash 230},
}

\bib{DeSilva}{article}{
      author={Silva, Daniela~De},
       title={Free boundary regularity for a problem with right hand side},
        date={2011},
     journal={Interfaces and Free Boundaries},
      volume={13},
      number={2},
       pages={223\ndash 238},
}

\bib{Weiss}{article}{
      author={Weiss, G.~S.},
       title={A singular limit arising in combustion theory: Fine properties of
  the free boundary},
        date={2003-07},
        ISSN={0944-2669},
     journal={Calculus of Variations and Partial Differential Equations},
      volume={17},
      number={3},
       pages={311\ndash 340},
         url={https://doi.org/10.1007/s00526-002-0171-z},
      review={\MR{1989835}},
}

\end{biblist}
\end{bibdiv}
